\documentclass[a4paper, 12pt, oneside, article]{memoir}
\usepackage{wandset}
\bibliography{Set.bib}
\title{Wand/Set Theories}
\subtitle{A realization of {Conway}'s mathematicians' liberation movement, with an application to {Church}'s set theory with a universal set}
\author{Tim Button}
\emailaddress{tim.button@ucl.ac.uk}
\begin{document}\midsloppy
\pagestyle{nicelypage}
\maketitletop\selectbodyfont
\begin{abstract}
	Consider a variant of the usual story about the iterative conception of sets. As usual, at every stage, you find all the (bland) sets of objects which you found earlier. But you also find the result of tapping any earlier-found object with any magic wand (from a given stock of magic wands).
	
	By varying the number and behaviour of the wands, we can flesh out this idea in many different ways. This paper's main Theorem is that any loosely constructive way of fleshing out this idea is synonymous with a \ZF-like theory.
	
	This Theorem has rich applications; it realizes John \citepossess{Conway:ONG} Mathematicians' Liberation Movement; and it connects with a lovely idea due to Alonzo \textcite{Church:STUS}.
	
	\textcolor{myblue}{This paper is forthcoming in \emph{Journal of Symbolic Logic}. Some of the proofs in the published version of this paper are quite compressed; in those cases, this document offers a link (in the margin) to further details (in appendices \ref{proofsforwev}--\ref{s:app:CUSstuff}).}
\end{abstract}
\begin{center}***\end{center}
\
\\Here is a template for introducing mathematical objects. 
\begin{storytime}\textbf{The Wand/Set Template.}
	Objects are found in stages. For every stage \stage{s}:
	\begin{listbullet}
		\item for any things found before \stage{s}, you find at \stage{s} the bland set whose members are exactly those things;
		\item for anything, $a$, which was found before \stage{s}, and for any magic wand, $w$, you find at \stage{s} the result of tapping $a$ with $w$ (if tapping $a$ with $w$ yields something other than a bland set).
	\end{listbullet}
	You find nothing else at \stage{s}.
\end{storytime}\noindent
As I will explain in \S\ref{s:Motivations}: this Template has rich applications; it realizes John \citepossess{Conway:ONG} Mathematicians' Liberation Movement; and it generalizes a lovely idea due to Alonzo \textcite{Church:STUS}. 

Some parts of this Template are familiar: the bit about bland sets is just the ordinary story we tell about the cumulative iterative conception of set. But this talk of ``magic wands'' is new and different. Moreover, parts of the story are left unspecified. We know very little about: the sequence of stages; what the wands are; the circumstances under which tapping something with a wand will yield anything; nor what will thereby be revealed. 

This under-specification is deliberate: the Template is left \emph{as} a template precisely so that we can flesh it out in umpteen different ways. Nevertheless, the Main Theorem of the paper is that \emph{any loosely constructive way of fleshing out the Template is synonymous with a \ZF-like theory}. (I define ``loosely constructive'' in \S\ref{s:Constructive}, and ``\ZF-like'' in \S\ref{s:Strategy}.)

\section{Motivations}\label{s:Motivations}
I will start by explaining why the Wand/Set Template is worthy of attention. Briefly, this is because it provides us with an extremely rich framework for introducing mathematical objects (see \S\ref{s:Motivations:Examples}). I hope this is enough to make it interesting. But I am not above appealing to authorities; in this case, Conway (\S\ref{s:Motivations:Conway}) and Church (\S\ref{s:Motivations:Church}).

\subsection{Simple examples}\label{s:Motivations:Examples}
The Wand/Set Template is very abstract. So it will help to start with some concrete instances of the Template. This will also illustrate the scope of my Main Theorem.
\begin{example}[Conway-games]\label{ex:ConwayGames}
	\textcite{Conway:ONG} famously provided a general theory of (hereditary) two-player games. A Conway-game has \emph{left-options} and \emph{right-options}, and Conway-games are identical iff they have the same left- and right-options. 
	
	Here is a sketch of how to provide a wand/set theory of Conway-games. 
	We have a single wand, \textsf{game}. Tapping an ordered-pair, $\tuple{a,b}$, of bland sets with  \textsf{game} (typically) yields a new object; we deem this to be a Conway-game whose left-options are exactly $a$'s members and whose right-options are exactly $b$'s members. By courtesy, a bland set, $a$, is the Conway-game with no right-options and whose left-options are exactly $a$'s members.
	
	Here is the idea in a bit more detail. Reserve curly-brackets for bland sets, and $\in$ for membership of the same. With $\tuple{a,b} \coloneq \{\{a\}, \{a,b\}\}$, stipulate that \textsf{game} only acts on objects of the form $\tuple{a,b}$, where both $a$ and $b$ are bland sets. Then $\textsf{game}(\tuple{a,\emptyset}) \coloneq a$, and if $b \neq \emptyset$ then $\textsf{game}(\tuple{a,b})$ is an object which is not a bland set; and games are governed by this rule:
		$$\textsf{game}(\tuple{a,b}) = \textsf{game}(\tuple{a',b'}) \lonlyif (a = a' \land b = b')$$
	Reading ``$x \circ_\textsc{l} y$'' as ``$x$ is a left-option of $y$'', and similarly for ``$x \circ_\textsc{r} y$'', explicitly define:
	\begin{align*}
		x \circ_\textsc{l} y &\colonequiv \exists a \exists b(y = \textsf{game}(\tuple{a,b}) \land x \in a)\\
		x \circ_\textsc{r} y &\colonequiv \exists a \exists b(y = \textsf{game}(\tuple{a,b}) \land x \in b)
	\end{align*}	
	Properly implemented, we have a theory whose objects are exactly the Conway-games. By my Main Theorem, this is synonymous with a \ZF-like theory.\footnote{\label{fn:CoxKaye}This essentially subsumes \citepossess{CoxKaye:AZF} result concerning \textsf{Amphi-ZF}. Note a difference: \citeauthor{CoxKaye:AZF} take $\circ_\textsc{l}$ and $\circ_\textsc{r}$ as primitives; I take $\blandpred$ as a primitive and (arbitrarily) treat bland sets as Conway-games with no right-options. \citeauthor{CoxKaye:AZF}'s approach is more natural; mine lets me press the theory of Conway-games into the mould of a wand/set theory (see Definition \ref{def:WandSet}) without positing bland sets as objects which are \emph{not} themselves Conway-games (see the end of \S\ref{s:Motivations:Conway}).}
\end{example}
\begin{example}[Partial functions]\label{ex:PartialFunctions}
	Consider a cumulative iterative hierarchy of one-place partial functions: we start with the function which is undefined for every input; at every layer,\footnote{Such \emph{layers} are not precisely the \emph{stages} mentioned in the Wand/Set Template. As we will see: $\emptyset$ is (by courtesy) the function which is undefined for any input; $\{\emptyset\}$ is (by courtesy) the identity function which just maps $\emptyset\mapsto\emptyset$; now let $f$ be the function which maps $\{\emptyset\} \mapsto \emptyset$ and is undefined for all other inputs. So $f$ is found at the \emph{third layer} in the explanation at the start of this example; but $f = \textsf{fun}(\{\tuple{\{\emptyset\}, \emptyset}\})$ is found at the \emph{sixth stage} of the Wand/Set implementation.
		
	This mismatch between layers and stages is ultimately unimportant. Provided that there is no last stage or layer (see \ref{WS:weakheight}), all the same objects will be found, whether we tell the story in terms of layers or stages; moreover, we can easily define the layer-ordering by a recursion on the stage-ordering, and vice versa (and any wand/set theory, in the sense of Definition \ref{def:WandSet}, has a suitable recursion theorem.)} 
 	we find all possible one-place partial functions whose domain-and-range are exhausted by objects we found earlier.\footnote{This idea was outlined by \textcite{Robinson:TCMVNS}.} 
	
	Here is a sketch of how to put this into the form of a wand/set theory. We have a single wand, \textsf{fun}. We use bland sets to encode functions in the usual way; tapping such a bland set with \textsf{fun} (typically) yields a new object; we deem this to be a function of the appropriate sort. By courtesy, a bland set, $a$, will be the identity function defined only over $a$'s members. Of course, more detail can be offered (as in Example \ref{ex:ConwayGames}); but this ultimately yields a theory whose objects are exactly the (hereditary) partial one-place functions. By the Main Theorem, this is synonymous with a \ZF-like theory.\footnote{This essentially subsumes \citepossess{Meadows:WSTCNBA} and my \parencite*{Button:ICFICS} independently-obtained synonymy results about iterative theories of functions.}
	
	We can easily overcome the restriction to \emph{one}-place partial functions by positing a countable infinity of wands. For each $n$, tapping appropriate bland sets with $\textsf{fun}_n$ yields the $n$-place partial functions whose domain-and-range are exhausted by what we have already found. Again: the result is synonymous with a \ZF-like theory. 
\end{example}
\begin{example}[Multisets]\label{ex:Multisets}
	Consider a cumulative iterative hierarchy of multisets: we start with the empty multiset; at every layer, we find all possible multisets whose members and associated multiplicities (when the multiplicity is $> 1$) were found earlier. So for any $x$: we first find multisets with just one copy of $x$ immediately after we have found $x$; and for each cardinal $\mathfrak{a} > 1$, we first find multisets with $\mathfrak{a}$ copies of $x$ immediately after we have found both $x$ and $\mathfrak{a}$. 
	
	To render this as a wand/set theory, just treat a multiset, $s$, as a one-place partial function, where $s(x) = \mathfrak{a}$ indicates that $s$ has $\mathfrak{a}$ copies of $x$. Then the theory of multisets is just a slight variant of Example \ref{ex:PartialFunctions}. The only differences are: (i) our functions' values must always be cardinals (in the bland sense), and (ii) each bland set, $a$, is treated as the multiset which contains exactly $1$ copy of each of $a$'s members. Again: the result is synonymous with a \ZF-like theory.\footnote{So this goes well beyond e.g.\ \citepossess{Blizard:MT} embedding result.}
\end{example}
\begin{example}[Accessible pointed graphs]\label{ex:Aczel}
	\textcite{Aczel:NWS} thinks of ``sets'' in terms of accessible pointed graphs (APGs), whose ``members'' are equivalents of those objects which the APG's point can ``see''. Skimming over details, this is easy to render as a wand/set theory. We can encode the definition of an APG in terms of bland sets. Then we have a single wand, \textsf{set}, which behaves thus: when $a$ is a bland set which encodes an APG, $\textsf{set}(a)$ is an object whose ``APG-members'' are those things, $x$, such that $a$ encodes an edge from $a$'s point to an entity equivalent to $x$. 
	
	This example is somewhat different from Examples \ref{ex:ConwayGames}--\ref{ex:Multisets}. Before, I explained how to treat any bland set as an object of the relevant sort (a Conway-game; a partial function; a multiset). But I have given no uniform way to treat any bland set \emph{as} an APG, and indeed this is problematic. To see why, consider this APG: 
	\begin{tikzpicture}[baseline=-1mm]
		\node[inner sep=2.5pt] (A) at (0,0) [circle, fill=black] {};
		\draw[->] (A.30) to [out=30, in=330, looseness=10] (A.330);
	\end{tikzpicture}. 
	It depicts a (unique) ill-founded set, $\Omega$, whose sole ``member'' is $\Omega$ itself.\footnote{See \textcite[6--7]{Aczel:NWS}.} But the bland singleton of $\Omega$---which we should find at the next stage after we have found $\Omega$---has exactly one member, $\Omega$; so, by extensionality, that bland set should \emph{be} $\Omega$ itself; and so some bland set would be self-membered, contrary to what the Template tells us.\footnote{Specifically, it contradicts Lemma \ref{lem:WS:levof}\eqref{levnotin}.} 
	
	Consequently, \emph{this} approach does not give us a theory whose objects are exactly the APGs. Instead, it gives us a theory whose objects are exactly the APGs and bland sets thereof. Nonetheless, the result is synonymous with a \ZF-like theory.\footnote{Thanks to Luca Incurvati for getting me thinking about treating Aczel's APGs in terms of a wand/set theory (see also his \cite*[ch.7]{Incurvati:CS}).}
\end{example}

\subsection{Conway's liberationism}\label{s:Motivations:Conway}
Further instances of the Wand/Set Template could be offered, but I hope that Examples \ref{ex:ConwayGames}--\ref{ex:Aczel} already illustrate the Template's power. I now want to move beyond piecemeal examples and explore a more principled reason for investigating the Template. 

Immediately after outlining his theory of games (see Example \ref{ex:ConwayGames}), Conway advocated for a Mathematicians' Liberation Movement. His liberationism held that:
	\begin{listr-0}
		\item\label{Conway:constructive} ``Objects may be created from earlier objects in any reasonably constructive fashion.
		\item\label{Conway:equiv} Equality among the created objects can be any desired equivalence relation.''\footnote{\textcite[66]{Conway:ONG}.}
	\end{listr-0}
He also expected that there would be a meta-theorem, which would show that any theory whose objects are created in such a ``reasonably constructive'' fashion can be embedded within (some extension of) \ZF.

As the examples of \S\ref{s:Motivations:Examples} suggest, my Wand/Set Template realizes Conway's proposed liberationism. Regarding \eqref{Conway:constructive}: to implement the idea of creating some object ``from earlier objects'', all you need to do is form their bland set, and then tap that set with an appropriate wand (or perhaps with the right wands in the right order). Regarding \eqref{Conway:equiv}: when we specify the behaviour of our wands, we can allow ``any desired equivalence relation'' to dictate when the result of tapping $a$ with wand $w$ is equal to the result of tapping $b$ with wand $u$. And my Main Theorem exceeds Conway's expectations: all loosely constructive implementations of the Wand/Set Template are not merely \emph{embeddable} in (some extension of) \ZF, but \emph{synonymous} with a \ZF-like theory.

Admittedly, my Template may not fully exhaust the scope of Conway's liberationism. This is because Conway's \eqref{Conway:constructive} and \eqref{Conway:equiv} give us \emph{permissions}, not \emph{obligations}; they say what we \emph{may} create, without dictating what we \emph{must} create. My Template, by contrast, \emph{insists} on creating\footnote{I prefer to speak of ``finding'' (than ``creating'') objects; it seems marginally preferable to tell a story which is compatible with platonism.} bland sets at each stage. However, my insistence is unlikely to impinge seriously on Conway's proposed liberationism. For one thing: whenever Conway creates entities of some kind, $K$, he is confident that they will be embeddable into a \ZF-style hierarchy; so it is no great burden to request that, when he creates the $K$s, he also create the \ZF-like bland-sets into which his $K$s will be embedded. For another: if we \emph{really} want our theory to be a theory of $K$s only, we can often do this by creatively (and courteously) regarding bland sets as $K$s of a particular sort: we saw this in Examples \ref{ex:ConwayGames}--\ref{ex:Multisets}, though not in 
Example \ref{ex:Aczel}.

\subsection{Church's set theory with a universal set}\label{s:Motivations:Church}
I have presented some examples of the Wand/Set Template, and explained how the Template relates to Conway's Liberationism. I now want to consider a more surprising connection: we can regard Church's \cite*{Church:STUS} \emph{set theory with a universal set} as both an instance of the Wand/Set Template and as an exercise of Conway's liberationism.\footnote{\label{fn:ChurchCaveats}I have \textcite{Forster:ICS} to thank for approaching Church in essentially this way, and for making the link to Conway. I should note that, whilst I think it is \emph{cleanest} to understand Church in this way, it is not exactly his own presentation; for details, and some important caveats, see \S\ref{s:app:CUS:caveats}.
	
	This section builds on my \parencite*{Button:LT3}, which essentially covers the case without any $\textsf{cardinal}_n$ wands, so that my Main Theorem subsumes the synonymy result in my \parencite*{Button:LT3}.} 

Church's \cite*{Church:STUS} idea is roughly this. We have a \ZF-like hierarchy of bland sets. We also have some wands: a \textsf{complement} wand and, for each positive natural $n$, a \textsf{cardinal}$_n$ wand. These wands should behave as follows: 
\begin{listbullet}
	\item[(0)] Tapping any object, $a$, with \textsf{complement} yields $a$'s absolute complement; that is, the set whose members are exactly the things not in $a$. 
	\item[($n$)] Tapping some bland set, $a$, with \textsf{cardinal}$_n$, yields the set of all bland sets \emph{$n$-equivalent} with $a$ (if there are any; if there are none, this yields nothing; and tapping anything other than a bland set with $\textsf{cardinal}_n$ yields nothing).
\end{listbullet} 
This explanation mentions \emph{$n$-equivalence}. This is an explicitly defined equivalence relation (defined recursively on $n$; see Definition \ref{def:cardinaln}). Its details are largely irrelevant for now; what matters is just that Church's theory can be thought of as invoking a countable infinity of wands with explicitly defined actions.

What is particularly interesting about \emph{Church's} theory, however, is that---unlike Examples \ref{ex:ConwayGames}--\ref{ex:Aczel}---it does not obviously seem like an instance of \emph{Conway's} liberationism. Conway's \eqref{Conway:equiv} says that ``Equality among the created objects can be any desired equivalence relation''; this is fine, since Church's wands are based around explicitly defined equivalence relations. But Conway's \eqref{Conway:constructive} says that ``Objects may be created from earlier objects in any reasonably constructive fashion''; Church's wands seem to break this condition. After all: following Church, tapping $\emptyset$ with \textsf{complement} must reveal the universal set, $V$; but things we have not \emph{yet} discovered must be members of $V$; so one might well doubt whether this is ``reasonably constructive''!

Crucially, this doubt can be addressed. The complement of a bland set is never itself bland. So, where $\in$ is the primitive membership relation we use for bland sets, we can explicitly define an extended sense of ``membership'' which is suitable for both bland sets and their complements, along these lines:
\begin{align*}
	x \varin y &\colonequiv (\blandpred(y)\land x \in y) \lor \relexists{c}{\blandpred}(y =\textsf{complement}(c) \land x \notin c)
\end{align*}
Now, $\emptyset$ is bland, and \emph{nothing} is in$_{\varin}$ it; and $\emptyset$'s complement, $V$, is not bland, and \emph{everything} is in$_{\varin}$ it. By definition, though, the latter claim is exactly as constructive as the former. Doubt resolved!

Of course, this is only a quick sketch of how to make Church's \cite*{Church:STUS} theory look ``reasonably constructive'': it says nothing about the \textsf{cardinal}$_n$ wands, and omits a lot of other details. Still, it gives us all the right ideas, and the rest really is \emph{just} details (for which, see \S\ref{s:app:CUS}). The surprising upshot is that Church's theory can be regarded as a loosely constructive implementation of the Wand/Set Template. By my Main Theorem, it is therefore synonymous with a \ZF-like theory.\footnote{Unlike Examples \ref{ex:ConwayGames}--\ref{ex:Aczel}, this synonymy result is genuinely novel to this paper.}

Having tackled Church's theory, it is natural to ask whether we can grapple with other theories which are---at least at first blush---deeply non-constructive. Specifically, we might observe that Quine's \NF can be given a finite axiomatization which is based on the idea that there are certain ``starter'' sets, and that the other sets are formed from these by simple operations, such as taking complements.\footnote{See e.g.\ \citepossess{Holmes:ESTUS} axiomatization. The ``starter sets'' are given by the axioms of: Empty Set, Diagonal, Projections, Inclusion. The operations are given by: Complements, (Boolean) Unions, Set Union, Singletons, Ordered Pairs, Cartesian Products, Converses, Relative Products, Domains, Singleton Images. 
	
	Admittedly, some of the operations require two inputs (e.g.\ ordered pairs), whereas the wands in my Template act on just one input. However, the Template could easily be tweaked to allow for this. Moreover, the adjustment is scarcely likely to be necessary: if we wanted to input $a$ and $b$ into some ``binary'' wand, we could instead input the bland set $\{\{a\}, \{a,b\}\}$ into some ``unary'' wand (cf.\ Example \ref{ex:ConwayGames}).} 
Regarding the operations as wand-actions, we might then hope to be able to regard \NF \emph{itself} as an implementation of the Wand/Set Template. This might even lead to a consistency proof for \NF, or---as Church suggested---a ``synthesis or partial synthesis of'' \ZF with \NF.\footnote{\textcite[308]{Church:STUS}; Church mentions \citepossess{Hailperin:SAL} axioms rather than \citepossess{Holmes:ESTUS}.}

Clearly, the Wand/Set Template opens up some giddying possibilities. Inspired by them, my aim is to set the Template on a formal footing, and then establish my Main Theorem: \emph{any loosely constructive way of fleshing out the Wand/Set Template is synonymous with a \ZF-like theory.} (For better or worse, this Theorem ultimately restores some sobriety to our giddy state, since it seems to undermine the possibility of saying much about \NF using the Wand/Set Template.\footnote{See \citepossess{ForsterHolmes:SQCQS} argument for Sadie Kaye's Conjecture, that ``No extension of \NF is synonymous with any theory of wellfounded sets''. Their argument uses my Main Theorem; indeed, I wrote this paper in part to provide fuel for their argument.})

I will explain this synonymy result more precisely as I go. After some simple preliminaries (\S\S\ref{s:Notation}--\ref{s:Wevels}), I will explicate what it takes to flesh out the Wand/Set Template in a \emph{loosely constructive} way (\S\S\ref{s:Constructive}--\ref{s:WandSet}). This will allow me to present a general definition of an arbitrary wand/set theory (Definition \ref{def:WandSet}). I will then explain how any wand/set theory mechanically generates a theory which is ``\ZF-like'' in a precise sense (\S\ref{s:Strategy}). Having done this, I will prove the Main Theorem (\S\S\ref{s:WS}--\ref{s:int:biint}). I conclude by applying my wand/set apparatus to Church's \cite*{Church:STUS} theory (\S\ref{s:app:CUS}).

\section{Notation}\label{s:Notation}
In what follows, I restrict my attention to first-order theories.\footnote{If desired, both the idea of wand/set theories and the Main Theorem easily carry over to the second-order case; we simply replace schemes with axioms in the most obvious way.} For readability, I concatenate infix conjunctions, writing e.g.\ $a \subseteq r \in s \in t$ for $a \subseteq r \land r \in s \land s \in t$. I also use some simple abbreviations (where $\Psi$ can be any predicate with $x$ free, and $\lhd$ can be any infix predicate):
\begin{align*}
	(\forall x : \Psi)\phi&\coloneq \forall x(\Psi(x) \lonlyif \phi) & 	(\forall x \lhd y)\phi&\coloneq\forall x(x \lhd y \lonlyif \phi)\\
	(\exists x : \Psi)\phi&\coloneq\exists x(\Psi(x) \land \phi)& 
	(\exists x \lhd y)\phi&\coloneq\exists x(x \lhd y \land \phi)
\end{align*}
When I announce a result or a definition, I will list in brackets the axioms which I am assuming. 

Over the next few sections, my aim is to define the general idea of a wand/set theory. All such theories will have the following signature:
\begin{listbullet}
	\item[$\blandpred$:] a one-place predicate; intuitively, this picks out the bland sets
	\item[$\wandpred$:] a one-place predicate; intuitively, this picks out the wands
	\item[$\in$:] a two-place (infix) predicate; intuitively, $x \in a$ says that $a$ is a bland set with $x$ as a member
	\item[$\tappred$:] a three-place relation; intuitively, $\tapto{w}{a}{c}$ says that $c$ is found by tapping $a$ with $w$. (Nothing special is indicated by separating arguments with a slash rather than a comma; it simply improves readability.)
\end{listbullet}
When working in this signature, I use set-builder notation exclusively for \emph{bland sets}; so  $\Setabs{x}{\phi}$ is to be the bland set (if there is one) whose members are exactly those $x$ such that $\phi$. Similarly, $\emptyset$ will be the empty bland set; $\{\emptyset\}$ will be the bland set whose unique member is $\emptyset$, etc. For brevity, I will often describe something as simply ``bland'', to mean that it is a bland set. 

I now introduce three quasi-notational axioms:\footnote{For readability, I often omit outermost universal quantifiers on axioms; these should be read as ``the universal closure of\ldots''.}
\begin{listaxiom}
	\labitem{InNB}{WS:notein} $x \in a \lonlyif \blandpred(a)$
	\labitem{TapNB}{WS:notetap} $\tapto{w}{a}{c} \lonlyif \big(\wandpred(w) \land \lnot\blandpred(c)\big)$
	\labitem{TapFun}{WS:funtap} $\big(\tapto{w}{a}{c} \land \tapto{w}{a}{d}\big) \lonlyif c = d$
\end{listaxiom}
Axiom \ref{WS:notein} ensures that we reserve $\in$ for the notion of membership according to which bland sets have members. It does not follow that non-bland things have \emph{no} members; it just means that, if we want to allow non-bland things to have members, then we will have to define some richer notion of membership. (I foreshadowed this point in \S\ref{s:Motivations:Church}, when I outlined  ``$\varin$'', and I develop this in detail in  \S\S\ref{s:app:CUS:caveats}--\ref{s:app:CCC:example1.3}.) 

Axiom \ref{WS:notetap} is similar. It accords with the Template's prescriptions that (1) we only ever tap things with \emph{wands}, and that (2) bland things are never \emph{found} by tapping something with a wand but are always (and familiarly) found just by waiting until the stage immediately after you have found all their members. Again: this does not force us to deny that tapping something with a wand might ever ``yield'' something bland; we can allow this by defining a richer notion of wand-tapping. (Indeed, Church's \cite*{Church:STUS} theory allows that tapping something with $\textsf{complement}$ can \emph{yield} a bland entity, but it will always be a bland entity which was \emph{found} at an earlier stage; I discuss this in detail in \S\ref{s:app:CUS:caveats}.) 

Finally, \ref{WS:funtap} makes explicit the implicit point that $\tappred$ is functional (where defined). This licenses another bit of helpful notation: in what follows, I write $\tapop{w}{a}$ for the (necessarily unique) object, $c$, if there is one, such that $\tapto{w}{a}{c}$.

\section{Core axioms governing wevels}\label{s:Wevels}
My first serious task is to formalize the key idea of the Wand/Set Template: that, at each stage, we find (1) all bland sets of earlier-found objects, and (2) the result of tapping any earlier-found object with a wand (if the result is not bland). 

This key idea mentions stages. So I could start by positing \emph{stages} as special objects of a primitive sort, and laying down axioms about them. However, the tools and techniques developed in my \parencite*{Button:LT1} presentation of \LT, Level Theory, allow me to skip this step. Instead, I can just define some special bland sets which will go proxy for the Template's stages. These stage-proxies will be called \emph{wand-levels}, or \emph{wevels} for short.

To build up to the definition of a wevel, I begin with some basic axioms, which are the obvious versions of Extensionality and Separation for bland sets:
\begin{listaxiom}
	\labitem{Ext}{WS:ext} $\relforall{a}{\blandpred}\relforall{b}{\blandpred}\big(\forall x(x \in a \liff x \in b) \lonlyif a = b\big)$
	\labitem{Sep}{WS:sep}
	$\relforall{a}{\blandpred}\relexists{b}{\blandpred}\forall x\big(x \in b \liff (x \in a \land \phi)\big)$\\
	schema for any  $\phi$ not containing `$b$'
\end{listaxiom}
The next idea requires more explanation. Recall that wevels will be bland sets which go proxy for stages. Suppose we can define a formal predicate ``$\wevpred(x)$'', for \emph{$x$ is a wevel}, and a formal predicate ``$x \cfoundat r$'', for \emph{$x$ is found at wevel $r$}. Then we will want a wevel to be precisely the bland set of all earlier-found things, i.e.: 
\begin{align*}
	s 
	&= \Setabs{x}{\relexists{r}{\wevpred}(x \cfoundat r \in s)}
\end{align*}
That is our target; it just remains to give the explicit definition. Here it is:\footnote{Here, I highlight the use of two quasi-notational axioms. These ensure that, in the second disjunct, $x$ is non-bland and $r$ is bland and $w$ is a wand.}
\begin{define}[\ref{WS:notein}, \ref{WS:notetap}]\label{def:wevel}
	Define:
	\begin{align*}
		x \cfoundat r &\colonequiv \big(\blandpred(x) \land x \subseteq r\big) \lor \exists w(\exists b \in r)\tapto{w}{b}{x}\\
		x \incpot a&\colonequiv \exists r(x \cfoundat r \in a)
	\end{align*}\noindent
	Let $\cpot{a}$ be $\Setabs{x}{x \incpot a}$.\footnote{We do not assume at the outset that $\cpot{a}$ exists; but Lemma \ref{lem:WS:levof} proves it does. Note that $\cpot{a}$ is bland by definition.} 
	Say that $h$ is a wistory, $\wistpred(h)$, iff $\blandpred(h) \land (\forall a \in h)a = \cpot{(a \cap h)}$.\footnote{As usual, we define $a \cap h = \Setabs{x\in a}{x \in h}$, i.e.\ it is bland by definition; by \ref{WS:sep}, it will always exist.} Say that $s$ is a wevel, $\wevpred(s)$, iff $\relexists{h}{\wistpred}s = \cpot{h}$. \defendhere
\end{define}\noindent
To see why this is the correct definition of $\cfoundat$, recall that, at any wevel, we should find both 
all bland subsets of earlier wevels, and
all wand-taps of members of earlier wevels (unless this yields a bland set), and nothing else. 
The definition of $\wevpred$ is less immediately intuitive,\footnote{For those familiar with \LT (see \cite{Button:LT1}), here is a quick explanation of why Definition \ref{def:wevel} works. In ordinary \LT, where we have no wands, we replace $\cfoundat$ with $\subseteq$, since a level should just comprise all (bland) subsets of earlier levels. The remaining definitions are the same, using $\subseteq$ in place of $\cfoundat$.} but it is vindicated by these results (see \S\ref{s:WS} for more):
\begin{listbullet}
	\item[\emph{Lemma \ref{lem:WS:acc}.}] $s$ is a wevel iff $s = \cpot{\Setabs{r \in s}{\wevpred(r)}}$
	\item[\emph{Theorem \ref{thm:WS:wo}.}] The wevels are well-ordered by $\in$
\end{listbullet}\noindent
Unpacking the definition of $\cpot{}$, Lemma \ref{lem:WS:acc} says that a wevel is \emph{exactly} the bland set of all earlier-found things, which was precisely our target a couple of paragraphs above. Then Theorem \ref{thm:WS:wo} tells us that wevels are very nicely behaved. 

The two results just mentioned can be proved using a very weak theory, whose axioms are just \ref{WS:notein}, \ref{WS:notetap}, \ref{WS:ext} and all instances of \ref{WS:sep}. I will call this theory \CoreWev, since it proves the core results about wevels. 

The Wand/Set Template itself, though, tells us to go beyond \CoreWev. It says that every object is found at some stage. Since wevels are proxies for stages, this becomes the following axiom:
\begin{listaxiom}
	\labitem{Strat}{WS:strat} 
	$\forall a \relexists{s}{\wevpred}a \cfoundat s$
\end{listaxiom}
Since wevels are well-ordered, we can define the wevel ``at which $a$ is first found'':
\begin{define}[\CoreWev, \ref{WS:strat}]
	\label{def:wevof}
	For each $a$, let $\wevof{a}$ be the $\in$-least wevel $s$ such that $a \cfoundat s$. So: $a \cfoundat \wevof{a}$ and $\lnot\relexists{r}{\wevpred}a \cfoundat r \in \wevof{a}$. Let $\ordrank{a}$ be the ordinal $\beta$ such that $\wevof{a} = \wevof{\beta}$.\footnote{Here and throughout I define ordinals in the usual, von Neumann, fashion, i.e.\ as transitive (bland) sets well-ordered by $\in$; and I reserve the first few lower-case Greek letters to refer to ordinals.} \defendhere
\end{define}\noindent
Using just these axioms, we can prove several unsurprising but reassuring results about wevels; for example, that the $\alpha^\text{th}$ wevel is always $\wevof{\alpha}$.

\section{Loose constructivism}\label{s:Constructive}
We have seen that \CoreWev + \ref{WS:strat} ensures that our objects are arranged in a well-ordered hierarchy. However, this says nothing, yet, about the hierarchy's height, nor about the effects of tapping various objects with our various wands. 

This is deliberate. Given my aims (see \S\ref{s:Motivations}), I want to consider many different instances of the Wand/Set Template. Since different instances may describe hierarchies of various heights, or employ different wands, I cannot afford to be \emph{too} specific. 

That said, I do want to be a \emph{bit} more specific. Following Conway's lead (see \eqref{Conway:constructive} of \S\ref{s:Motivations:Conway}), I will restrict my attention to \emph{loosely constructive} ways of fleshing out the Wand/Set Template. The idea of loose constructivism is both conceptually subtle and formally fiddly, so the point of this section is to motivate and explicate the idea.

\subsection{The Picture}\label{s:LooseConstruct:Picture}
I will start with some very soft imagery. Imagine we are about to tap some object, $a$, with a wand, $w$, for the first time. We are trying to calculate what, if anything, we will find as a result of this wand-tap, i.e.\ whether $\tapop{w}{a}$ will exist and what it will be like. Now: \emph{what resources should we be allowed to invoke in our calculations?}

Evidently, we should be allowed to refer to any objects which we have \emph{already} found. They are, after all, lying around for us to inspect.  
But if we attempt to consider objects which have not \emph{yet} been found, we start to run the risk of vicious circularity. After all, suppose we think that whether $\tapop{w}{a}$ will exist depends upon which objects will exist \emph{later}. Well, what exists \emph{later} depends upon whether $\tapop{w}{a}$ will exist after the wand-tap we are about to perform, and we have run into circularity.

Such concerns are completely familiar from discussions of (im)predicativity in set theory. But they need not force us to adopt a fully predicative standpoint. After all, most of us are perfectly relaxed about the amount of impredicativity licensed by \ZF. So I propose that we be allowed \emph{exactly} that amount of impredicativity. 

Let me spell out what this proposal amounts to. Our Template tells us that we find new bland sets at every stage. Moreover, these bland sets behave just like \ZF-style sets. So, whatever objects we have already found (via whatever esoteric use of wands), we know that there \emph{will be} a ``bland-hierarchy'' which is built from them, intuitively treating the already-found objects as urelements and collecting them together into (bland) sets without any \emph{further} deployment of our wands. This ``bland-hierarchy'' is exactly as impredicative as ordinary \ZF; so we should be allowed to make reference to it, but nothing outside it.

Since things are getting complicated, a picture may help. (Compare this with the big ``V'' we draw to describe the ``ordinary'' set-theoretic hierarchy.)
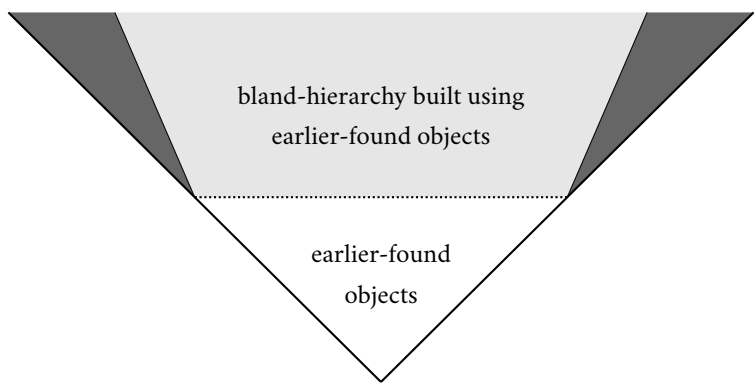
\begin{center}
	
	\
	
	\begin{tikzpicture}[scale=0.7]
		\fill[fill=black!60] (-3.5,3.5)--(-7,7)--(-5,7);
		\fill[fill=black!60] (3.5,3.5)--(7,7)--(5,7);
		\fill[fill=black!10] (-5,7)--(-3.5,3.5)--(3.5,3.5)--(5,7);
		\draw[thick] (0,0)--(-7,7);
		\draw[thick] (0,0)--(7,7);
		\draw (-5,7)--(-3.5,3.5); 
		\draw[thick, densely dotted] (-3.5,3.5)--(3.5,3.5);
		\draw (3.5,3.5)--(5,7);
		\node at (0,2)[align=center] {\selectfootfont earlier-found\\ \selectfootfont objects};
		\node at (0,5)[align=center] {\selectfootfont bland-hierarchy built using \\ \selectfootfont earlier-found objects};
	\end{tikzpicture}
	
	\
	
\end{center}
The unshaded area represents what we have already found; we are allowed to refer to this part of the hierarchy in our calculations, since it is perfectly transparent to us. The lightly-shaded area represents the bland-hierarchy constructed from the earlier-found objects. Although this bland-hierarchy is undiscovered as yet, its impredicativity is exactly of a piece with that of ordinary \ZF; we are allowed to refer to this part of the hierarchy in our calculations, since its behaviour has already been tamed. But the darkly-shaded area represents objects we will find (only) by future uses of wand-taps; we are barred from referring to this part of the hierarchy in our calculations since, as legend warns, here be viciously circular dragons. 

In what follows, I refer to this image as the \emph{Picture of Loose Constructivism}. This Picture guides me through the next couple of sections. But my immediate task is to make this (informal) Picture more precise, so that we can embed it within our formal theorizing. 

\subsection{Loosely constructive formulas}\label{s:LooseConstruct:Defined}
I start by explicitly defining the ``bland-hierarchy built from $\urbase$'', where $\urbase$ is any set of objects which we will treat as urelements. Intuitively, we treat $\urbase$ as given; we then resist the temptation to deploy any wands, but keep repeatedly forming all possible bland sets.

Evidently, this ``bland-hierarchy'' should be arranged into well-ordered layers, which we might call \emph{levels$_\urbase$}. Indeed, writing $\Tlev{\alpha}{\urbase}$ for the $\alpha^\text{th}$ level$_\urbase$, a little reflection on the ``intuitive'' idea just sketched indicates that the levels$_\urbase$ should behave thus:\footnote{Equivalently we could say: $\Tlev{0}{\urbase} = \urbase$; and $\Tlev{\alpha+1}{\urbase} = \powerset(\Tlev{\alpha}{\urbase}) \cup \urbase$; and $\Tlev{\alpha}{\urbase} = \bigcup_{\beta < \alpha}\Tlev{\beta}{\urbase}$ for limit $\alpha$.}
\begin{align*}
	\Tlev{\alpha}{\urbase} &= \Setabs{x}{\blandpred(x) \land (\exists \beta < \alpha)x \subseteq \Tlev{\beta}{\urbase}} \cup \urbase
\end{align*}
This identity would work perfectly well as a recursive \emph{definition} of the levels$_\urbase$.\footnote{See \S\ref{s:tolt:recursion} for a little more about the use of recursion definitions in the level-theoretic context.} However, it will be slightly easier for us if we instead, officially, define the levels$_\urbase$ non-recursively (noting that this is equivalent to the recursive approach, by Lemma \ref{lem:WS:Tlev}):\footnote{Compare this with Definition \ref{def:wevel} and my \parencite*[\S{}A]{Button:LT1} Level Theory with Urelements.}
\begin{define}\label{def:levelx}
	For any $\urbase$, say that:\footnote{Equivalently: $x \inpot{\urbase} a$ iff $x \incpot \urbase \lor \big(\blandpred(x) \land x \incpot a\big)$.}
	\begin{align*}
		x \inpot{\urbase} a&\colonequiv \big(\blandpred(x) \land \exists c(x \subseteq c \in a)\big) \lor x \in \urbase
	\end{align*}
	Let $\potur{\urbase}(a) \coloneq \Setabs{x}{x \inpot{\urbase} a}$. Say that $\histpred_\urbase(h)$ iff $\blandpred(h) \land (\forall x \in h)x = \potur{\urbase}(x \cap h)$. Say that $t$ is a level$_\urbase$, written $\levpred_\urbase(t)$, iff $(\exists h : \histpred_\urbase)t = \potur{\urbase}(h)$. Say that $\Vfrom{\urbase}(a)$ iff $\relexists{t}{\levpred_\urbase}a \in t$.
	\defendhere
\end{define}\noindent
Now the formal claim that $x$ satisfies  $\Vfrom{\urbase}$ precisely explicates the intuitive idea that $x$ falls somewhere in the ``bland-hierarchy built from $\urbase$''.

I can now explain how to embed the Picture of Loose Constructivism within our formal theorizing. Suppose we are at stage $s$, and some formula, $\phi$, is supposed to tell us something about what will happen later. So $\phi$ must deploy (only) resources which count as loosely constructive at $s$. In this case, I will say that $\phi$ is \emph{loosely bound} by $s$; and I will now define this notion precisely. 

Given the Picture, $\phi$ will be loosely bound by $s$ iff $\phi$ considers (only) the bland-hierarchy built from the objects found at $s$, i.e.\ the bland-hierarchy built from the objects \emph{in} $s^+$, where $s^+$ is the level immediately after $s$.\footnote{Recall: Theorem \ref{thm:WS:wo} tells us that the wevels are well-ordered. Lemma \ref{lem:WS:acc} tells us that $x \in s^+$ iff $x \cfoundat s$, i.e.\ iff $x$ is found at $s$. (And by Lemma \ref{lem:WS:levelpottrans}, if $x$ is found strictly before $s$, i.e.\ if there is a wevel $r$ such that $x \cfoundat r \in s$, then $x$ is found at $s$, i.e.\ $x \cfoundat s$.)} So $\phi$ should range (only) over $\Vfrom{s^+}$. Moreover, if $\phi$ mentions any wand-taps, i.e.\ if $\phi$ contains expressions of the form ``$\tapop{w}{x}$'', then, since $\tapop{w}{x}$ is never bland (by \ref{WS:notein}), $\phi$ must contain a guarantee that $\tapop{w}{x}$ is found no later than $s$; this is equivalent to insisting that $x$ itself is found \emph{strictly before} $s$,\footnote{See Lemma \ref{lem:WS:minnonbland}.} i.e.\ to insisting that $x \in s$. 

These ideas are rolled together in the following definition.
\begin{define}\label{def:looselybound}
	Let $s$ be any wevel and $\phi$ be any formula.
		
	An instance of a quantifier in $\phi$ is \emph{loosely bound} by $s$ iff it is restricted to $\Vfrom{{s^+}}$. So such instances are of the form $\relforall{x}{\Vfrom{s^+}}\psi$ or $\relexists{x}{\Vfrom{s^+}}\psi$. 
		
	An instance of a free variable $b$ is \emph{loosely bound} by $s$ in $\phi$ iff $\phi$ is a conjunction and $\Vfrom{s^+}(b)$ is one of its conjuncts.
	
	An instance of $\tapto{w}{a}{x}$ is \emph{loosely bound} by $s$ in $\phi$ iff it occurs in a binary conjunction whose other conjunct is $a \in s$, i.e.\ $(\tapto{w}{a}{x} \land a \in s)$.
	
	A formula is \emph{loosely bound} by $s$ iff all its quantifier-instances, free variables, and $\tappred$-instances are loosely bound by $s$.	\defendhere
\end{define}\noindent
To round off this discussion, it is worth giving particular consideration to the bland-hierarchy which is built by treating \emph{nothing} as an urelement, i.e.\ the denizens of $\Vfrom{\emptyset}$. These are bland sets, whose members are bland, whose members of members are bland\ldots and so on, all the way down. In a word, they are \emph{hereditarily bland}. These are important because, given the Picture, hereditarily bland sets are precisely the objects which are available to consider ``at the outset''. It turns out (Lemma \ref{lem:WS:V0}) that we can equivalently define hereditarily bland sets without mentioning $\Vfrom{\emptyset}$, and this will be my official definition:
\begin{define}\label{def:heredbland}
	Say that $a$ is \emph{hereditarily bland}, written  $\heblandpred(a)$, iff both $\blandpred(a)$ and $(\exists c \supseteq a)(\forall x\in c)(x \subseteq c \land \blandpred(x))$. \defendhere
\end{define}

\section{General definition of a wand/set theory}\label{s:WandSet}
In the last section, I formalized the notion of loose constructivism (Definition \ref{def:looselybound}). Using this, I can at last define what it takes for a theory to provide us with a loosely constructive implementation of the Wand/Set Template.

\subsection{The wands are hereditarily bland}\label{s:WandSet:Hebland}
I have said that any wand/set theory will need to have the three quasi-notational axioms (\ref{WS:notein}, \ref{WS:notetap} and \ref{WS:funtap} from \S\ref{s:Notation}) and the three axioms which arrange everything into well-ordered wevels (\ref{WS:ext}, \ref{WS:sep}, and \ref{WS:strat} from \S\ref{s:Wevels}). But I have said nothing, yet, about the wands. This obviously needs to addressed.

Different instances of the Template will have different wands. Still, we can still make a general comment. The Template tells us to tap everything we can, at every stage, with every wand. This means that every wand must be ``available'' to us at the outset. But only hereditarily bland sets are ``available'' to us at the outset (see the end of \S\ref{s:LooseConstruct:Defined}). So every wand must be hereditarily bland. We therefore adopt the axiom: 
\begin{listaxiom}
	\labitem{HebWands}{WS:wandsheb} $\relforall{w}{\wandpred}\heblandpred(w)$
\end{listaxiom}

\subsection{Weak claims about the hierarchy's height}\label{s:WandSet:Height}
Next, I lay down an axiom which makes a weak claim about the hierarchy's height: 
\begin{listaxiom}
	\labitem{WeakHeight}{WS:weakheight} $\Setabs{w}{\wandpred(w)}$ exists; and for every ordinal $\alpha$ both $\alpha \cup \{\alpha\}$ and  $\ordrank{(\Setabs{w}{\wandpred(w)})} + \alpha$ exist.\footnote{\label{fn:SpellOutSum}More fully spelled out, this would say: there is an ordinal which is isomorphic to the well-order obtained by taking a $0$-indexed copy of $\ordrank{(\Setabs{w}{\wandpred(w)})}$ and appending a $1$-indexed copy of $\alpha$.}
\end{listaxiom}
Given the ambient axioms, to say that $\Setabs{w}{\wandpred(w)}$ exists is to say that the hierarchy is sufficiently tall that all the wands are (eventually) found together at some wevel. The statement about an ordinal-sum existing just makes a further claim about the hierarchy's height. After all, each ordinal $\beta$ is first found at the $\beta^\text{th}$ wevel (see \S\ref{s:Wevels}). So we can control the hierarchy's height just by insisting that certain ordinals exist; and that is what \ref{WS:weakheight} does. 

Having understood roughly \emph{what} \ref{WS:weakheight} says, we need to ask \emph{why} we should adopt it. Candidly, \ref{WS:weakheight} is \emph{proof-generated}: it is the weakest height-principle which licenses my proof of synonymy!\footnote{Specifically, it is exactly what Definition \ref{def:conch} requires; see Lemmas \ref{lem:LTwand:Theta}--\ref{lem:LTwand:goodrecursions}.} I concede that it goes beyond anything mentioned in the Wand/Set Template---the Template itself allows that we might have precisely \emph{forty-two} stages, for example---but I note that it really is \emph{extraordinarily} weak for most mathematical purposes.\footnote{\ref{WS:weakheight} is immediately acceptable to anyone tempted by the vague but intuitive idea that the hierarchy ``shouldn't stop too early''. A good way to cash out this intuitive idea is via the precise claim: \emph{there is no unbounded map from any bland set into the hierarchy's levels}. This entails a version of Replacement (see \cites[323--6]{Shoenfield:AST}[93--95]{Incurvati:CS}[\S7]{Button:LT1}), and Replacement immediately entails that the sum of any two ordinals exists. It also yields a compelling argument for the existence of $\Setabs{w}{\wandpred(w)}$. Suppose, for reductio, that $\Setabs{w}{\wandpred(w)}$ does not exist; so, via our newly-minted Replacement-principle, there are proper-class many wands. Now consider an instance of the Template according to which tapping any object with any wand always yields a new object. So: at the outset, we find $\emptyset$. At the next stage, we find $\tapop{w}{\emptyset}$ for every wand $w$; by supposition, then, we find proper-class many entities. At the next stage we find all possible bland sets comprising only $\emptyset$ and the various $\tapop{w}{\emptyset}$'s; hence we find two-to-the-proper-class many bland sets; and this is a contradiction.} I adopt it going forward.

\subsection{Variable axioms for hereditarily bland sets}\label{s:WS:heredbland}
So far, I have specified some \emph{generic} axioms which will hold for any wand/set theory. I now need to start considering axioms which vary between theories.

To illustrate: Example \ref{ex:ConwayGames} only required a single wand (see \S\ref{s:Motivations:Examples}), but Church's \cite*{Church:STUS} theory treated $\omega$ as the set of wands (see \S\ref{s:Motivations:Church}). Saying that $\omega$ exists makes a specific demand on the height of the hierarchy, which is independent of the axioms we have laid down so far. But it is a demand which, with Church, we may well want to make.

We can make demands on the height of the hierarchy in other ways too. For example: suppose we add \ZF's Replacement scheme, restricting all the quantifiers in each scheme-instance to hereditarily bland sets; via perfectly familiar reasoning, this will ensure that our hierarchy is strongly inaccessible.

We might want to say this; we might not; it all depends on exactly how we want to implement the Wand/Set Template. In order to be maximally accommodating, I allow the following. A wand/set theory can take \emph{any} further sentences as axioms, provided that they are the result of taking a sentence, $\phi$, in the signature $\{\wandpred, \in\}$, and then restricting all of $\phi$'s quantifiers to hereditarily bland sets. (In what follows, I use the notation $\phi^\heblandpred$ for this.) 

Such axioms allow us to characterize the height of the hierarchy, precisely by characterizing the height of its hereditarily bland inner model. And any axiom which has been laid down in this way automatically conforms with loose constructivism. After all, the hereditarily bland sets are all available at the outset (see the end of \S\ref{s:LooseConstruct:Defined}). 

Such axioms will also allow us to say that the set of wands is $\omega$, or the least strongly inaccessible ordinal, or what-have-you. Of course, such axioms might lead us to inconsistency; there is no foolproof method for guarding against that. But nor should we expect such a method, any more than we expect one when we consider extensions of \ZF.

\subsection{Axioms for wand-taps}\label{s:WS:WandAxioms}
Finally, we need some axioms which tell us how wand-taps \emph{behave}. Different wand/set theories will, of course, say different things at this point; but we can still offer some useful general principles. 

Nothing I have said yet tells us whether (and when) we should find some object by tapping $a$ with $w$, i.e.\ whether $\tapop{w}{a}$ exists. We will address this by asking whether $a$ is in $w$'s ``domain of action'', written $\wanddom{w}{a}$. Of course, different theories will specify different ``domains'' for different wands. So, to specify a wand/set theory, we must explicitly define some formula, $\wanddom{w}{a}$, with all free variables displayed, and we then stipulate that $\wanddompred$ is a necessary and sufficient condition for finding something by wand-tapping:
\begin{listaxiom}
	\labitem{Making}{WS:make}
	$\wanddom{w}{a} \liff \exists c\ \tapto{w}{a}{c}$
\end{listaxiom}
We also need to know whether (and when) entities found by wand-tapping are identical, i.e.\ whether $\tapop{w}{a}=\tapop{u}{b}$. We will address this by asking whether $w$ and $a$ (in that order) are ``equivalent'' to $u$ and $b$ (in that order); we write this $\churchequiv{w}{a}{u}{b}$. Of course, different theories will use different notions of equivalence. So, to specify a wand/set theory, we must explicitly define some formula, $\churchequiv{w}{a}{u}{b}$, with all free variables displayed, and we then lay down a kind of quotienting principle:
\begin{listaxiom}
	\labitem{Equating}{WS:eq} 
	$\tapto{w}{a}{c} \lonlyif \big(\tapto{u}{b}{c} \liff \churchequiv{w}{a}{u}{b}\big)$ 
\end{listaxiom}
Otherwise put: if $\tapop{w}{a}$ exists, then $\tapop{w}{a} = \tapop{u}{b}$ iff $w,a$ is equivalent to $u,b$. 

As remarked: the axioms \ref{WS:make} and \ref{WS:eq} invoke make use of formulas, $\wanddompred$ and $\equivpred$, which a wand/set theory must explicitly define. So the last question is: \emph{What sorts of definitions are allowed?} Unsurprisingly, I want to maximize flexibility. So I want to say that any definitions are allowed, provided they meet these six constraints:
\begin{listdom}
	\item\label{dom:wand} $\wanddom{w}{a}$ entails $\wandpred(w)$;
	\item\label{dom:bound} $\wanddom{w}{a}$ is loosely bound by $\wevof{a}$;
	\item\label{dom:default} $\wanddompred$ is preserved under equivalence, in that we always have:
	\begin{align*}
		(\churchequiv{w}{a}{u}{b} \land \wanddom{w}{a}) &\lonlyif \wanddom{u}{b} 
	\end{align*}
\end{listdom}
\begin{listequiv} 
	\item\label{equiv:wand} $\churchequiv{w}{a}{u}{b}$  entails $\wandpred(w)$ and $\wandpred(u)$;
	\item\label{equiv:bound} $\churchequiv{w}{a}{u}{b}$ is loosely bound by $\wevof{a} \cup \wevof{b}$;
	\item\label{equiv:default} $\equivpred$ describes an equivalence relation, in that we always have:
	\begin{align*}
		\wandpred(w) &\lonlyif \churchequiv{w}{a}{w}{a}\\
		(\churchequiv{w}{a}{u}{b} \land \churchequiv{w}{a}{v}{c}) &\lonlyif \churchequiv{u}{b}{v}{c}
	\end{align*}
\end{listequiv}
These six constraints really just make explicit some of our implicit commitments. We insist on \domref{dom:wand} and \equivref{equiv:wand} because we are considering tapping things with wands, and only with wands. The Picture of Loose Constructivism forces \domref{dom:bound} and \equivref{equiv:bound} upon us: these constraints ensure that, when we (first) ask whether we should find anything by tapping $a$ with $w$, or whether $w,a$ is equivalent to $u,b$, we offer answers using only resources which are ``available'' to us there and then. The need for \equivref{equiv:default} is wholly obvious: we want $\equivpred$ to define a kind of \emph{equivalence}, which governs \emph{identity} between wand-taps via \ref{WS:eq}. The most interesting constraint is \domref{dom:default}. To see why we need it: suppose that $\wanddom{w}{a}$, so that tapping $\tapop{w}{a}$ should indeed exist, by \ref{WS:make}. Suppose also that $w,a$ is equivalent to $u,b$; then  $\tapop{w}{a}=\tapop{u}{b}$ by \ref{WS:eq}. But then, by \ref{WS:make}, we must also have that $\wanddom{u}{b}$.

This explains \emph{why} I have imposed conditions  \domref{dom:wand}--\domref{dom:default} and \equivref{equiv:wand}--\equivref{equiv:default}; it remains to explain \emph{how} I will impose them. 

Intuitively, the idea is as follows. We can lay down any ``attempted'' definitions of $\wanddompred$ and $\equivpred$ that we like; but if our ``attempted'' definitions ever stop behaving well, then we (thereafter) ``default'' to something which \emph{is} guaranteed to behave well. So, this does not limit our ability to define wand/set theories; it simply ensures that $\wanddompred$ and $\equivpred$ will satisfy \domref{dom:wand}--\domref{dom:default} and \equivref{equiv:wand}--\equivref{equiv:default}, even if our ``attempt'' misses the mark.

Here are the details. Let $D(w/a)$ be any formula which is loosely bound by $\wevof{a}$, and let $\dummyequiv{w}{a}{u}{b}$ be any formula which is loosely bound by $\wevof{a} \cup \wevof{b}$. (Think of these as our ``attempted'' definitions of $\wanddompred$ and $\equivpred$.) Now we define:\footnote{\label{fn:EquivRelationDefined}This is what it means, to say that $E$ restricted to $s$ is an equivalence relation: for any wands $x,y,z$, and for any $c,d,e\cfoundat s$, both $\dummyequiv{x}{c}{x}{c}$ and also $(\dummyequiv{x}{c}{y}{d} \land \dummyequiv{x}{c}{z}{e})\lonlyif \dummyequiv{y}{d}{z}{e}$. 

This is what it means, to say that $D$ restricted to $s$ is preserved under $E$: for any wands $x,y$, and for any $c,d \cfoundat s$, we have $(D(x/c) \land \dummyequiv{x}{c}{y}{d}) \lonlyif D(y/d)$.}
\begin{align}
	\wanddom{w}{a} &\colonequiv \wandpred(w) \land D(w/a)\nonumber
\intertext{and:}
	\churchequiv{w}{a}{u}{b} &\colonequiv 
	\wandpred(w) \text{ and }\wandpred(u) \text{ and either:}\nonumber\\
		&\phantom{{}\colonequiv{}}  w=u \land a = b\text{; or:}\label{eq:id}\\
		&\phantom{{}\colonequiv{}}  \dummyequiv{w}{a}{u}{b} \land {}\label{eq:E}\\
	&\hspace{3em} D\text{ restricted to } \wevof{a} \cup \wevof{b}\text{ is preserved under }E\land{}\nonumber\\
	&\hspace{3em} E\text{ restricted to }\wevof{a} \cup \wevof{b}\text{ is an equivalence relation}\nonumber
\end{align}
By \eqref{eq:id} and \eqref{eq:E}: if $D$ or $E$ stops behaving well at some level, then $\equivpred$ collapses to strict identity at all higher levels. This all works as required to ensure that \domref{dom:wand}--\equivref{equiv:default} hold (see Lemma \ref{lem:WS:wonderful}). So our only constraint is this: we insist that our theory's predicates $\wanddompred$ and $\equivpred$ are defined in this way, from suitable predicates $D$ and $E$.

\subsection{Definition of a wand/set theory}\label{s:WS:definition}
I am now in a position to define the general notion of a \emph{wand/set theory}. This makes precise the idea of a loosely constructive instance of the Wand/Set Template. 
\begin{define}\label{def:WandSet}
	A wand/set theory, \WS, is any theory in the signature $\{\blandpred$, $\wandpred$, $\mathord{\in}$, $\tappred\}$ which is axiomatized with all these axioms: 
	\begin{listn-0}
		\item \ref{WS:notein}, \ref{WS:notetap} and \ref{WS:funtap} (see \S\ref{s:Notation});
		\item \ref{WS:ext}, \ref{WS:sep} and \ref{WS:strat} (see \S\ref{s:Wevels});
		\item \ref{WS:wandsheb} and \ref{WS:weakheight} (see \S\S\ref{s:WandSet:Hebland}--\ref{s:WandSet:Height});
		\item \ref{WS:eq} and \ref{WS:make} (see \S\ref{s:WS:WandAxioms}), which use two explicitly defined formulas, $\wanddompred$ and $\equivpred$, that meet \domref{dom:wand}--\domref{dom:default} and \equivref{equiv:wand}--\equivref{equiv:default}.
	\end{listn-0}
	Moreover, if \WS is axiomatized with any \emph{further} axioms, then these are of the form $\phi^\heblandpred$, where $\phi$ is some sentence in the signature $\{\wandpred, \in\}$, and $\phi^\heblandpred$ simply involves restricting all quantifiers in $\phi$ to $\heblandpred$  (see \S\ref{s:WS:heredbland}).
		\defendhere
\end{define}\noindent
To recap: I have made some notational choices which impose a certain shape of formalism upon us; but these impose no substantial restrictions on our ability to flesh out the Template however we like (see e.g.\ \S\ref{s:Notation}). My method for arranging entities into a hierarchy of wevels does nothing more than eliminate the need to talk about stages in our formalism (see \S\ref{s:Wevels}). The Picture of Loose Constructivism yields a well-motivated explication of what it is to be ``reasonably constructive'' (see \S\ref{s:Constructive}). The constraints on the notion of a wand's ``domain of action'', and of ``equivalence'' between wand-taps, are as permissive as they can be (see \S\ref{s:WS:WandAxioms}). I admit that I have gone \emph{slightly} beyond the Template in making certain claims about the height of the hierarchy via \ref{WS:weakheight}, but these are \emph{extremely} weak claims (see \S\ref{s:WandSet:Height}). Summing all this up: when we study all wand/set theories, in the sense of Definition \ref{def:WandSet}, we study exactly the ``reasonably constructive'' ways to flesh out the Wand/Set Template (which accept a very weak principle about the hierarchy's height).

At this point, readers may wish to check that all of the examples discussed in \S\ref{s:Motivations} can be rendered \emph{as} wand/set theories. I leave Examples \ref{ex:ConwayGames}--\ref{ex:Aczel} to the reader, but exhaustively deal with Church's \cite*{Church:STUS} theory in \S\ref{s:app:CUS}.

\section{A \ZF-like theory, \LTwand }\label{s:Strategy}
The Main Theorem of this paper is that any wand/set theory is synonymous with a \ZF-like theory. I have explained what a wand/set theory is; I now need to articulate the ``\ZF-like theory'' in question. 

Of course, no single \ZF-like theory is synonymous with \emph{all} wand/set theories. Recall that wand/set theories can have any (arbitrary) axioms whose quantifiers are all restricted to hereditarily bland sets (see \S\ref{s:WS:heredbland}). What we need, then, is a canonical way to associate each wand/set theory, \WS, with a particular \ZF-like theory, which I will call \LTwand.

The choice of name references \LT, or Level Theory. The crucial point is that \LT is a pure theory of sets which says \emph{exactly} enough to ensure that the sets are arranged into well-ordered levels. This is achieved with a rather terse  definition:\footnote{See my \parencite*{Button:LT1}; and cf.\ Definition \ref{def:wevel} and \WS's axioms \ref{WS:ext}, \ref{WS:sep} and \ref{WS:strat}.}
\begin{define}\label{def:LTstuff}
	Define $x \blacktriangleleft a$ iff $\exists r(x \subseteq r \in a)$. Say that $a$ is potent iff $(\forall x \blacktriangleleft a)x \in a$. Let $\pot{a}$ be $\Setabs{x}{x \blacktriangleleft a}$. Say that $h$ is a history, $\histpred(h)$, iff $(\forall a \in h)a = \pot{(a\cap h)}$. Say that $s$ is a level, $\levpred(s)$, iff $\relexists{h}{\histpred}s = \pot{h}$. The axioms of (first-order) \LT are then:
	\begin{listaxiom}
		\labitem{\normalfont{Extensionality}}{lt:ext}$\forall a \forall b\big(\forall x(x \in a \liff x \in b) \lonlyif a =b\big)$
		\labitem{\normalfont{Separation}}{lt:scheme}$\forall a \exists b\forall x\big(x \in b \liff (\phi(x) \land x \in a)\big)$, for any $\phi$ not containing `$b$'
		\labitem{\normalfont{Stratification}}{lt:strat} $\forall a \relexists{s}{\levpred} a \subseteq s$ \defendhere
	\end{listaxiom}
\end{define}\noindent
Now, recall that I want to associate each wand/set theory, \WS, with some theory, \LTwand. This will be an extension of \LT, which makes the same claims about the height of the hierarchy as \WS, whilst using a new predicate to single out certain special objects as proxies for \WS's wands; we might as well just re-use the predicate ``$\wandpredlt$'' for this purpose. Here is the formal definition:
\begin{define}\label{def:LTwand}
	Where \WS is any wand/set theory, the theory \LTwand is defined as the theory with signature $\{\in, \wandpredlt\}$ and exactly these axioms:
	\begin{listn-0}
		\item\label{LTwand:LTbasis} each of the axioms of \LT (allowing \ref{lt:scheme}-instances to include ``$\wandpredlt$'').
		\item\label{LTwand:weakheight} an axiom, Pseudo-\ref{WS:weakheight}, which states: $\Setabs{w}{\wandpredlt(w)}$ exists; and for any ordinal $\alpha$, both $\alpha\cup\{\alpha\}$ and $\ordrank{(\Setabs{w}{\wandpredlt(w)})}+\alpha$ exist.\footnote{This is not \emph{verbatim} \WS's axiom \ref{WS:weakheight}, since (for example) \WS's axiom uses the predicate ``$\blandpred$''.}
		\item\label{LTwand:anyfurther} an axiom, $\phi$, for every \WS-axiom of the form $\phi^\heblandpred$ (see Definition \ref{def:WandSet}). \defendhere
	\end{listn-0}
\end{define}\noindent
My target result is then that every wand/set theory is synonymous with its associated \ZF-like theory. More precisely: 
\theoremstyle{generic}
\newtheorem*{maintheorem}{\textbf{Main Theorem}}
\begin{maintheorem}
	\LTwand and \WS are synonymous, for any wand/set theory, \WS.
\end{maintheorem}\noindent
What follows in the next four sections is nothing but a proof of the Main Theorem, without any further pauses for reflection. Throughout, I take \WS to be some fixed, arbitrary, wand/set theory. Here is my proof strategy.

In \S\ref{s:WS}, I will provide some characteristic results concerning \WS. 

In \S\ref{s:int:toc}, I will define a translation $\toc : \LTwand \functionto \WS$ and prove that it is an interpretation. 

In \S\ref{s:tolt}, I will define a translation $\tolt : \WS \functionto \LTwand$ and prove that it is an interpretation.

In \S\ref{s:int:biint}, I will show that these two interpretations constitute a bi-interpretation. By the Friedman--Visser Theorem, it follows that \LTwand and \WS are synonymous. 

\section{Elementary results within \WS}\label{s:WS}
The aim of this section is to prove generic results about \WS. I start with basic results about wevels (\S\ref{s:WS:wevels}), then consider some similar results about levels$_\urbase$ (\S\ref{s:WS:levelsur}). 

\subsection{Results about wevels}\label{s:WS:wevels}
The following definition is helpful for reasoning about wevels:
\begin{define}\label{def:pottranswand-transitivebland}Say that $a$ is \emph{wand-potent} iff $(\forall x \incpot a)x \in a$. Say that $a$ is \emph{wand-transitive} iff $(\forall x \in a)x \cfoundat a$. \defendhere
\end{define}\noindent 
Using this,\appref{proofsforwev} we can replicate in \WS a sequence of elementary results which hold for \LT (with tiny adjustments).
\footnote{\textcite[Results 3.2--3.10]{Button:LT1}.} I leave the details to the reader, but they culminate in two central results:
\begin{lem}[\CoreWev]\applabel{lem:WS:acc}
	$s$ is a wevel iff $s = \cpot{\Setabs{r \in s}{\wevpred(r)}}$.
\end{lem}
\begin{thm}[\CoreWev]\label{thm:WS:wo} 
	The wevels are well-ordered by $\in$.
\end{thm}\noindent
We can also show that the wevels behave nicely:
\begin{lem}[\CoreWev, \ref{WS:strat}]\applabel{lem:WS:levof}
	For all $a, b$, and all wevels $r, s$: 
	\setcounter{ncounts}{0}
	\begin{listn}
		\item\label{levofexists} $\wevof{a}$ and $\cpot{a}$ exist, with $\cpot{a} \subseteq \wevof{a}$
		\item\label{levofnotin} $a \notin \wevof{a}$
		\item\label{levofquick} $r\subseteq s$ iff $s\notin r$
		\item\label{levofidem} $s = \wevof{s}$
		\item\label{levofsubs} if $b \subseteq a$ and both are bland, then $\wevof{b} \subseteq \wevof{a}$
		\item\label{levofin} if $b \in a$, then $\wevof{b} \in \wevof{a}$
		\item\label{levnotin} $a \notin a$.
	\end{listn}
\end{lem}\noindent 
Arranging everything into well-ordered wevels has many nice consequences. For example, recall that wand-tapping is controlled by the defined predicates $\wanddompred$ and $\equivpred$ (see 
\ref{s:WS:WandAxioms}); now simple reasoning about wevels shows that my method for ensuring the ``good behaviour'' of $\wanddompred$ and $\equivpred$ works just as intended (I leave the proof to the reader):
\begin{lem}[\CoreWev, \ref{WS:strat}]\label{lem:WS:wonderful}
	Since $\wanddompred$ and $\equivpred$ are defined as in \S\ref{s:WS:WandAxioms}, all of \domref{dom:wand}--\domref{dom:default} and \equivref{equiv:wand}--\equivref{equiv:default} hold.
\end{lem}\noindent
Here is another important but simple consequence. If $c$ is not bland, there is some rank-minimal $a$ such that $c$ is obtained by tapping $a$ with some wand:
\begin{lem}[\CoreWev, \ref{WS:strat}]\label{lem:WS:min}
	If $c$ is not bland, there are $w$ and $a$ such $\tapto{w}{a}{c}$ and $\forall u \forall b(\tapto{u}{b}{c} \lonlyif \ordrank{a} \leq \ordrank{b})$. In that case, $\ordrank{a}+1=\ordrank{c}$.
\end{lem}\noindent 
It will later be helpful to restate Lemma \ref{lem:WS:min} in terms of $\equivpred$ rather than $\tappred$:
\begin{lem}[\CoreWev, \ref{WS:strat}, \ref{WS:eq}]\label{lem:WS:minnonbland}
	Say $\minirank{w}{a} \colonequiv \forall u \forall b\big(\churchequiv{w}{a}{u}{b} \lonlyif \ordrank{a}\leq \ordrank{b}\big)$. If $c$ is not bland, there are $w$ and $a$ such $\tapop{w}{a} = c $ and $\minirank{w}{a}$. In such case, $\ordrank{a}+1=\ordrank{c}$.
\end{lem}\noindent
This also yields an expected result: anything non-bland can be obtained from something bland using only finitely many wand taps: 
\begin{lem}[\CoreWev, \ref{WS:strat}, \ref{WS:funtap}]\applabel{lem:WS:fundamental}
	For any $a$, there is some $n$ such that $\relexists{b}{\blandpred}a = \tapop{w_n}\ldots\tapop{w_1}b$.
\end{lem}
\begin{proof}
	Here is the key idea. If $a$ is bland, we are done. Otherwise, $a = \tapop{w}{b}$ for some $b$ with $\ordrank{b} + 1 = \ordrank{a}$, using Lemma \ref{lem:WS:min}. If some such $b$ is bland, we are done; otherwise, keep going\ldots 
	this involves a strictly descending chain of ordinals, and hence we will find some bland set after finitely many steps.
	
	The only complication, in fact, comes from the need to define  ``$\tapop{w_n}\ldots\tapop{w_1}b$''.\footnote{\label{fn:WS:fundamental}This needs an explicit definition, since we are reasoning about ``finitely many steps'' in the \emph{object} language. See \S\ref{s:tolt:recursion} for comments on recursive definitions in this sort of setting. 
		
	What follows is a ``first pass'' definition, because its success in general requires that there be no last wevel. This can be secured using \ref{WS:weakheight}; but we can avoid the need for this by coding a few iterated wand-taps ``manually'' (cf.\ my \cite*[footnote 37]{Button:LT3}). For example, we can ask directly whether there are $v_0, \ldots, v_4$ such that $\tapop{v_4}\tapop{v_3}\tapop{v_2}\tapop{v_1}v_0 = a$ and $\ordrank{v_0} + 4 = \ordrank{a}$; and then search for a function $e$ and some $n$ with $\emph{BigTap}(e, n) = v_0$. (Note also that the proof requires no version of Choice.)} As a first pass: where $e$ is a function whose domain is a positive natural, define $\emph{BigTap}(e, 0) \coloneq e(0)$ and $\emph{BigTap}(e, n+1) \coloneq \tapop{e(n+1)}(\emph{BigTap}(e, n))$. Writing $e(i)$ as $e_i$, we can then think of $\emph{BigTap}(e, n)$ as $\tapop{e_n}\ldots\tapop{e_1}e_0$. 
\end{proof}\noindent
Note that there is no \emph{uniqueness} condition associated with Lemma \ref{lem:WS:fundamental}; there may be many (equally short) paths to the same object.

\subsection{Results about levels$_\urbase$ and hereditary-blandness}\label{s:WS:levelsur}
I now move from considering wevels to considering levels$_\urbase$. Unsurprisingly, given their similar definitions, wevels and levels$_\urbase$ behave very similarly, provided we assume that $\urbase$ has no self-membered elements.\footnote{I could ensure this via Lemma \ref{lem:WS:levof}\eqref{levnotin}. However, that would require the axiom \ref{WS:strat}, and it is worth noting that we do not need the full strength of that axiom.} In detail, let \CoreLev be the theory whose axioms are \ref{WS:ext}, all instances of \ref{WS:sep}, and $(\forall x \in \urbase)\urbase \notin \urbase$. Next, we use a new definition (cf.\ Definition \ref{def:pottranswand-transitivebland}):
\begin{define}
	Say that $a$ is \emph{potent}$_\urbase$ iff $\relforall{x}{\blandpred}((\exists c \notin \urbase)(x \subseteq c \in a) \lonlyif x \in a)$. Say that $a$ is \emph{transitive$_\urbase$} iff $(\forall x \notin \urbase)(x \in a \lonlyif x\subseteq a)$. \defendhere
\end{define}\noindent
The enthusiastic reader will now find it easy to prove these key facts:
\appref{proofsforlevurbase}
\begin{lem}[\CoreLev]\applabel{lem:C0:acc}  $s$ is a level$_\urbase$ iff 
	$s = \potur{\urbase}{\Setabs{r \in s}{\levpred_\urbase(r)}}$.
\end{lem}
\begin{lem}[\CoreLev]\label{lem:C0:wo} The levels$_\urbase$ are well-ordered by $\in$.
\end{lem}\noindent
%
%
These results allow us to talk about \emph{the} $\alpha^\text{th}$ level$_\urbase$, or $\Tlev{\alpha}{\urbase}$,\footnote{We do not assume this exists for every $\alpha$ and $\urbase$; the point is that $\Tlev{\alpha}{\urbase}$ is unique \emph{if} it exists.} which immediately vindicates the recursive characterization of the levels$_\urbase$ offered in \S\ref{s:LooseConstruct:Defined}:
\begin{lem}[\CoreLev]\label{lem:WS:Tlev}
	$\Tlev{\alpha}{\urbase} = \Setabs{x}{\blandpred(x) \land (\exists \beta < \alpha)x \subseteq \Tlev{\beta}{\urbase}} \cup \urbase$, for any $\alpha$.
\end{lem}\noindent
In \S\ref{s:LooseConstruct:Defined}, I also ``unofficially'' defined the hereditarily bland sets as the denizens of $\Vfrom{\emptyset}$, i.e.\ as members of any $\Tlev{\alpha}{\emptyset}$, whilst ``officially'' offering Definition \ref{def:heredbland}.  To see that these definitions are equivalent, we just need to introduce the idea of the hereditarily bland part of a set; the rest is simple:
\begin{define}\label{def:heblandpart}
	Let $\breve{a} \coloneq \Setabs{x \in a}{\heblandpred(x)}$. \defendhere
\end{define}
\begin{lem}[\CoreWev, \ref{WS:strat}]\applabel{lem:WS:hered}
	$a$ is hereditarily bland iff $a$ is bland and every member of $a$ is hereditarily bland.
\end{lem}
\begin{proof}
	\emph{Left to right.} When $c$ witnesses that $a$ is hereditarily bland, also $c$ witnesses that $x \in a$ is hereditarily bland. 
	\emph{Right to left.} In this case, $\breve{s}$ witnesses that $a$ is hereditarily bland, for any wevel $s \supseteq a$, using Lemma \ref{lem:WS:levelpottrans}.
\end{proof} 
\begin{lem}[\CoreWev, \ref{WS:strat}]
	\label{lem:WS:V0}
	$a$ is hereditarily bland iff $\Vfrom{\emptyset}(a)$. 
\end{lem}
\begin{proof}
	By induction, using Lemmas \ref{lem:WS:Tlev} and \ref{lem:WS:hered}. 
\end{proof}

\section{Interpreting \LTwand in \WS}\label{s:int:toc}
With some understanding of how \WS behaves, I will now show how to interpret \LTwand in \WS. We defined \LTwand precisely so that \WS's hereditarily bland sets will interpret \LTwand, but here are the details. I start by defining the translation:
\begin{define}
	Let $\toc$ be an identity-preserving translation, whose domain formula is $\heblandpred$, and with atomic clauses:\footnote{Note that \WS proves $\wandpred_\toc(x) \lonlyif \heblandpred(x)$ by \ref{WS:wandsheb}, and $x \in_\toc y \lonlyif \heblandpred(x)$ by Lemma \ref{lem:WS:hered}.} 
	\begin{align*}
		x \in_\toc y &\colonequiv x \in y \land \heblandpred(y)\\
		\wandpredlt_\toc(x) &\colonequiv \wandpred(x) 
	\end{align*}\defendnoskip
\end{define}\noindent
By \ref{WS:wandsheb}, \ref{WS:weakheight} and Lemma \ref{lem:WS:hered}, \WS delivers Pseudo-\ref{WS:weakheight}$^\toc$ (see Definition  \ref{def:LTwand}\eqref{LTwand:weakheight}). Moreover, where $\phi^\heblandpred$ is any \WS-axiom (see Definition \ref{def:LTwand}\eqref{LTwand:anyfurther}), $\phi$ is an \LTwand-axiom, and $\WS \proves \phi^\toc \liff \phi^\heblandpred$ by definition and Lemma \ref{lem:WS:hered}. So it just remains to show that \WS proves \ref{lt:ext}$^\toc$, \ref{lt:scheme}$^\toc$, and \ref{lt:strat}$^\toc$ (see Definition \ref{def:LTwand}\eqref{LTwand:LTbasis}). 
The first step is easy: by repeatedly applying Lemma \ref{lem:WS:hered} to \WS's axioms \ref{WS:ext} and \ref{WS:sep}, we find:
\begin{lem}[\WS]\applabel{lem:blt:helow:es} 
	Both \ref{lt:ext}$^\toc$ and \ref{lt:scheme}$^\toc$ hold.
\end{lem}\noindent
The key step to establishing Stratification$^\toc$ is to show that the levels$^\toc$ are the hereditarily bland parts of wevels (recall Definition \ref{def:heblandpart}); I leave this to the reader:\footnote{The proof is almost exactly as in my \parencite*[Lemmas C.4--C.5]{Button:LT3}.}
\begin{lem}[\WS]\applabel{lem:blt:levhelow} 
	$a$ is a level$^\toc$ iff $\relexists{s}{\wevpred} a = \breve{s}$. 
\end{lem}\noindent
Now Stratification$^\toc$ follows immediately via \ref{WS:strat}, so that:
\begin{prop}[\WS]\applabel{lem:WS:helowinterpret}
	$\toc: \LTwand \functionto \WS$ is an interpretation.
\end{prop}

\section{Interpreting \WS in \LTwand}\label{s:tolt}
I now show how to interpret \WS in \LTwand. This is much more challenging. The hardest part is to define the right translation (see \S\ref{s:tolt:def}); then it is just a matter of laboriously checking that the definition works as intended (see \S\S\ref{s:tolt:recursion}--\ref{s:tolt:wevlev}). 

\subsection{Defining the translation}\label{s:tolt:def}
It will help to start with an informal overview of how I will define the translation. 
I start by defining the things which will serve as the domain of interpretation; I call these things \emph{conches}, as a rough portmanteau of ``\emph{Con}way'' and ``\emph{Ch}urch''. Here are basic ideas about conches.
\begin{listbullet}
	\item[--] We need a way record the result of ``tapping'' some earlier-found object, $a$, with some ``wand'', $w$. (Scare-quotes are needed, since this is all under interpretation.) We do this by considering the pair $\tuple{w, a}$. 
	\item[--] We need to quotient such pairs under an equivalence relation (the interpretation of $\equivpred$); so a conch will typically be a set of such equivalent pairs. To ensure that we have a \emph{set} of such pairs, we use a version of Scott's trick, taking only the equivalent pairs of minimal rank (for some suitable notion of rank).
	\item[--] Finally, we need a way to consider ``bland'' entities; we take them to have the form $\{\tuple{\emptyset, a}\}$, which I write $\carrier{a}$. (So I will insist that $\emptyset$ is not a ``wand'', but a marker of ``blandness''.) 
\end{listbullet}
Where $\tolt$ is the translation that I will define, I will set things up so that $x \in_\tolt \carrier{a}$ iff $x \in a$. This allows me to compress most of the above into the following statement of intent (which I ultimately prove in \S\ref{s:tolt:wevlev}): 
\begin{listbullet}
	\item[\emph{Corollary \ref{cor:LTS:target}.}] $c$ is a conch iff either: 
	\begin{listn-0}
		\item[\eqref{conch:taxonomy:bland}] $c = \carrier{s}$, for some $s$ whose members are all conches; or
		\item[\eqref{conch:taxonomy:nonbland}] $c = \Setabs{\tuple{w,a}}{\taptolt{w}{a}{c}\land \miniranktolt{w}{a}} \neq \emptyset$
	\end{listn-0}
\end{listbullet}\noindent
Condition \eqref{conch:taxonomy:bland} will hold iff $c$ is bland$_\tolt$. Condition \eqref{conch:taxonomy:nonbland} will hold iff $c$ is non-bland$_\tolt$, whereupon $\tuple{w,a} \in c$ iff $c$ is found by ``tapping'' $a$ with ``wand'' $w$ (and not by ``tapping'' anything of lower rank than $a$).

It should be clear that, if we are careful, we will indeed have interpreted \WS in \LTwand. It just remains to thrash through the details. 

In \WS, the hereditarily bland sets form a sort of ``spine'' around which the rest of the hierarchy is constructed. So I will start by explaining how to simulate \WS's hereditarily bland sets within \LTwand. In fact, this simulation is dictated by the choices I made above:
\begin{define}[\LTwand]\label{def:Theta}
	Let $\carrier{a} \coloneq \{\tuple{\emptyset, a}\}$. When $a = \{\tuple{\emptyset, c}\} = \carrier{c}$, let $\uncarrier{a} \coloneq c$, so that $\uncarrier{\carrier{c}} = c$. Recursively let $\Theta a \coloneq \carrier{\Setabs{\Theta x}{x \in a}}$.
	\defendhere
\end{define}\noindent
Note that \LTwand proves both that $\Theta$ is injective and that $\emptyset \neq \Theta x$ for any $x$ (see Lemma \ref{lem:LTwand:Theta}). So, having used $\emptyset$ to code ``bland'' conches, I can say that a ``wand'' is any $\Theta$-image of a $\wandpredlt$. This motivates the following:
\begin{define}\label{def:tostep}
	For any $w,x,y$:
	\begin{align*}
		\blandpred_\tostep(y)&\colonequiv \exists z\ y = \carrier{z}\\
		x \in_\tostep y&\colonequiv \blandpred_\tostep(y) \land x \in \uncarrier{y}\\
		\wandpred_\tostep(w) &\colonequiv \relexists{v}{\wandpredlt}w=\Theta v 
	\end{align*} \defendnoskip
\end{define}\noindent
This definition forms the heart of my translation $\tolt : \LTwand \functionto \WS$. However, I cannot complete that translation until I have defined a suitable domain formula, i.e.\ until I have said what the conches will be. Unfortunately, this will takes some time. 

In \S\ref{s:Constructive}, I explicated the Picture of Loose Constructivism. The general idea was to treat the objects found at some stage of the hierarchy as urelements, and then to build a ``bland-hierarchy'' from those urelements. My first job is to translate this idea into \LTwand; I do this by recursively defining $\Ufrom{\urbase}$ (cf.\ Definition \ref{def:levelx} and Lemma \ref{lem:WS:Tlev}):
\begin{define}[\LTwand]\label{def:internalhierarchy}
	For any set $\urbase$, recursively define (where this exists):
	\begin{align*}
		\Ulev{\alpha}{\urbase} &\coloneq \Setabs{\carrier{x}}{(\exists \beta < \alpha)x \subseteq \Ulev{\beta}{\urbase}} \cup \urbase
	\end{align*}
	Say that $\Ufrom{\urbase}(x)$ iff $\exists \alpha\ x \in \Ulev{\alpha}{\urbase}$. 
	\defendhere
\end{define}\noindent
%
The next task is to specify what the stages of our ``\WS-hierarchy'' will look like. Specifically, I must define the conches and associate them with some notion of \emph{rank} according, intuitively, to the stage at which they are first found.

To effect this, I define, for each ordinal $\sigma$, the sets $\conchstage{\sigma}$, $\domstage{\sigma}$, and $\equivstage{\sigma}$. Here are the guiding ideas:
\begin{listbullet}
	\item[--] $\conchstage{\sigma}$ is a set whose members are the conches of rank $\leq \sigma$ (mnemonically: Conch) 
	\item[--] $\belowstage{\sigma}$ is a whose members are the conches of rank $< \sigma$ (mnemonically: Below)
	\item[--] I write $\conchrank{a} =  \sigma$ to indicate that $a$ is a conch of rank $\sigma$.
	\item[--] $\domstage{\sigma}$ and $\equivstage{\sigma}$ respectively simulate $\wanddompred$ and $\equivpred$ for conches of rank $\leq \sigma$. Indeed, we will ultimately set things up so that:
	\begin{align*}
		\tuple{w, a} \in \domstage{\sigma}&\text{ iff $\conchrank{a} \leq \sigma$ and $\wanddomtolt{w}{a}$}\\
		\tuple{w,a,u,b} \in \equivstage{\sigma}&\text{  iff $\conchrank{a},\conchrank{b} \leq \sigma$ and $\churchequivtolt{w}{a}{u}{b}$}
	\end{align*}
\end{listbullet}
The guiding principles are quite clear; alas, the formal definition itself is quite long. Indeed, it is long enough for me to need to interrupt it with commentary to make it intelligible. 
Here is how it starts.
\begin{define}[\LTwand]\label{def:conch}
	For each ordinal $\sigma$, let: 
	\begin{align*}
		\conchstage{\sigma} &\coloneq \Setabs{\carrier{s}}{s \subseteq \belowstage{\sigma}} \cup \Setabs{c}{(\exists \tau < \sigma)(c\text{ is a $\tau$-tap})}
	\end{align*}
	where I write $\belowstage{\sigma} \coloneq \bigcup_{\tau < \sigma}\conchstage{\tau}$. Say that $x$ is a conch, written $\conchpred(x)$, iff $\exists \sigma\ x \in \conchstage{\sigma}$. For each conch $x$, let  $\conchrank{x}$ be the least $\sigma$ such that $x \in \conchstage{\sigma}$. \tbchere
\end{define}\noindent
This formal definition is incomplete (hence the ``\textsc{t}o \textsc{b}e \textsc{c}ontinued''). In particular, I have used the phrase ``a $\tau$-tap'', which I have not yet defined. However, I should pause to explain the definition for $\conchstage{\sigma}$. The component $\Setabs{\carrier{s}}{s \subseteq \belowstage{\sigma}}$ collects ``bland'' conches: when $s$ is a set of conches, $\carrier{s}$ will act as a ``bland'' conch with exactly the members of $s$ as ``members'', and $\carrier{s}$'s rank will be the supremum of its ``members''. Then $\Setabs{c}{(\exists \tau < \sigma)(c\text{ is a }\tau\text{-tap})}$ will collect the ``non-bland'' conches. 

To complete Definition \ref{def:conch}, then, I must define what I mean by a $\tau$-tap; or, changing index for convenience, a $\sigma$-tap. Roughly speaking, a $\sigma$-tap should be a ``non-bland'' conch of rank $\sigma+1$ which was obtained by tapping some conch of rank $\sigma$ with some ``wand''. More precisely:
\begin{define*}[\textbf{\ref{def:conch}}, cont.]
	For each ordinal $\sigma$, say that $c$ is a $\sigma$-tap iff $c$ is a non-empty set of ordered pairs such that, for all $\tuple{w, a} \in c$, these three conditions hold:
	\begin{align}
		\eqtag{\textsf{tapr}}\label{conch:c:rank}
		&\conchrank{a} = \sigma \text{ and }\text{if }\tuple{w,a, u,b} \in \equivstage{\sigma}\text{, then }\conchrank{b} = \sigma \\
		\eqtag{\textsf{tapd}}\label{conch:c:dom} 
		&\tuple{w, a} \in \domstage{\sigma}\\
		\eqtag{\textsf{tape}}\label{conch:c:equiv} 
		&\tuple{w,a, u,b} \in \equivstage{\sigma}\text{ iff }\tuple{u, b} \in c
	\end{align}\tbchere
\end{define*}\noindent
Again, I will pause to explain this. Roughly: 
\eqref{conch:c:rank} says that $a$ has rank $\sigma$ and that this rank is as small as possible; 
\eqref{conch:c:dom} says that $a$ is within $w$'s ``domain of action''; and
\eqref{conch:c:equiv} says that $c$ is an equivalence class of minimally-ranked entities. But of course these clauses can only have this effect if we have appropriately defined $\domstage{\sigma}$ and $\equivstage{\sigma}$. So we offer:
\begin{define*}[\textbf{\ref{def:conch}}, cont.]
	For each ordinal $\sigma$:
	\begin{align*}
		\domstage{\sigma}& \coloneq \Setabs{\tuple{w, a}}{
			\conchrank{a} \leq \sigma \text{ and }\wandpred_{\paramint{\sigma}}(w) \text{ and }
			\wanddomint{w}{a}{\paramint{\sigma}}}\\
		\equivstage{\sigma}&\coloneq \{\tuple{w,a,u,b}: 
			\conchrank{a}, \conchrank{b} \leq \sigma\text{ and }\wandpred_{\paramint{\sigma}}(w)\text{ and }\wandpred_{\paramint{\sigma}}(u)\text{ and }\\
			&\phantom{{}\coloneq{}\{\tuple{w,a,u,b}:{} }
			\churchequivint{w}{a}{u}{b}{\paramint{\sigma}} \lor (w=u \land a =b)\}
	\end{align*}\tbcnoskip
\end{define*}\noindent
Interrupting again: whatever the exact meaning of the $\paramint{\sigma}$-translation, it should be clear how the definitions of $\domstage{\sigma}$ and $\equivstage{\sigma}$ are supposed to fit with our intentions (see the comments just before the start of Definition \ref{def:conch}). 
So it just remains to define each $\paramint{\sigma}$:
\begin{define*}[\textbf{\ref{def:conch}}, finished.]
	For each ordinal, $\sigma$, define an identity-preserving translation $\paramint{\sigma}$ as follows (using $\Delta_I$ for $I$'s domain formula):
	\begin{align*}
			\Delta_{\paramint{\sigma}}(x) &\colonequiv \Ufrom{\conchstage{\sigma}}(x)\\
			\blandpred_{\paramint{\sigma}}(y)&\colonequiv \blandpred_\tostep(y) \land \Delta_{\paramint{\sigma}}(y)\\
			x \in_{\paramint{\sigma}} y&\colonequiv x \in_\tostep y \land \Delta_{\paramint{\sigma}}(y)\\
			\wandpred_{\paramint{\sigma}}(w) &\colonequiv \wandpred_\tostep(w) \land \Ufrom{\emptyset}(w)\\
		\taptosigma{w}{a}{c} &\colonequiv 
		\tuple{w,a} \in \domstage{\conchrank{a}} \land \conchrank{a} < \sigma \land
		\conchrank{c}\leq \sigma \land {}\\
		&\phantom{{} \colonequiv {}}
		\exists u \exists b\big(
			\tuple{w,a,u,b} \in \equivstage{\conchrank{a}}  \land
			\tuple{u,b} \in c)\big)  \nonumber 
	\end{align*}\defendnoskip
\end{define*}\noindent
The domain formula, $\Ufrom{\conchstage{\sigma}}$, ensures that $\paramint{\sigma}$ quantifies only over entities which are ``available'' at the $\sigma^\text{th}$ ``wevel''. We then translate $\blandpred$ and $\in$ as suggested by Definition \ref{def:tostep}, whilst ensuring that we are dealing with entities in the intended ``domain''. The clause for $\wandpred_{\paramint{\sigma}}$ is similar but more restrictive, and reflects the intention that ``wands are hereditarily bland'' (see Proposition \ref{prop:LTwand:iso} for more.) Finally, the idea for $\tappred_\paramint{\sigma}$ is as follows: if $c$ is found by ``tapping'' $a$ with $w$, then $c$ is a $\tau$-tap, for some $\tau < \sigma$. 

Definition \ref{def:conch} is finally complete. And I can use it to define the $\tolt$-translation as the ``limit'' of the $\paramint{\sigma}$-translations.
\begin{define}[\LTwand]\label{def:toltsecond}
	Define $\tolt$ as an identity-preserving translation as follows:
	\begin{align*}
		\Delta_\tolt(x) &\colonequiv \conchpred(x)\\
		\blandpred_{\tolt}(y)&\colonequiv \blandpred_\tostep(y) \land \Delta_{\tolt}(y)\\
		x \in_{\tolt} y&\colonequiv x \in_\tostep y \land \Delta_{\tolt}(y)\\
		\wandpred_{\tolt}(w) &\colonequiv \wandpred_\tostep(w) \land \Ufrom{\emptyset}(w)\\
		\taptolt{w}{a}{c} &\colonequiv 
		\tuple{w,a} \in \domstage{\conchrank{a}} \land \Delta_\tolt(c) \land {}\\
		&\phantom{{} \colonequiv {}} 
			\exists u \exists b\big(
			\tuple{w,a,u,b} \in \equivstage{\conchrank{a}}  \land
			\tuple{u,b} \in c)\big) 
	\end{align*}\defendnoskip
\end{define}\noindent
I hope that my informal explanations make it plausible that $\tolt$ is, indeed, an interpretation. The remainder of this section simply verifies this fact. 

\subsection{Confirming the recursions}\label{s:tolt:recursion}
I must start by commenting on my use of recursive definitions. In brief: recursive definitions make good sense in \LTwand, but the defined term-function may be \emph{partial}. Here is a slightly fuller explanation of this point. 

Given a term, $\textbf{t}$, we can stipulate (as usual) that $f$ is an $\alpha$-approximation to $\textbf{t}$ iff $\domain{f} = \alpha$ and $(\forall \beta < \alpha)f(\beta) = \textbf{t}(f\restrictto{\beta})$. Ordinal induction will establish that $\alpha$-approximations agree for any shared arguments, so we can write $f^\textbf{t}_\alpha$ for the unique $\alpha$-approximation, if any exists. Where we then explicitly define
	\begin{align*}
		\textbf{r}(\alpha) &\coloneq f^\textbf{t}_{\alpha+1}(\alpha)
\intertext{we can prove a version of the general recursion theorem:}
	\textbf{r}(\alpha) &= \textbf{t}(\Setabs{\tuple{\beta, \textbf{r}(\beta)}}{\beta< \alpha})&&\text{if $f^\textbf{t}_{\alpha+1}$ exists}\\
	\textbf{r}(\alpha)&\text{ is undefined} &&\text{otherwise}  
	\end{align*}
All of this is as in \ZF. The sole difference is that \ZF proves that $f^\textbf{t}_{\alpha+1}$ must exist for any $\alpha$, but \LTwand might not prove this. So $\textbf{r}$ may be \emph{partial}.

It is, then, important to confirm that many of the notions defined in \S\ref{s:tolt:def} are \emph{total}. To do this, we can take a slight shortcut. Pseudo-\ref{WS:weakheight} ensures that there is no last level. So, where $\textbf{r}$ is recursively-defined from $\textbf{t}$ as above, to show that every $\textbf{r}(\alpha)$ exists, we need only show that there must always be an ordinal corresponding to ``what $\textbf{r}(\alpha)$'s rank would have to be''.

We will first use this strategy to show that $\Ulev{\alpha}{\emptyset}$ and $\Theta$ are total.
\begin{lem}[\LTwand]\label{lem:LTwand:Theta}
	For any $\alpha$, $\Ulev{\alpha}{\emptyset}$ exists.  
	For any $a$, both $\Theta a \in \Ulev{\ordrank{a}+1}{\emptyset}$ and $\ordrank{a} = \conchrank{\Theta a}$. 
	Indeed, the map $\Theta : V \functionto \Ufrom{\emptyset}$ is a total bijection.
\end{lem}
\begin{proof}
	Where $\alpha$ is any ordinal, write it as $\alpha = \lambda_\alpha + n_\alpha$, with $\lambda_\alpha$ either a limit ordinal or $0$, and $n_\alpha$ a natural number. By induction, $\ordrank{(\Ulev{\alpha}{\emptyset})} = \lambda_\alpha + 4n_\alpha$, so each $\Ulev{\alpha}{\emptyset}$ exists.
	
	A simple $\in$-induction shows that $\Theta$ is injective. Suppose for induction that $b \in \Ulev{\ordrank{b}+1}{\emptyset}$ whenever $b \in a$. As $\ordrank{a} = \lsub_{b\in a}{\ordrank{b}}$, we have $\uncarrier{\Theta a}\subseteq \Ulev{\ordrank{a}}{\emptyset}$, so $\Theta a \in \Ulev{\ordrank{a}+1}{\emptyset}$. Suppose also that $\ordrank{b} = \conchrank{\Theta b}$ whenever $b \in a$; then:
	\begin{align*}
		\ordrank{a}
		= \lsub_{b \in a}{\ordrank{b}} 
		= \lsub_{b \in a}{\conchrank{\Theta b}} 
		= \lsub_{x \in \uncarrier{\Theta a}} \conchrank{x}
		=\conchrank{\Theta a}
	\end{align*}	
	It follows that $\Theta : V \functionto \Ufrom{\emptyset}$ is total. 
	For surjectivity, suppose for induction that if $x \in \Ulev{\gamma}{\emptyset}$ then $\Theta^{-1}x$ exists, whenever $\gamma < \beta$. Let $b \in \Ulev{\beta}{\emptyset}$; then $a = \Setabs{\Theta^{-1}x}{x \in \uncarrier{b}}$ exists by Separation, and and $\Theta a = \carrier{\uncarrier{b}} = b$.
\end{proof}\noindent
Note that \LTwand might \emph{not} prove that $\Ulev{\alpha}{\urbase}$ exists for every $\alpha$ and every $\urbase$. This is not a problem; we do not in general require that it does. (Similarly: \WS might not prove that $\Tlev{\alpha}{\urbase}$ exists for every $\urbase$ and $\alpha$.) Unsurprisingly, though, we \emph{do} require that $\conchstage{\alpha}$ exists for every $\alpha$. And it provably does. Indeed (Pseudo-)\ref{WS:weakheight} was chosen precisely to secure this fact.
\begin{lem}[\LTwand]\label{lem:LTwand:goodrecursions}
	For any $\alpha$, each of  $\belowstage{\alpha}, \conchstage{\alpha}, \domstage{\alpha}$ and $\equivstage{\alpha}$ exists. 
\end{lem}
\begin{proof}
	It suffices to establish the result for $\conchstage{\alpha}$. Let $\Omega = \lub_{w\text{ a wand}}\ordrank{(\Theta w)}$. Writing each $\alpha$ as $\lambda_\alpha + n_\alpha$, I claim that: 
		$$\ordrank{\conchstage{\alpha}} \leq \Omega + \lambda_\alpha + 4n_\alpha + 4$$
	By Lemma \ref{lem:LTwand:Theta}, there is $m$ such that $\Omega = \ordrank{(\Setabs{w}{\wandpred(w)})} + m$; so $\conchstage{\alpha}$'s required rank exists by my claim and Pseudo-\ref{WS:weakheight}. 
	
	It just remains to establish my claim. This is by induction. Trivially, $\ordrank{\conchstage{0}} = 4$. 
	
	\emph{Successor case.} Suppose $\ordrank{\conchstage{\alpha}} \leq \Omega + \lambda_\alpha + 4n_\alpha + 4$. If $s \subseteq \belowstage{\alpha+1}=\conchstage{\alpha}$, then:
	\begin{align*}
		\ordrank{(\carrier{s})} 
			&\leq \Omega + \lambda_\alpha + 4n_\alpha + 7
	\intertext{If $c$ is a $\beta$-tap with $\beta \leq \alpha$, then $c$'s members are of the form $\tuple{\Theta w, d}$, with $w$ a wand and $\ordrank{d} \leq \Omega + \lambda_\alpha + 4n_\alpha + 3$, and since $\ordrank{(\tuple{\Theta w, d})} = \max(\ordrank{\Theta w}, \ordrank{d}) + 2$, we obtain:}
		\ordrank{c} 
			&\leq 
			\Omega + \lambda_\alpha +4n_\alpha+6
	\end{align*}
	So $\ordrank{\conchstage{\alpha}} \leq \Omega + \lambda_\alpha + 4n_\alpha + 8$, as required. 
		
	\emph{Limit case.} Using the induction hypothesis,  
	if $s \subseteq \belowstage{\alpha}$ then $\ordrank{(\carrier{s})} \leq \lub_{\beta < \alpha}(\Omega + \lambda_\beta + 4n_\beta + 4) + 3 = \Omega + \alpha + 3$. So $\ordrank{(\conchstage{\alpha})} \leq \Omega + \alpha + 4$, as required. 
\end{proof}

\subsection{$\star$-translations and  blandness$_\star$}\label{s:tolt:star}
I will now look at the $\tolt$ and $\paramint{\sigma}$ translations. For readability, I use $\star$ to stand ambiguously for either $\tolt$ or any $\paramint{\sigma}$. So if I assert that $\phi^\star$, I mean that $\phi^\tolt$ and $\phi^{\paramint{\sigma}}$ for any $\sigma$. Recall that $\Delta_\star$ is always $\star$'s domain-formula. 

I start with some elementary results about bland$_\star$ conches. These results are pretty immediate consequences of the definitions (and previous results) and I leave their proofs to the reader:
\begin{lem}[\LTwand]\applabel{lem:LTwand:note}
	Both \ref{WS:notein}$^\star$ and \ref{WS:notetap}$^\star$ hold.
\end{lem}
\begin{lem}[\LTwand]\applabel{lem:LTwand:xor}
	If $\Delta_\star(c)$, then: $c$ is not bland$_\star$ iff $c$ is some $\sigma$-tap.
\end{lem}
\begin{lem}[\LTwand]\applabel{lem:LTwand:staristrans} If $\blandpred_\star(a)$ then $\Delta_\star(a)$. If $b \in_\star a$, then $\Delta_\star(b)$ and $\Delta_\star(a)$. If $\taptostar{w}{a}{c}$, then $\Delta_\star(w)$, $\Delta_{\star}(a)$ and $\Delta_\star(c)$. So $\star$ is a translation.  
\end{lem}
\begin{lem}[\LTwand]\applabel{lem:LTwand:transint}
	If $\Delta_\star(\carrier{a})$, then:
	\begin{listn-0}
		\item\label{n:rank} $\conchrank{\carrier{a}} = \lsub_{x\in a}\conchrank{x}$
		\item\label{n:in} $b \in a$ iff $b \in_\star \carrier{a}$
		\item\label{n:sub} $c \subseteq a$ iff $\carrier{c} \subseteq^\star \carrier{a} \land \Delta_\star(\carrier{c})$
	\end{listn-0}
\end{lem}
\begin{lem}[\LTwand]\applabel{lem:LTwand:extsepstar}
	\ref{WS:ext}$^\star$ and \ref{WS:sep}$^\star$ (scheme) hold
\end{lem}\noindent
I now want to consider hereditarily bland$^\star$ conches. In particular, I want to show that  $\Theta$ is an ``isomorphism'' to the hereditarily bland$^\star$ conches:
\begin{prop}[\LTwand]\applabel{prop:LTwand:iso}
	For any $a, b, w$:
	\begin{listn-0}
		\item\label{Theta:Ulev} $\heblandpred^\star(b)$ iff $\Ufrom{\emptyset}(b)$, so that $\Theta: V \functionto \heblandpred^\star$ is a total bijection
		\item\label{Theta:iso:in} $a \in b$ iff $\Theta a \in^{\toc\star} \Theta b$ iff $\Theta a \in_\star \Theta b$.
		\item\label{Theta:iso:wandpredlt} $\wandpredlt(w)$ iff $\wandpredlt^{\toc\star}(\Theta w)$ iff $\wandpred_{\star}(\Theta w)$
		\item\label{Theta:ord} $b$ is an ordinal$^\star$ iff $\exists \alpha\ b = \Theta \alpha$, whereupon $\alpha = \conchrank{b}$
	\end{listn-0}
\end{prop}
\begin{proof}
	\emphref{Theta:Ulev}	
	If $\heblandpred^\star(b)$ then $b \in \Ulev{\conchrank{b}+1}{\emptyset}$ by induction on rank using Lemma \ref{lem:LTwand:transint}. For the converse, suppose $\Ufrom{\emptyset}(b)$. By construction, $\blandpred_\star(b)$ and $\uncarrier{b} \subseteq \Ulev{\beta}{\emptyset}$ for some $\beta$. Now $\carrier{\Ulev{\beta}{\emptyset}}$ witnesses that $\heblandpred^\star(b)$. 
	After all, $b \subseteq^\star \carrier{\Ulev{\beta}{\emptyset}}$ by Lemma \ref{lem:LTwand:transint}. And if 
	$x \in_\star \carrier{\Ulev{\beta}{\emptyset}}$ i.e.\ $x \in \Ulev{\beta}{\emptyset}$, then $\blandpred_\star(x)$ and $\uncarrier{x} \subseteq \Ulev{\beta}{\emptyset}$ so that $x \subseteq^\star \carrier{\Ulev{\beta}{\emptyset}}$ by Lemma \ref{lem:LTwand:transint}. 
%
	
	This establishes the biconditional; ``so that'' holds via Lemma \ref{lem:LTwand:Theta}.
	
	\emphreffrom{Theta:iso:in}\emphref{Theta:iso:wandpredlt} Immediate from \eqref{Theta:Ulev}. 


	\emphref{Theta:ord} By induction, using $\Theta$'s injectivity and \eqref{Theta:iso:in}, if $\alpha$ is any ordinal, then $\Theta \alpha$ is an ordinal$^\star$. Furthermore, $\alpha = \ordrank{\alpha} = \conchrank{\Theta \alpha}$ by Lemma \ref{lem:LTwand:Theta}. Now let $c$ be an ordinal$^\star$ and let $d = \Theta(\conchrank{c})$. So $\conchrank{c} = \conchrank{d}$ as above, so that $c \notin_\star d$ and $d \notin_\star c$ by Lemma \ref{lem:LTwand:transint}; so $c = d = \Theta(\conchrank{c})$. 
\end{proof}\noindent
It follows that hereditary blandness$^\star$ does not depend on the choice of $\star$. We also obtain the $\star$-translations of any \WS-axioms which only concern hereditarily bland sets:
\begin{lem}[\LTwand]\label{lem:LTwand:hevanstar} 	
	\ref{WS:wandsheb}$^\star$ holds, as does $(\phi^{\heblandpred})^\star$, for any \WS-axiom $\phi^{\heblandpred}$.
\end{lem}
\begin{proof}
	By Proposition \ref{prop:LTwand:iso}\eqref{Theta:Ulev}, every wand$_\star$ is hereditarily bland$^\star$.	 Next, by Definition \ref{def:LTwand}: whenever $\phi^{\heblandpred}$ is an \WS-axiom, $\phi$ is an \LTwand-axiom. Now note that, for any sentence $\phi$ in \LTwand's signature, $\phi \liff (\phi^{\heblandpred})^\star$; this holds by an induction on complexity using Proposition \ref{prop:LTwand:iso}. 
\end{proof}

\subsection{Wevels$^\star$ and levels$_\urbase^\star$}\label{s:tolt:wevlev}
I will now move on to considering ``wand taps'', starting with some easy facts:
\begin{lem}[\LTwand]\applabel{lem:LTwand:nonbland}
	If $\Delta_\star(c)$ and $c$ is an $\alpha$-tap, then:  
	\begin{listn-0}
		\item\label{n:tap:alpha+1} $\conchrank{c}=\alpha+1$, and $\conchrank{a} = \alpha$ for every $\tuple{w,a} \in c$; if also $\star = \paramint{\sigma}$ then $\alpha + 1 \leq \sigma$.
		\item\label{n:tap:ranklower} if $\taptostar{w}{a}{c}$, then $\conchrank{c} \leq \conchrank{a} + 1$.
		\item\label{n:tap:memtap} if $\tuple{w,a} \in c$, then $\taptostar{w}{a}{c}$.
	\end{listn-0}
\end{lem}\noindent
This allows us to describe the moment when a conch is ``first found'':
\begin{lem}[\LTwand]\applabel{lem:LTwand:helpercum}
	If $\Delta_\star(c)$ and $\Delta_\star(\carrier{\belowstage{\alpha}})$, then:
	$\conchrank{c} \leq \alpha$ iff $c \cfoundat^\star \carrier{\belowstage{\alpha}}$.
\end{lem}
\begin{proof} 
	If $c$ is bland$_\star$, use Lemma \ref{lem:LTwand:transint}.
	Otherwise, $c$ is some $\sigma$-tap by Lemma \ref{lem:LTwand:xor}; now use Lemmas \ref{lem:LTwand:transint} and \ref{lem:LTwand:nonbland}. 
\end{proof}\noindent 
The next challenge is characterize the wevels$^\star$. This is easier than we might fear. The basic facts about wevels follow from \CoreWev (see \S\ref{s:Wevels} and \S\ref{s:WS:wevels}); so the $\star$-translations of those basic facts hold by Lemmas \ref{lem:LTwand:note} and \ref{lem:LTwand:extsepstar}. Specifically: $\in_\star$ well-orders the wevels$^\star$ (by Theorem \ref{thm:WS:wo}$^\star$), and a wevel$^\star$ is exactly the bland$_\star$ set of everything found-earlier$^\star$ (by Lemma \ref{lem:WS:acc}$^\star$). Once we combine this with Lemma \ref{lem:LTwand:helpercum}, we can characterize the wevels$^\star$ in terms of the $\belowstage{\alpha}$'s:
 \begin{lem}[\LTwand]\label{lem:LTwand:wevelstar}
	If $\Delta_\star(\carrier{\belowstage{\alpha}})$, then $\carrier{\belowstage{\alpha}}$ is the $\alpha^\text{th}$ wevel$^\star$.
\end{lem}
\begin{proof}
	Whenever $\Delta_\star(c)$, Lemma \ref{lem:LTwand:helpercum} tells us that $\conchrank{c} < \alpha$ iff $(\exists \beta < \alpha)c \cfoundat^\star \carrier{\belowstage{\beta}}$. So suppose for induction that $\carrier{\belowstage{\beta}}$ is the $\beta^\text{th}$ wevel$^\star$, whenever $\beta < \alpha$. Then by Lemma \ref{lem:LTwand:transint}, the earlier biconditional becomes
	\begin{align*}
		c \in_\star \carrier{\belowstage{\alpha}}&\text{ iff }
		\relexists{r}{\wevpred^\star}c \cfoundat^\star r \in_\star \carrier{\belowstage{\alpha}}
		\end{align*}
	So $\carrier{\belowstage{\alpha}}$ is the $\alpha^\text{th}$ wevel$^\star$, by Lemma \ref{lem:WS:acc}$^\star$. Now use induction, Theorem \ref{thm:WS:wo}$^\star$.
%
\end{proof}\noindent
Combining Lemmas \ref{lem:LTwand:helpercum}--\ref{lem:LTwand:wevelstar} allows us to characterize $\wevof^\star$:
\begin{lem}[\LTwand]\label{lem:LTwand:wevelneat}
	If $\Delta_\star(c)$ and 
	$\Delta_\star(\carrier{\belowstage{\conchrank{c}}})$, then $\wevof^\star{(c)} =\carrier{\belowstage{\conchrank{c}}}$ and $\ordranktoint{(c)}{\star} = \Theta(\conchrank{c})$. 
\end{lem}
\begin{proof}
	In this case, $c \cfoundat^\star \carrier{\belowstage{\conchrank{c}}}$, and $\conchrank{c}$ is minimal here by Lemma \ref{lem:LTwand:helpercum}. By Lemma \ref{lem:LTwand:wevelstar}, $\carrier{\belowstage{\conchrank{c}}}$ is the $\conchrank{c}^\text{th}$ wevel$^\star$. By Proposition \ref{prop:LTwand:iso}\eqref{Theta:ord}, the $\conchrank{c}^\text{th}$ ordinal$^\star$ is $\Theta(\conchrank{c})$, i.e.\ $\ordranktoint{(c)}{\star} = \Theta(\conchrank{c})$. 
\end{proof}\noindent
This delivers the $\tolt$-translation of two of \WS's axioms:
\begin{lem}[\LTwand]\label{lem:LTwand:stratheight}
	\ref{WS:strat}$^\tolt$ and \ref{WS:weakheight}$^\tolt$ hold.
\end{lem}
\begin{proof}
	\emph{\ref{WS:strat}$^\tolt$.} By Lemma \ref{lem:LTwand:wevelneat}, as $\Delta_\tolt(\carrier{\belowstage{\conchrank{c}}})$ for any conch $c$. 
	
	\emph{\ref{WS:weakheight}$^\tolt$.} Using Pseudo-\ref{WS:weakheight} and Lemma \ref{lem:LTwand:Theta}, let $a = \Setabs{w}{\wandpredlt(w)}$, and let $b = \Theta a = (\Setabs{w}{\wandpred(w)})^\tolt$. 
	Now $\Theta(\ordrank{a}) = \Theta(\conchrank{b}) = \ordranktoint{b}{\tolt}$ by Lemmas \ref{lem:LTwand:Theta} and \ref{lem:LTwand:wevelneat}. 
	Also, $\ordrank{a} + \alpha$ exists for any ordinal $\alpha$, by Pseudo-\ref{WS:weakheight}. So $\Theta(\ordrank{a} + \alpha) = \Theta(\ordrank{a}) +^\tolt \Theta\alpha$ is an ordinal$^\tolt$ by Proposition \ref{prop:LTwand:iso}; so $\ordranktoint{b}{\tolt} +^\tolt x$ exists for any ordinal$^\tolt$ $x$ by Proposition \ref{prop:LTwand:iso}\eqref{Theta:ord}.
\end{proof}\noindent
Note that this only gives us the $\tolt$-translations of \ref{WS:strat} and \ref{WS:weakheight}. Indeed, we should \emph{not} expect their $\paramint{\sigma}$-translations to hold: denizens of $\Delta_{\paramint{\sigma}}$ must sit at some level$_{\carrier{\conchstage{\sigma}}}^\paramint{\sigma}$ but  may come after any wevel$^\paramint{\sigma}$. Our next job, then, must be to characterize the levels$_{\carrier{\urbase}}^\star$. This can be done along the same line as the wevels$^\star$, mimicking the proof of Lemma \ref{lem:LTwand:wevelstar}:
\begin{lem}[\LTwand]\applabel{lem:LTwand:levelnu}
	If $\Delta_\star(\carrier{\urbase})$, then  $\carrier{\Ulev{\alpha}{\urbase}}$ is the $\alpha^\text{th}$ level$_{\carrier{\urbase}}^\star$.
\end{lem}
\begin{proof}
	By Lemmas \ref{lem:LTwand:transint}, \ref{lem:LTwand:helpercum}, \ref{lem:C0:acc}$^\star$, and \ref{lem:C0:wo}$^\star$.
\end{proof}\noindent
Now that we understand wevels$^\star$ and levels$_{\carrier{\urbase}}^\star$, we can $\star$-translate the central notion used in the Picture of Loose Constructivism, i.e.\ $\Vfrom{\urbase}$:
\begin{lem}[\LTwand]\label{lem:LTwand:UV}
	If $\Delta_\star(\carrier{\belowstage{\alpha}})$, these are equivalent predicates:
	\begin{align*}
		\Ufrom{\belowstage{\alpha}}(x) && \exists w\big(w = 
		\carrier{\belowstage{\alpha}} \land (\Vfrom{w}(x))^\star\big)&&
		\exists w\big(\conchrank{w} = 
		\alpha \land (\Vfrom{\wevof{w}}(x))^\star\big)
	\end{align*}
\end{lem}
\begin{proof}
	Use Lemmas \ref{lem:LTwand:levelnu} and  
%
	\ref{lem:LTwand:wevelneat} respectively for the two equivalences.
\end{proof}\noindent
This delivers a key fact: satisfaction of loosely bound formulas is preserved as we expand our interpretations, moving (with $\sigma < \tau$) from $\paramint{\sigma}$ to $\paramint{\tau}$ and ultimately to $\tolt$:
\begin{prop}[\LTwand, scheme]\label{prop:LTwand:starlt}
	Let $\phi(v_1, \ldots, v_n)$ be loosely bound by $\max(\wevof{v_1}, \ldots, \wevof{v_n})$. When $\max(\conchrank{c_1}, \ldots, \conchrank{c_n}) \leq \sigma < \tau$, these are equivalent:
	\begin{align*}
		\phi^\paramint{\sigma}(c_1, \ldots, c_n) & & \phi^\paramint{\tau}(c_1, \ldots, c_n) & & \phi^\tolt(c_1, \ldots, c_n)
	\end{align*}
\end{prop} 
\begin{proof}
	Let $n = 1$, and consider $v = v_1$ and $c = c_1$ with $\conchrank{c} = \gamma \leq \sigma$. (This is for readability; the same proof works for $n > 1$.) 
	Since $\phi(v)$ is loosely bound by $\wevof{v}$:
	\begin{listn-0}
		\item\label{restrict:V} $\phi(v)$'s free and bound variables are restricted to $\Vfrom{(\wevof{v})^+}$
		\item\label{restrict:conjoin} any instance of $\tapto{w}{a}{x}$ in $\phi(v)$ is conjoined to $a \in \wevof{v}$.
	\end{listn-0}
	By \eqref{restrict:V},  $\phi^{\paramint{\sigma}}(v)$'s free and bound variables are restricted to $(\Vfrom{(\wevof{v})^+})^\paramint{\sigma}$. Since $\gamma\leq\sigma$, we have $\Delta_\paramint{\sigma}(\carrier{\conchstage{\gamma}})$; 
	so, substituting $c$ for $v$, by Lemma \ref{lem:LTwand:UV},  $\phi^{\paramint{\sigma}}(c)$'s free and bound variables are, equivalently, restricted to $\Ufrom{\conchstage{\gamma}}$. Now, $\paramint{\sigma}$-translation also restricts all quantifiers to $\Ufrom{\conchstage{\sigma}}$; but since $\gamma \leq \sigma$, this imposes no \emph{additional} restriction. So we have (up to equivalence):
	\begin{listn} 
		\item[(\ref{restrict:V}$c$)] $\phi^\paramint{\sigma}(c)$'s free and bound variables are restricted to $\Ufrom{\conchstage{\gamma}}$
	\end{listn}
	Similar reasoning using Lemma \ref{lem:LTwand:wevelneat} gives us (up to equivalence):
	\begin{listn}
		\item[(\ref{restrict:conjoin}$c$)] any instance of $\taptosigma{w}{a}{x}$ in $\phi^\paramint{\sigma}(c)$ is conjoined to $\conchrank{a} < \gamma$.
	\end{listn}
%
	Finally, recall that $\paramint{\sigma}$, $\paramint{\tau}$ and $\tolt$ explicitly agree on everything, except on any mentions of $\Delta_\star$ and on their rendering of wand-taps. But the differences concerning $\Delta_\star$ are irrelevant here, since $\gamma \leq \sigma < \tau$, so that $\paramint{\tau}$ and $\tolt$ yield the same  restrictions (\ref{restrict:V}$c$)--(\ref{restrict:conjoin}$c$). Then Lemma \ref{lem:LTwand:nonbland} and (\ref{restrict:conjoin}$c$) ensure they agree about wand-taps. 
	So they fully agree. 
\end{proof}
\begin{cor}[\LTwand]\label{cor:LTwand:preserve}
	Lemma \ref{lem:WS:wonderful}$^\tolt$ holds, so that e.g.\ $\equivpred^\tolt$ is an equivalence relation. Furthermore, these formulas are equivalent (row-wise), when $\conchrank{a}, \conchrank{b} \leq \sigma$:
	\begin{align*}
		\wanddomtolt{w}{a}&&
		\wanddomint{w}{a}{\paramint{\sigma}} &&
		\tuple{w,a} \in \domstage{\conchrank{a}}\\
		\churchequivtolt{w}{a}{u}{b}&&
		\churchequivint{w}{a}{u}{b}{\paramint{\sigma}} &&
		\tuple{w,a,u,b} \in \equivstage{\max(\conchrank{a}{, \conchrank{b}})}
	\end{align*}
\end{cor}
\begin{proof}
	Lemma \ref{lem:WS:wonderful}$^\tolt$ holds by Lemmas \ref{lem:LTwand:extsepstar} and \ref{lem:LTwand:wevelneat}. 
	This gives us \domref{dom:wand}$^\tolt$--\domref{dom:default}$^\tolt$ and \equivref{equiv:wand}$^\tolt$--\equivref{equiv:default}$^\tolt$. 
	Now Proposition \ref{prop:LTwand:starlt} yields the equivalences. 
%
\end{proof}\noindent
Using Corollary \ref{cor:LTwand:preserve} freely and frequently from this point onwards, we are at last in a position to show that $\tolt$ is an interpretation:
\begin{lem}[\LTwand]\applabel{lem:LTwand:lastbits}
	\ref{WS:funtap}$^\tolt$, \ref{WS:make}$^\tolt$, and \ref{WS:eq}$^\tolt$ hold
\end{lem}
\begin{proof}
	\emph{\ref{WS:funtap}$^\tolt$.} Suppose $\taptolt{w}{a}{c}$ and $\taptolt{w}{a}{d}$. So there is $\tuple{w_1,a_1} \in c$ with $\churchequivtolt{w}{a}{w_1}{a_1}$ and $\tuple{w_2, a_2} \in d$ with $\churchequivtolt{w}{a}{w_2}{a_2}$. 
	Now  $\churchequivtolt{w_1}{a_1}{w_2}{a_2}$, 
	so $\conchrank{a_2} = \conchrank{a_1}$ by \eqref{conch:c:rank}. Hence $\tuple{w_1, a_1} \in d$ and $\tuple{w_2, a_2} \in c$ 
	by \eqref{conch:c:equiv}. Since $c$ and $d$ are overlapping equivalence classes of minimal rank, it is easy to show that $c = d$. 
	
	\emph{\ref{WS:make}$^\tolt$.} Similar but easier.

	\emph{\ref{WS:eq}$^\tolt$.} Similar, using \ref{WS:make}$^\tolt$.
\end{proof}\noindent
Assembling all the $\tolt$-translations of \WS's axioms, we have what we required:
\begin{prop}\label{prop:LTwand:Ctolt}
	$\tolt: \WS \functionto \LTwand$ is an interpretation.
\end{prop}
\noindent
We can also obtain the Corollary which I stated as a target in \S\ref{s:tolt:def}, via an easy lemma:
\begin{lem}[\LTwand]\label{lem:LTwand:teeny}
	If $\taptolt{w}{a}{c}$, then $c = \Setabs{\tuple{u,b}}{\churchequivtolt{w}{a}{u}{b} \land \minirankint{u}{b}{\tolt}}$.
\end{lem}
\begin{proof}
	Suppose $\taptolt{w}{a}{c}$ and let 
	\begin{align*}
		e &= \Setabs{\tuple{u,b}}{\churchequivtolt{w}{a}{u}{b} \land \minirankint{u}{b}{\tolt}}
	\end{align*}
	To see $c \subseteq e$, fix $\tuple{u,b} \in c$. Now $\taptolt{u}{b}{c}$ by Lemma \ref{lem:LTwand:nonbland}\eqref{n:tap:memtap}, and so $\churchequivtolt{w}{a}{u}{b}$ by \ref{WS:eq}$^\tolt$. Also, if $\churchequivtolt{u}{b}{v}{c}$, then 
	$\conchrank{b} \leq \conchrank{c}$ by \eqref{conch:c:rank} i.e.\ $\ordranktoint{b}{\tolt} \leq^\tolt \ordranktoint{c}{\tolt}$ by Lemma \ref{lem:LTwand:wevelneat}, so $\minirankint{u}{b}{\tolt}$. To see $e \subseteq c$, keep $\tuple{u,b} \in c \neq \emptyset$ as before, and suppose $\tuple{u', b'} \in e$. Now $\churchequivtolt{u}{b}{u'}{b'}$, and since $\minirankint{u'}{b'}{\tolt}$ we have $\tuple{u',b'} \in c$ by \eqref{conch:c:equiv}. 
\end{proof}
\begin{cor}[\LTwand]\label{cor:LTS:target} $c$ is a conch iff either: 
	\begin{listn-0}
		\item\label{conch:taxonomy:bland} $c = \carrier{s}$, for some $s$ whose members are all conches; or
		\item\label{conch:taxonomy:nonbland} $c = \Setabs{\tuple{w,a}}{\taptolt{w}{a}{c}\land \miniranktolt{w}{a}} \neq \emptyset$
	\end{listn-0}%
\end{cor}

\section{Bi-interpretation and synonymy}\label{s:int:biint}
In the past two sections, I have laid down two interpretations
	$$\toc : \LTwand \mathrel{\stackanchor[-2pt]{\functionto}{\longleftarrow}} \WS : \tolt$$
Furthermore, Proposition \ref{prop:LTwand:iso} tells us that $\toc\tolt$ is a self-embedding. I now want to show that $\tolt\toc$ is too, so that $\toc$ and $\tolt$ form a bi-interpretation. 

Inspired by Corollary \ref{cor:LTS:target}, I will recursively define the map, $\Xi$, which will witness that $\tolt\toc$ is a self-embedding:
\begin{define}[\WS]\label{def:WS:Xi}
	For each $c$, define $\Xi c$ as a bland set of ordered pairs thus:
	\begin{align*}
		\Xi c &\coloneq
		\carrier{\Setabs{\Xi x}{x \in c}}&&\text{if }\blandpred(c)\\
		\Xi c &\coloneq \Setabs{\tuple{\Xi w, \Xi a}}{\tapto{w}{a}{c} \land \minirank{w}{a}}&&\text{if }\lnot\blandpred(c)
	\end{align*}
	where $\carrier{s} = \{\tuple{\emptyset, s}\}$, as before. \defendhere
\end{define}\noindent
This recursive definition has a twist. We \emph{first} define $\Xi c$, just for hereditarily bland $c$; this ensures that $\Xi w$ is well-defined for each wand $w$. \emph{Then} we augment this to a definition which covers everything else. This second step does not disrupt the definition for hereditarily bland entities, so it is well-defined. We can immediately establish some useful properties (using induction and Lemma \ref{lem:WS:hered}):
\begin{lem}[\WS]\label{lem:WS:Xiinj}
	$\Xi$ is an injective map $V \functionto \heblandpred$.
\end{lem}\noindent
My goal, though, is to show that $\Xi$ is an isomorphism $V \functionto \conchpred^\toc$. My strategy is to use ordinal induction within \WS to show that, for each ordinal $\alpha$, $\Xi$'s restriction is a structure-preserving bijection from $\Vfrom{\wevof{\alpha}}$ to $(\Vfrom{\wevof{\alpha}})^{\tolt\toc}$. 

Let me start by unpacking the last expression. By Proposition \ref{prop:LTwand:iso}\eqref{Theta:ord}$^\toc$ and Lemma \ref{lem:LTwand:UV}$^\toc$, we know that $(\Vfrom{\wevof{\alpha}}(x))^{\tolt\toc}$ 
is equivalent to $(\Ufrom{\belowstage{\alpha}}(x))^\toc$. Next, note that $\toc$ translates all predicates verbatim, but restricts our purview to $\heblandpred$; but $\Ufrom{\belowstage{\alpha}}(x)$, as defined in \WS, is clearly solely concerned with hereditarily bland entities in any case; so we can just \emph{ignore} $\toc$'s action.\footnote{Strictly speaking, there is a \emph{tiny} abuse of notation here: to characterize $\emptyset$ in \WS, we say ``the \emph{bland} thing with no members'', rather than ``the thing with no members'' (as in \LTwand). So, strictly speaking, a slightly different definition of ``$\Ufrom{\belowstage{\alpha}}$'' is required. But the difference would evidently make no real difference. I repeatedly rely on this harmless abuse of notation throughout this section.} My aim, then, is to show that $\Xi$'s restriction is a structure-preserving bijection $\Vfrom{\wevof{\alpha}} \functionto \Ufrom{\belowstage{\alpha}}$. 

In fact, it is easy to show (without using induction) that $\Xi$ preserves $\wandpred$:
\begin{lem}[\WS]\label{lem:IH:wandpreserve}
	$\wandpred^{\tolt\toc}(x)\lonlyif \relexists{w}{\wandpred}x = \Xi w$, and $\wandpred(w) \liff \wandpred^{\tolt\toc}(\Xi w)$
\end{lem}
\begin{proof} 
	Using \ref{WS:wandsheb}, note that  
	$\wandpred^{\tolt\toc}(x)$ 
	iff
	$\relexists{v}{\wandpred}x = \Theta^\toc v \land \Ufrom{\emptyset}(x) \land \conchpred^\toc(x)$. Now $\Xi$, restricted to hereditarily bland things, is $\Theta^\toc$. This gives the first claim. For the second:  $\Ufrom{\emptyset}(x)$ entails $\conchpred^\toc(x)$; and $\Xi$ is injective; so $\wandpred^{\tolt\toc}(\Xi w)$ iff $\wandpred(w) \land \Ufrom{\emptyset}(\Xi w)$. To complete the proof, just note that when $\wandpred(w)$, so that $\heblandpred(w)$, we have $\Ufrom{\emptyset}(\Xi w)$.
\end{proof}\noindent
However, it will take more effort to secure the preservation of \LTwand's other primitives. Indeed, it requires some thought about what we \emph{want} to preserve. I will set up my induction hypothesis to (try to) preserve $\blandpred$ and $\in$ in the most obvious way. However, recalling the Picture of Loose Constructivism from \S\ref{s:Constructive}, when considering entities of rank $< \alpha$, we must remember that wand-tapping an entity whose rank is \emph{immediately} below rank $\alpha$ will take us \emph{outside} of $\Vfrom{\wevof{\alpha}}$. For this reason, my induction hypothesis considers preservation of $\tappred$ for entities found at least one clear rank before $\alpha$. 

With all that in mind, here is my big induction hypothesis, BIH, for all $\beta < \alpha$:
	\begin{align*}
		\text{If }\Vfrom{\wevof{\beta}}(a)\text{ then:}&&\Ufrom{\belowstage{\beta}}(\Xi a)&\text{ and }\ordrank{a} = \naughtyrank{\Xi a}\\
		\text{If }\Ufrom{\belowstage{\beta}}(b)\text{ then:}&&b = \Xi a&\text{ for some $a$ such that $\Vfrom{\wevof{\beta}}(a)$}\\
		\text{if }\Vfrom{\wevof{\beta}}(a)\text{ then:}&&	\blandpred(a) & \liff \blandpred^{\tolt\toc}(\Xi a)\\
		\text{If }\Vfrom{\wevof{\beta}}(a)\text{ and }\Vfrom{\wevof{\beta}}(b)\text{ then:}&&	a \in b & \liff \Xi a \in^{\tolt\toc} \Xi b\\
		 \text{If }\ordrank{a} + 1 < \beta\text{ then:}	&&\tapto{w}{a}{c} & \liff \taptoint{\Xi w}{\Xi a}{\Xi c}{\tolt\toc}
	\end{align*}
\emph{\textbf{Note:} from here until the end of the section, I reserve $\alpha$ and $\beta$ for use in this induction hypothesis, BIH.} Assuming BIH, I will show that each of the claims holds with $\alpha$ in place of $\beta$. By ordinal induction, the claims will hold for every ordinal, so that $\Xi : V \functionto \conchpred^\toc$ is isomorphism, as required. 

My first result ensures preservation of loosely bound formulas:
\begin{lem}[\WS, BIH]
	\label{lem:WS:ind:preserves}
	Let $\phi(v_1, \ldots, v_n)$ be loosely bound by $\max(\wevof{v_1}, \ldots, \wevof{v_n})$. Let $\max(\ordrank{a_1}, \ldots, \ordrank{a_n}) < \beta$. Then $\phi(a_1, \ldots, a_n) \liff \phi^{\tolt\toc}(\Xi a_1, \ldots, \Xi a_n)$. So in particular, this holds for $\wanddompred$, $\equivpred$ and $\minirankpred$.
\end{lem}
\begin{proof}
	Where $\gamma = \max(\ordrank{a_1}, \ldots, \ordrank{a_n}) < \beta$, bound and free variables in $\phi(a_1, \ldots, a_n)$ are restricted to $\Vfrom{(\wevof{\gamma})^+}$, and any instance of ``$\tapto{x}{y}{z}$'' is conjoined with ``$x \in \wevof{\gamma}$'', i.e.\ in effect the claim that $\ordrank{x} < \gamma$. By BIH, $\Xi$'s restriction is a structure-preserving bijection $\Vfrom{\wevof{\beta}} \functionto \Ufrom{\belowstage{\beta}}$, with the caveat that $\tappred$ is preserved when tapping entities of rank $< \gamma < \beta$, which aligns perfectly with the conjunctive-restriction. So the biconditional holds by induction on complexity.
\end{proof}\noindent
Using this, I can show that $\tappred$ is preserved in the required way:
\begin{lem}[\WS, BIH]\phantom{.}
	\label{lem:WS:ind:taplower} 
	\begin{listn-0}
		\item\label{n:IH:tapup} If $\ordrank{a} + 1 < \alpha$ and $\tapto{w}{a}{c}$, then $\taptoint{\Xi w}{\Xi a}{\Xi c}{\tolt\toc}$ and $\ordrank{c} = \naughtyrank{\Xi c}$.
		\item\label{n:IH:nonbland} If $\naughtyrank{e} < \alpha$ and $\lnot\blandpred^{\tolt\toc}(e)$, then $e = \Xi c$ for some non-bland $c$, and $\ordrank{c} = \naughtyrank{\Xi c}$.
		\item\label{n:IH:tapdown} If $\ordrank{a} + 1 < \alpha$ and $\taptoint{\Xi w}{\Xi a}{\Xi c}{\tolt\toc}$, then $\tapto{w}{a}{c}$.
	\end{listn-0}
\end{lem}
\begin{proof}
	\emphref{n:IH:tapup} 
	Let $\ordrank{a} + 1 < \alpha$ and $\tapto{w}{a}{c}$. 
	So $\wanddom{w}{a}$ by \ref{WS:make}, and 
	$\wanddomint{\Xi w}{\Xi a}{\tolt\toc}$ by Lemma \ref{lem:WS:ind:preserves}. 
	By \ref{WS:make}$^{\tolt\toc}$, there is $e$ such that $\taptoint{\Xi w}{\Xi a}{e}{\tolt\toc}$. By Lemma \ref{lem:LTwand:teeny}$^\toc$:
	\begin{align*}
		e &= \Setabs{\tuple{x,d}}{\churchequivint{\Xi w}{\Xi a}{x}{d}{\tolt\toc} \land \minirankint{x}{d}{\tolt\toc}}
	\intertext{Fix $\tuple{x,d} \in e$; so $\naughtyrank{d} \leq \naughtyrank{\Xi a}$; but $\naughtyrank{\Xi a} = \ordrank{a}$ by BIH; so $\naughtyrank{d} +1 < \alpha$.  By BIH again, $d = \Xi b$ for some $b$ with $\ordrank{b} = \naughtyrank{d}$. Also, $x = \Xi u$ for some $u$ by Lemma \ref{lem:IH:wandpreserve}. So:}
		e &= \Setabs{\tuple{\Xi u,\Xi b}}{\churchequivint{\Xi w}{\Xi a}{\Xi u}{\Xi b}{\tolt\toc} \land \minirankint{\Xi u}{\Xi b}{\tolt\toc}}\\
	\intertext{and therefore, by Lemma \ref{lem:WS:ind:preserves} and \ref{WS:eq}:}
		e &= \Setabs{\tuple{\Xi u,\Xi b}}{\tapto{u}{b}{c} \land \minirank{u}{b}} = \Xi c
	\end{align*}
	To compute rank, fix $\tuple{\Xi u, \Xi b} \in e$. Since $\minirank{u}{b}$, both $\ordrank{b} + 1 = \ordrank{c}$ by Lemma \ref{lem:WS:minnonbland}, and also $\ordrank{b} \leq \ordrank{a}$ and hence $\ordrank{b} = \naughtyrank{\Xi b}$ by BIH. Since $\naughtyrank{\Xi b} + 1= \naughtyrank{e}$, we get $\ordrank{c} = \naughtyrank{e}$.

	\emphref{n:IH:nonbland}
	With $\conchrank{e} < \alpha$ and $e$ non-bland$^{\tolt\toc}$, Corollary \ref{cor:LTS:target}$^\toc$ tells us:
	\begin{align*}
		e &= \Setabs{\tuple{x,d}}{\taptoint{x}{d}{e}{\toc\tolt}\land \minirankint{x}{d}{\tolt\toc}} \neq \emptyset 
	\end{align*}
	Fix $\tuple{x,d} \in e$. So $\taptoint{x}{d}{e}{\tolt\toc}$ and $\naughtyrank{d} < \naughtyrank{e} < \alpha$ so $\wanddomint{x}{d}{\tolt\toc}$ by \ref{WS:make}$^{\tolt\toc}$. By BIH, $d = \Xi a$ for some $a$ with $\ordrank{a} = \naughtyrank{d}$. Using Lemma \ref{lem:IH:wandpreserve}, let $x = \Xi w$; now $\wanddom{w}{a}$ and $\minirank{w}{a}$ by Lemma \ref{lem:WS:ind:preserves}; so by \ref{WS:make} there is $c$ such that $\tapto{w}{a}{c}$. Since $\ordrank{a} + 1 < \alpha$,  now $\ordrank{c} = \naughtyrank{\Xi c}$ and $\taptoint{\Xi w}{\Xi a}{\Xi c}{\tolt\toc}$ by \eqref{n:IH:tapup}, i.e.\ $\taptoint{x}{d}{\Xi c}{\tolt\toc}$; so $e = \Xi c$ by \ref{WS:funtap}$^{\tolt\toc}$. 
	
	\emphref{n:IH:tapdown}
	Let $\ordrank{a} + 1 < \alpha$ and $\taptoint{\Xi w}{\Xi a}{\Xi c}{\tolt\toc}$. 
	Now $\naughtyrank{\Xi c} \leq \naughtyrank{\Xi a} +1$ by Lemma \ref{lem:LTwand:nonbland}\eqref{n:tap:ranklower}$^\toc$. Also $\ordrank{a} = \naughtyrank{\Xi a}$ by BIH; so $\naughtyrank{\Xi c} < \alpha$ and hence $c$ is not bland and $\ordrank{c} = \naughtyrank{\Xi c} < \alpha$ by \eqref{n:IH:nonbland}. Using Lemma \ref{lem:WS:min}, fix $u,b$ with $\tapto{u}{b}{c}$ and $\ordrank{b} + 1 = \ordrank{c} < \alpha$; now $\taptoint{\Xi u}{\Xi b}{\Xi c}{\tolt\toc}$ by \eqref{n:IH:tapup}, so $\churchequivint{\Xi w}{\Xi a}{\Xi u}{\Xi b}{\tolt\toc}$ by \ref{WS:eq}$^{\tolt\toc}$, and $\tapto{w}{a}{c}$ by Lemma \ref{lem:WS:ind:preserves} and \ref{WS:eq}.
\end{proof}\noindent
With wand-taps taken care of, I can now show that $\Xi$'s restriction is the bijection we wanted, and that it preserves blandness as hoped:
\begin{lem}[\WS, BIH]\label{lem:IH:blandrank}
	For any $\gamma$:
	\begin{listn-0}
		\item\label{n:IH:LtoU} If $a \in \Tlev{\gamma}{\wevof{\alpha}}$, then $\ordrank{a} = \conchrank{\Xi a}$ and $\Xi a \in \Ulev{\gamma}{\belowstage{\alpha}}$ and $\blandpred(a) \liff \blandpred^{\tolt\toc}(\Xi a)$.
		\item\label{n:IH:UtoL} If $b \in \Ulev{\gamma}{\belowstage{\alpha}}$, then $b = \Xi a$ for some $a \in \Tlev{\gamma}{\wevof{\alpha}}$.
	\end{listn-0}
	Hence $\Xi$'s restriction is a bijection $\Vfrom{\wevof{\alpha}} \functionto \Ufrom{\belowstage{\alpha}}$, and if $\Vfrom{\wevof{\alpha}}(a)$ then $\blandpred(a) \liff \blandpred^{\tolt\toc}(\Xi a)$.
\end{lem}
\begin{proof}
	Both claims hold by induction on $\gamma$; ``hence'' then follows from $\Xi$'s injectivity and the definitions of $\Vfrom{\wevof{\alpha}}$ and $\Ufrom{\belowstage{\alpha}}$. In what follows, I use the fact that  $\blandpred^{\tolt\toc}(v)$
	iff
%
%
	$\relexists{e}{\heblandpred}v = \carrier{e} \land \conchpred^\toc(v)$. 
	
	\emph{Induction case when $\gamma = 0$.} For \eqref{n:IH:LtoU}, suppose $a \in \Tlev{0}{\wevof{\alpha}} =  \wevof{\alpha}$, i.e.\ $\ordrank{a} < \alpha$. First, suppose $a$ is not bland. Using Lemma \ref{lem:WS:min}, fix $w$ and $b$ such that $\tapto{w}{b}{a}$ and $\ordrank{b} < \ordrank{a}$; now  $\naughtyrank{\Xi a} = \ordrank{a}$ and $\lnot\blandpred^{\tolt\toc}(\Xi a)$ by Lemma \ref{lem:WS:ind:taplower}\eqref{n:IH:tapup}. 
	
	Next, suppose $a$ is bland. If $x \in a$, then $\ordrank{x} < \ordrank{a} < \alpha$, so $\conchpred^\toc(\Xi x)$ by BIH. Let $e = \Setabs{\Xi x}{x \in a}$; this witnesses that $\blandpred^{\tolt\toc}(\Xi a)$, because $\heblandpred(e)$ by Lemma \ref{lem:WS:Xiinj}, and $\conchpred^\toc(\carrier{e})$, and $\carrier{e} = \Xi a$. To compute $\naughtyrank{\Xi a}$, using BIH and Lemma \ref{lem:LTwand:transint}\eqref{n:rank}$^\toc$: 
		$$\ordrank{a} = \lsub_{x \in a}\ordrank{x} = \lsub_{x \in a} \naughtyrank{\Xi x} = \naughtyrank{\Xi a}$$

	For \eqref{n:IH:UtoL}, suppose $b \in \Ulev{0}{\belowstage{\alpha}} = \belowstage{\alpha}$, i.e.\ $\naughtyrank{b} < \alpha$. If $\lnot\blandpred^{\tolt\toc}(b)$, use Lemma \ref{lem:WS:ind:taplower}\eqref{n:IH:nonbland}; so suppose instead that $\blandpred^{\tolt\toc}(b)$; let $e$ witness this, so that $b = \carrier{e}$. If $z \in e$ then $\naughtyrank{z} < \naughtyrank{b}$ by Lemma \ref{lem:LTwand:transint}\eqref{n:rank}$^\toc$, so that $\ordrank{(\Xi^{-1}z)} = \naughtyrank{z}$ by BIH. Let $a = \Setabs{\Xi^{-1} z}{z \in e}$; now $b = \carrier{e} = \Xi a$,
	and indeed $\ordrank{a} = \naughtyrank{\Xi a}$ by the same computation as before. 
	
	\emph{Induction case when $\gamma > 0$.} Similar, using the induction hypothesis for $\gamma$ rather than BIH. 
\end{proof}\noindent
It follows straightforwardly that membership is preserved:
\begin{lem}[\WS, BIH]
	If $\Vfrom{\wevof{\alpha}}(a)$ and $\Vfrom{\wevof{\alpha}}(b)$, then $a \in b \liff \Xi a \in^{\tolt\toc} \Xi b$.
\end{lem}
\begin{proof}
	Note that 
	$\Xi a \in^{\tolt\toc} \Xi b$ 
	iff 	
	$\blandpred^{\tolt\toc}(\Xi b) \land \Xi a \in \uncarrier{\Xi b}$. 
	Now use $\Xi$'s injectivity and Lemma \ref{lem:IH:blandrank}. 
\end{proof}\noindent
And this was the last piece of the puzzle. We have proved the induction step; now by ordinal induction and Lemma \ref{lem:WS:Xiinj}, we immediately obtain what we wanted:
\begin{prop}[\WS]\label{prop:WS:embed}
	$\Xi: V \functionto \conchpred^{\toc}$ witnesses that the map $\tolt\toc$ is a self-embedding.
\end{prop}\noindent 
Assembling all of this, we obtain the Main Theorem, that \LTwand and \WS are synonymous.
\begin{proof}[\textbf{Proof of Main Theorem}]
	Propositions \ref{lem:WS:helowinterpret} and \ref{prop:LTwand:Ctolt} tell us that we have interpretations $\toc : \LTwand \mathrel{\stackanchor[-2pt]{\functionto}{\longleftarrow}} \WS : \tolt$. By Propositions \ref{prop:LTwand:iso} 
	and \ref{prop:WS:embed}, these witness an identity-preserving bi-interpretation. Since \LTwand is sequential, \LTwand and \WS are synonymous by the Friedman--Visser Theorem.\footnote{i.e.\ \textcite[Corollary 5.5]{FriedmanVisser:WBIS}.} 
\end{proof}

\section{Church's \cite*{Church:STUS} theory as a wand/set theory}\label{s:app:CUS}
My Main Theorem deals with wand/set theories in full generality. I will close with something much more specific. In \S\ref{s:Motivations:Church} and \S\ref{s:WS:definition}, I claimed that Church's \cite*{Church:STUS} set theory with a universal set can be presented as a wand/set theory (so that it is synonymous with a \ZF-like theory). The point of this section is to make good on this claim. 

Once you have defined the notion of a wevel, Church's \cite*{Church:STUS} theory can be given a perspicuous axiomatization; I cover this in \S\S\ref{s:app:CUS:axioms}--\ref{s:app:CCC:example1.3}. However, it is not immediately obvious that the perspicuous axiomatization is a wand/set theory, in the sense of Definition \ref{def:WandSet}; I show that it is in \S\ref{s:app:CUS:wandset}. 

\subsection{Axiomatizing Church's \cite*{Church:STUS} theory as \CUS}\label{s:app:CUS:axioms}
In this sub-section, I define \CUS, my version of Church's Universal Set theory. The basic idea is to help myself to the notion of a wevel, and then write down the axioms that would be obviously suggested by the discussion in \S\ref{s:Motivations:Church}.

\CUS's signature is the same as any wand/set theory. It has the quasi-notational axioms \ref{WS:notein}, \ref{WS:notetap}, and \ref{WS:funtap}. It also has the axioms which arrange everything into well-ordered wevels---\ref{WS:ext}, \ref{WS:sep}, and \ref{WS:strat}---and the axiom \ref{WS:weakheight}. 

As indicated in \S\ref{s:Motivations:Church}, I want to treat $0$ as the \textsf{complement} wand, and each $n > 0$ as the \textsf{cardinal}$_n$ wand. 
So Church's wands will be the natural numbers:
\begin{listaxiom}
	\labitem{Wands$_\textsc{c}$}{CUS:omega} $\Setabs{w}{\wandpred(w)} = \omega$
\end{listaxiom}
Note that this entails \ref{WS:wandsheb}, and it is formulated solely in terms of hereditarily bland sets and $\wandpred$ (see \S\ref{s:WS:heredbland}). So far, then, \CUS is on track to be a wand/set theory. 

I now need to provide axioms which specify how the wands behave. To explain the $\textsf{cardinal}_n$ wands, I need to introduce Church's notion of $n$-equivalence (though see \emph{Caveat 5} of \S\ref{s:app:CUS:caveats}):\footnote{\textcite[302]{Church:STUS}. Sheridan provides a version of Church's set theory which invokes a slight variant of $\approx_n$ (see \cite[109--11]{Sheridan:VCST}).} 
\begin{define}[\CoreWev, \ref{WS:strat}, \ref{WS:weakheight}]\label{def:cardinaln}
	Recursively define $\bigcup^0 a = a$ and $\bigcup^{n+1} a = \bigcup \bigcup^n a$.\footnote{This  is well-defined, for the reasons sketched in \S\ref{s:tolt:recursion}; and it is total since the axioms mentioned prove that $\bigcup x$ exists for any $x$.} For any $n$, say that $a$ and $b$ are $n$-equivalent, written $a \approx_{n} b$, iff: 
	\begin{listn-0} 
		\item both $\bigcup^i a$ and $\bigcup^i b$ are non-empty, and all their members are bland, for each $0 \leq i < n$; and
		\item there is a bijection $f_{n-1} : \bigcup^{n-1} a \functionto \bigcup^{n-1} b$; and
		\item setting $f_{i}(x) \coloneq \Setabs{f_{i+1}(y)}{y \in x}$ defines a bijection $f_i : \bigcup^i a \functionto \bigcup^i b$, for each $0 \leq i < n-1$.\footnote{Talk of bijections here is understood using the usual Kuratowski implementation for (bland) sets.} \defendhere
	\end{listn-0}
\end{define}\noindent
Note that $a \approx_1 b$ iff $a$ and $b$ are (non-empty, bland) sets which are equinumerous in the usual sense.

I now need axioms to say exactly when $\tapop{n}{a}$ should exist, for any wand $n$ and any object $a$. We want to say: tapping anything with \textsf{complement} yields an object; tapping anything bland with any \textsf{cardinal}$_n$ yields an object (if the $n$-equivalence class would be non-empty); and tapping anything non-bland with any \textsf{cardinal}$_n$ yields nothing. But there is a wrinkle. Suppose $b$ is bland; since we are treating $0$ as \textsf{complement}, we expect that $\tapop{0}{\tapop{0}{b}} = b$;\footnote{Note: $\tapop{w}{\tapop{u}{a}}$ is always to be read as $\tapop{w}(\tapop{u}a)$.} but this would violate \ref{WS:notetap}. 

To retain \ref{WS:notetap}, we will say that we \emph{find} an object by tapping $a$ with \textsf{complement} iff $a$ is not itself the complement of anything bland. (As foreshadowed in \S\ref{s:Notation}, I will eventually define a more expansive notion, according to which tapping $\tapop{0}{b}$ with $0$ will \emph{yield} $b$ itself; see Definition \ref{def:CUS:widetap}.) We effect this by stipulating:
\begin{listaxiom}
	\labitem{Making$_\textsc{c}$}{CUS:make} $\exists c\ \tapto{n}{a}{c} \liff \big((n = 0 \land  \relforall{b}{\blandpred}a \neq \tapop{0}{b}) \lor {}\\
	\phantom{\exists c\ \tapto{n}{a}{c} \liff \big(} (0 < n < \omega \land a \approx_n a)\big)$
\end{listaxiom} 
Note that this does not have the same shape as \ref{WS:make}, so we have departed from the shape of axioms stated in Definition \ref{def:WandSet}. (I return to this in \S\ref{s:app:CUS:wandset}.)

Our next axioms govern the \textsf{complement} wand:
\begin{listaxiom}
	\labitem{Comp1}{CUS:0:inj} if $\tapop{0}{a} = \tapop{0}{b}$, then $a = b$
	\labitem{Comp2}{CUS:0} if $\lnot\blandpred(a) \land \relforall{b}{\blandpred}a \neq \tapop{0}{b}$, then $a = \tapop{0}{\tapop{0}{a}}$
\end{listaxiom}
These say that $\tapop{0}{}$ is injective, and that  double-$\tapop{0}{\tapop{0}}$ is identity (provided this does not require taking $\tapop{0}{\tapop{0}{b}}$ for any bland $b$). Last, we have axioms governing the \textsf{cardinal}$_n$ wands:
\begin{listaxiom}
	\labitem{}{card:blank} if $0 < m < \omega$ and $0 <  n < \omega$, then:
	\labitem{Card1}{CUS:n} \hspace{2em}$\tapop{m}{a} = \tapop{n}{b} \liff (m=n \land a \approx_m b)$, and
	\labitem{Card2}{CUS:n:0b} \hspace{2em}$\relforall{b}{\blandpred}\tapop{m}{a} \neq \tapop{0}{b}$, and
	\labitem{Card3}{CUS:n:0ub}\hspace{2em}$\tapop{m}{a} \neq \tapop{0}{\tapop{n}{b}}$
\end{listaxiom}
So: $a$ and $b$ have the same cardinality$_n$ iff $a$ and $b$ are $n$-equivalent; if $n \neq m$ then no cardinal$_n$ is a cardinal$_m$; and no cardinal$_n$ is the complement of anything bland or of any cardinal$_m$.

We now define \CUS as the theory whose axioms are: \ref{WS:notein}, \ref{WS:notetap}, \ref{WS:funtap}, \ref{WS:ext}, \ref{WS:sep}, \ref{WS:strat}, \ref{WS:weakheight}, \ref{CUS:omega}, \ref{CUS:make}, \ref{CUS:0:inj}--\ref{CUS:0}, and \ref{CUS:n}--\ref{CUS:n:0ub}.

\subsection{Caveats and more expansive notions}\label{s:app:CUS:caveats}
With \CUS, we have a cleanly presented version of Church's \cite*{Church:STUS} theory. Or so I claim; I should now offer some small caveats.  

\emph{Caveat 1: Church had \ZF in mind.} More specifically, Church wanted (what I call) the hereditarily bland sets to obey \ZF. This causes no concern: if we want to follow Church here, we can augment \CUS with the axiom scheme of Replacement, restricted to hereditarily bland sets. This will not affect our synonymy result.; rather, the resulting theory will be synonymous with \ZF itself.

\emph{Caveat 2: Church did not speak about wands.} In a sense this is immaterial; talking about wands is only a helpful heuristic. But there is an important formal point: \CUS's signature has primitives $\blandpred$, $\in$, $\wandpred$ and $\tappred$; Church's own theory has just one primitive, $\in$. Fortunately, this is incidental: in  \S\ref{s:app:CCC:example1.3}, I show that we can reformulate \CUS using just $\in$.

\emph{Caveat 3: I have a restrictive notion of tapping.} As explained: if $b$ is bland, then \CUS proves that $\tapop{0}{\tapop{0}{b}}$ does not exist, but we intuitively want tapping $b$ twice with $0$ to yield $b$. To address this, I simply introduce a more expansive notion of tapping (whose good behaviour is confirmed in Lemma \ref{lem:CUS:compin}):
\begin{define}[\CUS]\label{def:CUS:widetap}
	Let $\widetapop{n}{a} \coloneq \tapop{n}{a}$ if $\tapop{n}{a}$ exists. Let $\widetapop{0}{\tapop{0}{b}} \coloneq b$ if $\tapop{0}{b}$ exists but $\tapop{0}{\tapop{0}{b}}$ does not. In all other cases, let $\widetapop{n}{a}$ be undefined. \defendhere
\end{define}

\emph{Caveat 4: I have a restrictive notion of membership.} I have been glossing $\tapop{0}{a}$ as ``$a$'s complement''. However, I have said nothing, yet, which indicates that $\tapop{0}{a}$'s members should be exactly those $x$ such that $x \notin a$. Worse: whilst $\tapop{0}{\emptyset}$ should be ``the universal set'', \CUS's axioms \ref{WS:notetap} and \ref{WS:notein} entail that $\forall x\ x \notin \tapop{0}{\emptyset}$. Evidently, I need to define a more expansive notion of membership; here it is:
\begin{define}[\CUS]\label{def:CUS:varin}
	When $c$ is bland and $0 < n < \omega$, define:
	\begin{align*}
		x \varin c &\colonequiv x \in c & 
		x \varin\tapop{0}{c} &\colonequiv x \notin c\\
		x \varin \tapop{n}{c} &\colonequiv x \approx_n c &
		x \varin \tapop{0}{\tapop{n}{c}} &\colonequiv x \napprox_n c 
	\end{align*}
	We confirm this covers all cases in Lemma \ref{lem:CUS:kinds}.
	\defendhere
\end{define}

\emph{Caveat 5: Church's definition of $\approx_n$ is not quite my own.} Relatedly: Church defined $\approx_n$ using an ``expansive'' notion of membership, according to which everything is a member of the (non-bland) universal set. I defined $\approx_n$ using the ``restrictive'' notion of membership, $\in$, according to which only bland sets have members. (I am forced to do things this way, because I invoke $\approx_n$ when I define the ``expansive'' notion of membership, $\varin$, in Definition \ref{def:CUS:varin}.) 

Again, the difference is unimportant. Church's definition of $\approx_j$ only shows up in his axiomatic theory in the following axiom:\footnote{\textcite[305]{Church:STUS}.}
	\begin{listclean}
		\item[$\text{L}_j$:] if $a$ is well founded, then $\Setabs{x}{a \approx_j x}$ is a set
	\end{listclean}
Now suppose that $a$ is well-founded and that $a \approx_j x$. Then $x$'s members are bland,\footnote{\textcite[298]{Church:STUS} would have said ``low'' rather than ``bland'', but the point remains.}  and so are the members of members, and so on, digging $j$-levels deep through membership chains. So, as it occurs in $\text{L}_j$, either definition would have the same effect. 

However, there is a more subtle difference. Church's principle $\text{L}_j$ is not a single axiom; instead, he has one axiom for each natural number $j > 0$, with numbers supplied in the metalanguage. This suggests that he has defined $j$-equivalence by \emph{meta}linguistic recursion rather than \emph{object}-linguistic recursion. By contrast, my Definition \ref{def:cardinaln} is fully  object-linguistic, and the success of Definition \ref{def:CUS:varin} requires this fact. 

I conjecture that Church relied on metalinguistic recursion because he had not proved a suitable object-linguistic recursion theorem.\footnote{\textcite[305]{Church:STUS} expresses the desire for a ``possibly more general axiom or axiom schema'' than $\text{L}_j$.} Fortunately, I have such a theorem (cf.\ \S\ref{s:tolt:recursion}), and I am not afraid to use it. 

\emph{Caveat 6: Church had no Beschr\"anktheitsaxiom.}\footnote{Thanks to Thomas Forster for emphasizing the importance of \emph{Beschr\"anktheitsaxiome}.} The axiom \ref{WS:strat} acts as an axiom of restriction, arranging everything into a well-founded hierarchy. Indeed, \CUS proves Lemma \ref{lem:WS:fundamental}, which states that anything non-bland is obtained by starting with something bland and tapping it finitely many times with some wands (i.e.\ some object-linguistic natural numbers). By contrast, Church's theory is open-ended: for all we know, its universe contains objects obtained by using hitherto-unmentioned wands, or beasts which float around outside of any hierarchy.

Deleting \ref{WS:strat} from \CUS would certainly destroy my proof-strategy for obtaining synonymy. (Indeed, I suspect it would actually surrender synonymy.) So: \emph{including \ref{WS:strat} genuinely changes Church's theory}.

That said: this change to the \emph{theory} does not not flag a departure from Church's \emph{approach}. Church established the consistency of his 
\cite*{Church:STUS} theory by outlining what is now called a ``Church--Oswald model'',\footnote{See \textcites[305--7]{Church:STUS}{Oswald:PhD}[\S2]{Forster:CSTUS}{Sheridan:VCST}.} and the denizens of Church--Oswald models are always arranged into a nicely well-founded hierarchy. So the model theory of Church-style approaches has always implicitly invoked a \emph{Beschr\"anktheitsprinzip}, even if his formal theory lacked an explicit \emph{Beschr\"anktheitsaxiom}. Indeed: my framework of wand/set theories really just provides a framework for turning a description of a Church--Oswald model into a formal theory with a \emph{Beschr\"anktheitsaxiom}, which is then (by my Main Theorem) synonymous with a \ZF-like theory.

\subsection{\CUS as a theory of nothing but sets}\label{s:app:CCC:example1.3}
I now want to show that \CUS does what we want it to do. I begin with two trivial results concerning $\tapop{0}$, whose proof I leave to the reader:
\begin{lem}[\ref{CUS:make}, \ref{CUS:0}]\applabel{lem:CUS:00elim}
	If $\tapop{0}{\tapop{0}{a}}$ exists, then both $a = \tapop{0}{\tapop{0}{a}}$.
\end{lem}
\begin{lem}[\ref{CUS:make}, \ref{CUS:n:0b}, \ref{CUS:n:0ub}]\applabel{lem:CUS:ma-0ma} 
	If $m > 0$ and $\tapop{m}{a}$ exists, then $\tapop{0}{\tapop{m}{a}}$ and $\tapop{0}{\tapop{0}{\tapop{m}{a}}}$ exist. 
\end{lem}\noindent
The next result says that, in \CUS, everything is (exclusively) either: bland, the complement of something bland, some cardinality$_n$ of something bland, or the complement of some such cardinality$_n$. 
\begin{lem}[\CUS]\label{lem:CUS:kinds}
	Any $a$ belongs to exactly one of these three kinds:
	\begin{listn-0}
		\item\label{CCC:bland} $a$ is bland;
		\item\label{CCC:one} $a = \tapop{n}{c}$ for some bland $c$ and $n < \omega$;
		\item\label{CCC:compn} $a = \tapop{0}{\tapop{n}{c}}$ for some bland $c$ and $0 < n < \omega$
	\end{listn-0}
	Moreover, if \eqref{CCC:one} or \eqref{CCC:compn} obtains, then $n$ is unique.
\end{lem}
\begin{proof}
	Using Theorem \ref{lem:WS:fundamental}, let $a = \tapop{n_i}\ldots\tapop{n_1}c$, with $c$ bland and $i < \omega$ as small as possible. I now note three facts.

	 \emph{First: there are no adjacent zeros}; i.e.\ there is no $j$ with $n_{j+1} = n_{j} = 0$. 
		For if $\tapop{0}{\tapop{0}{b}}$ exists then $b = \tapop{0}{\tapop{0}{b}}$ by Lemma \ref{lem:CUS:00elim}, but $i$ is minimal.
		
	\emph{Second: if $n_1 > 0$, then either $i=1$, or $i=2$ and $n_2 = 0$.} 
		Suppose $n_2$ exists; since $\tapop{n_1}{c}$ is not bland, $\tapop{n_1}{c} \napprox_m \tapop{n_1}{c}$ for any wand $m$; so $n_2 = 0$ by \ref{CUS:make}. Repeating this reasoning, the only wand $m$ such that $\tapop{m}{\tapop{0}{\tapop{n_1}{c}}}$ exists is $m = 0$; now recall the First fact.
	
	\emph{Third: if $n_{j+1} > 0$, then $j=0$.} Since if $m > 0$ and $\tapop{m}{x}$ exists, $x$ is bland by \ref{CUS:make}. 

	These facts entail that kinds \eqref{CCC:bland}--\eqref{CCC:compn} are exhaustive. It is now easy to show that kinds \eqref{CCC:bland}--\eqref{CCC:compn} are exclusive, and to establish the ``Moreover'' remark. 
\end{proof}\noindent
This result vindicates my expansive notion of ``membership''. Roughly, we want to say that $x$ is in $a$'s cardinal$_n$ iff $x$ and $a$ are $n$-equivalent; and that $x$ is in $a$'s complement iff $x$ is not in $a$; we stipulate this with Definition \ref{def:CUS:varin}, and can see by Lemma \ref{lem:CUS:kinds} that this exhausts all cases. 

It is obvious from Definition \ref{def:CUS:varin} that each $n > 0$ behaves \emph{as} the \textsf{cardinal}$_n$ wand. It is only slightly less obvious that $0$ behaves \emph{as} the \textsf{complement} wand. Specifically, using Definition \ref{def:CUS:widetap} and writing $x \varnotin y$ for $\lnot (x \varin y)$:
\begin{lem}[\CUS]\applabel{lem:CUS:compin} $x \varin a$ iff $x \varnotin \widetapop{0}{a}$.
\end{lem}
\begin{proof}
	By Lemma \ref{lem:CUS:kinds}, we only need to check a few cases. This is easy, using \ref{CUS:make}, \ref{WS:notetap}, \ref{CUS:0:inj}, and \ref{CUS:n:0b}.
\end{proof}\noindent
This yields a generalized version of extensionality:
\begin{lem}[\CUS]\label{lem:CUS:generalext}
	If $\forall x(x \varin a \liff x \varin b)$, then $a = b$.
\end{lem}
\begin{proof}
	By Lemma \ref{lem:CUS:kinds}, we need only check a few cases. This is easy, using Lemma \ref{lem:CUS:compin} and elementary (\CUS-provable) facts about $\approx_n$. 
\end{proof}\noindent
Since the objects of \CUS obey this generalized extensionality principle, we can treat \CUS as a theory of \emph{nothing but sets}. Indeed, as claimed in \emph{Caveat 2} of \S\ref{s:app:CUS:caveats}, we can use Lemmas \ref{lem:CUS:kinds}--\ref{lem:CUS:generalext} to reformulate \CUS using \emph{only} a membership-predicate.

Here is one way to do this. Start by defining an identity-preserving translation with these atomic clauses (and an unrestricted domain-formula):
\begin{align*}
	\blandpred_\bullet(a) &\colonequiv 
	a \notin a \land \exists b \forall x \big(x \in b \liff (x \in a \lor \forall z\ z\notin x)\big)\\
	\wandpred_\bullet(n) &\colonequiv n \in \omega\\
	b \in_\bullet a &\colonequiv b \in a \land \blandpred_\bullet(a)\\
	\taptointprim{n}{a}{c}{\bullet} &\colonequiv 
	\big(n = 0 \land \relforall{d}{\blandpred_\bullet}\exists x(x \in d \liff x \in a) \land \forall x(x \in c \liff x \notin a)\big) \lor {}\\
	&\phantom{{} \colonequiv{}} \big(0 < n < \omega \land a \approx_n^\bullet a \land \forall x(x \in c \liff x \approx_n^\bullet a)\big)
\end{align*}
The translation is well-defined, since $\approx_n$ is defined without using the primitive $\tappred$. Now $\in$ is the only primitive of $\CUS^\bullet$, i.e.\ the $\bullet$-translation of \CUS. 

Trivially, $\bullet$ is an interpretation $\CUS\functionto\CUS^\bullet$. We can also define a translation, $\circ$, in the opposite direction, via $x \in^\circ a \colonequiv x \varin a$. By running these translations back-to-back, it is easy to show that \CUS and $\CUS^\bullet$ are synonymous. Hence \CUS can be reformulated, without gain or loss, as a theory whose only primitive is $\in$, and which obeys Extensionality. 

\appref{app:CUSCUSbullet}
I leave the proof of this claim to keen readers, but here are the main steps. Using Lemmas \ref{lem:CUS:kinds} and  \ref{lem:CUS:generalext}, show that \CUS proves: $\blandpred(a) \liff \blandpred^{\bullet\circ}(a)$; and $x \in a \liff x \in^{\bullet\circ} a$; and $\wandpred(w) = \wandpred^{\bullet\circ}(w)$; and $\tapto{w}{a}{c} \liff \taptoint{w}{a}{c}{\bullet\circ}$. Then use Lemmas \ref{lem:CUS:kinds}$^\bullet$ and \ref{lem:CUS:generalext}$^\bullet$ to show that $\CUS^\bullet$ proves $x \in a \liff x \in^{\circ\bullet} a$. This establishes synonymy, and $\CUS^\bullet$ proves Extensionality by Lemma \ref{lem:CUS:generalext}$^\bullet$.

\subsection{\CUS is a wand/set theory}\label{s:app:CUS:wandset}
We have seen that \CUS does everything we wanted it to. But my last task is to show that \CUS is a wand/set theory, in the sense of Definition \ref{def:WandSet}. To this end, I will close by presenting a wand/set theory, \WSC, for Wand/Set Church theory, and showing that it is just an alternative axiomatization of \CUS. 

\WSC's first few axioms are: \ref{WS:notein}, \ref{WS:notetap}, \ref{WS:funtap}, \ref{WS:ext}, \ref{WS:sep}, \ref{WS:strat}, \ref{WS:weakheight}, and \ref{CUS:omega}. These are shared with \CUS. \WSC's remaining axioms, \ref{WS:make} and \ref{WS:eq}, arise by explicitly defining two formulas, $D$ and $E$ (see \S\ref{s:WS:WandAxioms}). To define $D$, I simply insert a loose bound on the formula in the right-hand-side of \ref{CUS:make}:\footnote{Note: whilst the definition of $x \approx_n y$ conceals some quantifiers, the axioms used in its definition (see Definition \ref{def:cardinaln}) entail that the quantifiers can be restricted to $\Vfrom{(\wevof{x}\cup\wevof{y})^+}$ without loss. So both $D$ and $E$ are provably equivalent (using those axioms) to suitably loosely-bounded formulas.}
\begin{align*}
	\dummydom{n}{a} & \colonequiv (n = 0 \land  (\forall x \in \wevof{a})(\blandpred(x) \lonlyif a \neq \tapop{0}{x})) \text{ or }(0 < n < \omega \land a \approx_n a)
\intertext{To define $E$, I simply list the few cases under which we would want to identify wand-taps, whilst inserting some loose bounds when required:}
	\dummyequiv{m}{a}{n}{b} &\colonequiv m, n \in \omega\text{ and:}\\
	&\hspace{3em} (m = n \land a = b) \text{ or} \eqtag{\textsf{e1}}\label{eq:a=b}\\
	&\hspace{3em} \big(0 < m = n \land a \approx_m b\big) \text{ or}\eqtag{\textsf{e2}}\label{eq:n}\\
	&\hspace{3em}  \big(0 = m  < n\land (\exists d \approx_n b)(a = \tapop{0}\tapop{n}d \land d, \tapop{n}{d} \in \wevof{a})\big) \text{ or}\eqtag{\textsf{e3}}\label{eq:a=0nb}\\
	&\hspace{3em}  \big(0 = n < m\land (\exists d \approx_m a)(b = \tapop{0}\tapop{m}d \land d, \tapop{m}{d} \in \wevof{b})\big)\eqtag{\textsf{e4}}\label{eq:b=0ma}
\end{align*}
\emph{\textbf{Note:} from here onwards, I reserve $D$ and $E$ for these defined predicates, and I reserve $\wanddompred$ and $\equivpred$ for the predicates we obtain from $D$ and $E$ as in \S\ref{s:WS:WandAxioms}.}

It just remains to show that $\WSC = \CUS$. The crucial insight is that, because everything is arranged neatly into wevels, the explicitly stipulated loose bounds have no real effect. Concerning $D$: if $x$ is bland and $a = \tapop{0}{x}$, then we ``know'' that $x$ is found before $a$. Concerning $E$: if there is some $d \approx_n b$ with $a = \tapop{0}{\tapop{n}{d}}$, then there is some $e \approx_n d$ of minimal rank; so $a = \tapop{0}{\tapop{n}{e}}$ and we ``know'' that both $e$ and $\tapop{n}{e}$ are found before $a$. 

It just remains to prove what I have claimed to ``know''. My approach is to build up a reservoir of principles shared by \WSC and \CUS, until the theories come to coincide.
\begin{lem}
	[\WSC/\CUS]
	\label{lem:WSC:inj} \ref{CUS:0:inj} holds
\end{lem}
\begin{proof}
	This is an axiom of \CUS. For \WSC: if $\tapop{0}{a} = \tapop{0}{b}$, then  $\churchequiv{0}{a}{0}{b}$ by \ref{WS:eq}; so $a = b$, either directly via $\equivpred$'s definition or via \eqref{eq:a=b} if $\dummyequiv{0}{a}{0}{b}$.
\end{proof}
\begin{lem}
	[\WSC/\CUS]
	\applabel{lem:WSCCUS:DE} For any wevel $s$:
	\begin{listn-0}
		\item\label{n:DE:E} $E$ restricted to $s$ is an equivalence relation. 
		\item\label{n:DE:D} $D$ restricted to $s$ is preserved under $E$.
	\end{listn-0}
	Hence: $\wanddom{n}{a} \liff \dummydom{n}{a}$ and $\churchequiv{m}{a}{n}{b} \liff \dummyequiv{m}{a}{n}{b}$. 
\end{lem}
\begin{proof}
	\emphref{n:DE:E}
	Fix $m,n,i \in \omega$, and $a,b,c \cfoundat s$. Evidently $\dummyequiv{m}{a}{m}{a}$ by \eqref{eq:a=b}. Now suppose $\dummyequiv{m}{a}{n}{b}$  and $\dummyequiv{m}{a}{i}{c}$; we show that $\dummyequiv{n}{b}{i}{c}$. Most of this holds by elementary reasoning concerning $n$-equivalence; only two cases merit comment. 
	
	\emph{Case when $\dummyequiv{m}{a}{n}{b}$ and $\dummyequiv{m}{a}{i}{c}$, via disjunct \eqref{eq:a=0nb} both times.} 
	So $0 = m < n, i$, there is $d \approx_n b$ such that $a = \tapop{0}{\tapop{n}{d}}$, and there is $e \approx_i c$ such that $a = \tapop{0}{\tapop{i}{e}}$. Now $\tapop{n}{d} = \tapop{i}{e}$ by \ref{CUS:0:inj}, so that $n = i$ and $d \approx_n e$. (In \CUS, this holds by \ref{CUS:n}; in \WSC, note that $\churchequiv{n}{d}{i}{e}$ by \ref{WS:eq}, so $n=i$ and either $d = e$ or $d \approx_n e$, either directly via $\equivpred$'s definition, or via \eqref{eq:a=b}--\eqref{eq:n}.) Hence $b \approx_n c$, so that $\dummyequiv{n}{b}{i}{c}$ via \eqref{eq:n}.
	
	\emph{Case when $\dummyequiv{m}{a}{n}{b}$ and $\dummyequiv{m}{a}{i}{c}$, via disjunct \eqref{eq:b=0ma} both times.} Similar. 
	
	\emphref{n:DE:D}
	Reason through the cases \eqref{eq:a=b}--\eqref{eq:b=0ma}, using \ref{CUS:0:inj} in the last. 
	
	\emph{Hence:} the claim for $\wanddompred$ is trivial; the claim for $\equivpred$ holds by \eqref{n:DE:E}--\eqref{n:DE:D}.
\end{proof}\noindent
From here on, I will invoke Lemma \ref{lem:WSCCUS:DE} without further comment.
\begin{lem}
	[\WSC/\CUS]
	\applabel{lem:cus:CUS:n}
	\ref{CUS:n}--\ref{CUS:n:0ub} hold.
\end{lem}
\begin{proof}
	These are axioms of \CUS; so work in \WSC, and fix positive $m$ and $n$.
	
	\emph{For \ref{CUS:n}.}
	Use \ref{WS:make}, \ref{WS:eq}, and conditions \eqref{eq:a=b}--\eqref{eq:n}.
	
	\emph{For \ref{CUS:n:0b}.} Use \ref{WS:eq} and  \eqref{eq:b=0ma}.
		
	\emph{For \ref{CUS:n:0ub}.} For reductio, suppose $\tapop{m}{a} = \tapop{0}{\tapop{n}{b}}$. Now $\churchequiv{m}{a}{0}{\tapop{n}{b}}$ by \ref{WS:eq}, so by \eqref{eq:b=0ma} there is $d$ with $\tapop{n}{b} = \tapop{0}{\tapop{m}{d}}$ and $\ordrank{\tapop{m}{d}} < \ordrank{\tapop{n}{b}}$. Now $\churchequiv{n}{b}{0}{\tapop{m}{d}}$, so by \eqref{eq:b=0ma} there is $e \approx_n b$ with $\ordrank{\tapop{n}{e}} < \ordrank{\tapop{m}{d}}$. By \ref{CUS:n}, $\tapop{n}{b} = \tapop{n}{e}$. Now $\ordrank{\tapop{n}{b}} < \ordrank{\tapop{n}{b}}$, a contradiction.\end{proof}
\begin{lem}[\WSC/\CUS]
	\label{lem:cus:nicelyrank}
	If $b$ is bland, then $\ordrank{b} < \ordrank{\tapop{0}{b}}$.
\end{lem}
\begin{proof}
	In this case, if $\tapop{0}{b} = \tapop{n}{c}$, then $n = 0$ by \ref{CUS:n:0b}, so $b = c$ by \ref{CUS:0:inj}. Hence $\ordrank{b}+1=\ordrank{\tapop{0}{b}}$ by Lemma \ref{lem:WS:min}. 
\end{proof}
\begin{lem}
	[\WSC/\CUS]
	\label{lem:WSC:maker}
	\ref{CUS:make} and \ref{WS:make} hold
\end{lem}
\begin{proof}
	Use \ref{WS:make}/\ref{CUS:make} and Lemmas \ref{lem:WSCCUS:DE} and \ref{lem:cus:nicelyrank}. 
\end{proof}
\begin{lem}[\WSC/\CUS]\label{lem:WSC:ma}
	If $m > 0$ and $a \approx_m a$, then:
	\begin{listn-0} 
		\item\label{n:maL0ma} $\ordrank{\tapop{m}{a}} < \ordrank{\tapop{0}{\tapop{m}{a}}}$
		\item\label{n:ma00elim} if $\tapop{m}{a} = \tapop{0}{\tapop{0}{c}}$, then $\tapop{m}{a} = c$.
		\item\label{n:aLma} if $(\forall c \approx_m a)\ordrank{a} \leq \ordrank{c}$, then $\ordrank{a} < \ordrank{\tapop{m}{a}}$.
		\item\label{n:ma00intro} $\tapop{m}{a} = \tapop{0}{\tapop{0}{\tapop{m}{a}}}$. 
	\end{listn-0}
\end{lem}
\begin{proof}
	\emphref{n:maL0ma} 
	Now $\tapop{0}{\tapop{m}{a}}$ exists by Lemma \ref{lem:CUS:ma-0ma}. Using Lemma \ref{lem:WS:min}, fix $n$ and $b$ such that $\tapop{n}{b} = \tapop{0}{\tapop{m}{a}}$ and $\ordrank{b} < \ordrank{\tapop{0}{\tapop{m}{a}}}$. Since $a$ is bland, $n = 0$ by \ref{CUS:n:0ub}, so $b = \tapop{m}{a}$ by \ref{CUS:0:inj}.

	\emphref{n:ma00elim} 
	In \CUS, use Lemma \ref{lem:CUS:00elim}. In \WSC: by \ref{WS:eq},  $\churchequiv{m}{a}{0}{\tapop{0}{c}}$, so via \eqref{eq:b=0ma} there is $d\approx_m a$ with $\tapop{0}{c} = \tapop{0}{\tapop{m}{d}}$, hence $c = \tapop{m}{d} = \tapop{m}{a}$ by \ref{CUS:n} and \ref{CUS:0:inj}.
	
	\emphref{n:aLma} 
	Using Lemma \ref{lem:WS:min}, fix $b$ and $n$ such that $\tapop{n}{b} = \tapop{m}{a}$ and $\ordrank{b} < \ordrank{\tapop{m}{a}}$. Suppose for reductio that $b$ is not bland. Now $n = 0$ by \ref{CUS:make}. Since $b$ is not bland, use Lemma \ref{lem:WS:min} again to fix $c$ and $i$ such that $\tapop{i}{c} = b$ and $\ordrank{c} < \ordrank{b}$. Now $\tapop{m}{a} = \tapop{0}{b} = \tapop{0}{\tapop{i}{c}}$, so $i = 0$ by \ref{CUS:n:0ub}, and therefore $\tapop{m}{a}= c$ by \eqref{n:ma00elim}. So $\ordrank{c} < \ordrank{b} < \ordrank{c}$, a contradiction. So $b$ must be bland after all. Now $n\neq 0$ by \ref{CUS:n:0b}, so that $m=n$ and $a \approx_m b$ by \ref{CUS:n}; hence $\ordrank{a} \leq \ordrank{b}$ given our assumptions about $a$.
	
	\emphref{n:ma00intro}
	In \CUS, use Lemmas \ref{lem:CUS:00elim}--\ref{lem:CUS:ma-0ma}. In \WSC: fix $d \approx_m a$ of minimal rank, so that $\tapop{m}{a} = \tapop{m}{d}$ by \ref{CUS:n} and $\ordrank{d} <  \ordrank{\tapop{m}{d}} < \ordrank{\tapop{0}{\tapop{m}{d}}}$ by \eqref{n:aLma} and \eqref{n:maL0ma}. Now $d$ witnesses that $\churchequiv{m}{a}{0}{\tapop{0}{\tapop{m}{d}}}$ via \eqref{eq:b=0ma}. So $\tapop{m}{a} = \tapop{0}{\tapop{0}{\tapop{m}{d}}} = \tapop{0}{\tapop{0}{\tapop{m}{a}}}$ by \ref{WS:eq}.
\end{proof}
\begin{lem}
	[\WSC/\CUS]
	\label{lem:cus:idem}
	\ref{CUS:0} holds
\end{lem}
\begin{proof}
	This is an axiom of \CUS. For \WSC, we argue by induction. Suppose \ref{CUS:0} holds whenever $\ordrank{a} < \alpha$; now let $\ordrank{a} = \alpha$, and suppose $\lnot\blandpred(a)$ and $\relforall{x}{\blandpred}a \neq \tapop{0}{x}$. 
	Using Lemma \ref{lem:WS:min}, fix $n$ and $b$ with $\tapop{n}{b} = a$ and $\ordrank{b} < \ordrank{a}$.
	If $n > 0$, then $a = \tapop{0}{\tapop{0}{a}}$ by Lemma \ref{lem:WSC:ma}\eqref{n:ma00intro}. 
	So suppose instead that $n = 0$, i.e.\ $a = \tapop{0}{b}$.  
	Now $\lnot\blandpred(b)$, by \ref{CUS:0:inj} as $\relforall{x}{\blandpred}\tapop{0}{b} \neq \tapop{0}{x}$. 
	Also $\relforall{x}{\blandpred}b \neq \tapop{0}{x}$, by \ref{CUS:make} since $\tapop{0}{b}$ exists. So $b = \tapop{0}{\tapop{0}{b}}$ by the induction hypothesis, and $a = \tapop{0}{b} = \tapop{0}{\tapop{0}{\tapop{0}{b}}} = \tapop{0}{\tapop{0}{a}}$.
\end{proof}
\begin{lem}[\WSC/\CUS]
	\ref{WS:eq} holds
\end{lem}
\begin{proof}
	This is an axiom of \WSC. For \CUS, fix $\tapop{m}{a}$ and $\tapop{n}{b}$.
	
	\emph{Right to left.} Suppose $\churchequiv{m}{a}{n}{b}$. I claim that $\tapop{m}{a} = \tapop{n}{b}$. If \eqref{eq:a=b}, this is immediate. If \eqref{eq:n}, use \ref{CUS:n}. 
	If \eqref{eq:a=0nb} then $0 = m < n$ and $a = \tapop{0}{\tapop{n}{d}}$ for some $d \approx_n b$; now 
	$\tapop{m}{a} = \tapop{0}{\tapop{0}{\tapop{n}{d}}} = \tapop{n}{d}  = \tapop{n}{b}$ by Lemma \ref{lem:CUS:00elim} and \ref{CUS:n}. If \eqref{eq:b=0ma}, reason similarly. 
	
	\emph{Left to right.} Suppose $\tapop{m}{a} = \tapop{n}{b}$. I claim that $\churchequiv{m}{a}{n}{b}$. 
	If either $m=n=0$, or both $m,n >0$, this is easy; so suppose $m=0 < n$ (the case when $0=m < n$ is similar). 
	Now $b \approx_n b$ by \ref{CUS:make}. Fix $d \approx_n b$ of minimal rank; now $\tapop{n}{b} = \tapop{n}{d}$ by \ref{CUS:n}, and $\ordrank{d} < \ordrank{\tapop{n}{d}} < \ordrank{\tapop{0}{\tapop{n}{d}}}$ by Lemma \ref{lem:WSC:ma}, and $\tapop{0}{\tapop{n}{d}} = \tapop{0}{\tapop{0}{a}} = a$ by Lemma \ref{lem:CUS:00elim}. So $d$ witnesses that $\churchequiv{m}{a}{n}{b}$ via \eqref{eq:a=0nb}. 
\end{proof}\noindent
Assembling everything, $\WSC = \CUS$. So \CUS is a wand/set theory.

\section*{Acknowledgments}
Enormous thanks to Thomas Forster and Randall Holmes: Thomas's enthusiasm for both Conway's liberationism and Church's set theory has proved infectious, and I have had many wonderful conversations with Thomas and Randall about the ideas in this paper. Thanks also to Neil Barton, Hazel Brickhill, Luca Incurvati, \O{}ystein Linnebo, Carlo Nicolai, Chris Scambler, Flash Sheridan, Rob Trueman, and Giorgio Venturi.

\startappendix
\textcolor{myblue}{\begin{center}***\end{center}What came before is the published version of this paper. Some of the proofs in the published version are quite compressed. In case more expansive versions might be helpful, here are some further details.}

\section{Proofs related to \S\ref{s:WS}}\label{proofsforwev}
I start with the results leading to Lemma \ref{lem:WS:acc}.
\begin{lem}[\ref{WS:notein}, \ref{WS:notetap}]\label{lem:WS:cfoundat}
	If $x \cfoundat c \cfoundat a$, then $x \cfoundat a$
\end{lem}
\begin{proof}
	Suppose $x$ is bland, so $x \subseteq c$ by \ref{WS:notetap}. If $c$ is also bland, then $c \subseteq a$ so that $x \subseteq a$. If $c$ is not bland, then $c$ is empty by \ref{WS:notein}; so vacuously $x \subseteq a$. 
	
	Suppose instead $x$ is not bland; so $\tapto{w}{b}{x}$ for some $w$ and some $b \in c$. So $c$ is bland by \ref{WS:notein}, and hence $c \subseteq a$. So $b \in a$ and hence $x \cfoundat a$.
\end{proof}
\begin{lem}[\ref{WS:notein}, \ref{WS:notetap}]\label{lem:WS:cpotpotent}
	If $\cpot{a}$ exists, then  $\cpot{a}$ is wand-potent. 
\end{lem}
\begin{proof}
	Suppose $x \incpot \cpot{a}$. So $x \cfoundat c \in \cpot{a}$ for some $c$. So $c \cfoundat b \in a$ for some $b$; so $x \cfoundat b \in a$ by  Lemma \ref{lem:WS:cfoundat}, i.e.\ $x \in \cpot{a}$.
\end{proof}
\begin{lem}[\CoreWev]\label{lem:WS:levelpottrans}
	Every wevel is wand-potent and wand-transitive. 
\end{lem}
\begin{proof}
	Fix a wevel, $s$; so $s = \cpot{h}$ for some wistory $h$, and $s$ is therefore wand-potent by Lemma \ref{lem:WS:cpotpotent}. To see $s$ is wand-transitive, fix $a \in s = \cpot{h}$; so there is $c$ such that $a \cfoundat c \in h$. 
	Now $c = \cpot{(c \cap h)} \subseteq \cpot{h} = s$, and $c$ is bland, so $c \cfoundat s$. Hence $a \cfoundat s$ by Lemma \ref{lem:WS:cfoundat}.
\end{proof}
\begin{lem}[\CoreWev, scheme]\label{lem:WS:minimal}	
	If something is $\phi$ and every $\phi$ is wand-potent, then  some $\phi$ is $\in$-minimal. Hence: if some wevel is $\phi$, then some $\in$-minimal wevel is $\phi$.
\end{lem}
\begin{proof}[\appproof{lem:WS:minimal}]
	Suppose $\phi(e)$. If $e$ is not bland, it is $\in$-minimal by \ref{WS:notein}. Suppose $e$ is bland; then use \ref{WS:sep} twice to obtain two bland sets: 
	\begin{align*}
		c &\coloneq \Setabs{x \in e}{\forall y(\phi(y) \lonlyif x \in y)} = \Setabs{x}{\forall y(\phi(y) \lonlyif x \in y)}\\
		d &\coloneq \Setabs{x \in c}{x \notin x}
	\end{align*}
	Clearly $d \notin c$; so there is some $a$ such that $\phi(a)$ and $d \notin a$. Since $a$ is wand-potent, also $d \mathrel{\slashed{\incpot}} a$. Now suppose $\phi(b)$; so $d \subseteq c \subseteq b$, so that $d \cfoundat b \in a$, i.e.\ $d \incpot b$; now $b \notin a$. ``Hence'' follows via Lemma \ref{lem:WS:levelpottrans}.
\end{proof}
\begin{lem}[\CoreWev]\label{lem:WS:histin}
	Every member of a wistory is a wevel.
\end{lem}
\begin{proof}
	For reductio, let $h$ be a wistory with some non-wevel in it. Since $c = \cpot(c \cap h)$ for every $c \in h$, every member of $h$ is both bland and wand-potent by Lemma \ref{lem:WS:cpotpotent}. Using Lemma \ref{lem:WS:minimal}, let $a$ be an $\in$-minimal non-wevel in $h$. So $a = \cpot(a \cap h)$, and $a \cap h$ exists by \ref{WS:sep}. So, to obtain our desired contradiction, it suffices to show that $a \cap h$ is itself a wistory, since then $a$ is a wevel after all.
	
	Fix $b \in a \cap h$. Now $b$ is a wevel, by $a$'s minimality. Fix $x \in b$; so $x \cfoundat b \in a$ by Lemma \ref{lem:WS:levelpottrans}, and hence $x \in a$ as $a$ is wand-potent (see above). Generalizing, $b \subseteq a$. Accordingly, $b = \cpot(b \cap h) = \cpot(b \cap (a \cap h))$. Generalizing, $a \cap h$ is a wistory.
\end{proof}
\begin{proof}[\appproof{lem:WS:acc}]
	Throughout, let $k = \Setabs{r \in s}{\wevpred(r)}$. 
	
	\emph{Left to right.} Let $s$ be a wevel; I must show that $s = \cpot{k}$. 
	If $a \in \cpot{k}$, then $a \cfoundat r \in k \subseteq s$ for some $r$, i.e.\ $a \incpot s$, so that $a \in s$ by Lemma \ref{lem:WS:cpotpotent}. 	
	Conversely, suppose $a \in s$. Fix a wistory $h$ such that $s = \cpot{h}$; so $a \incpot h$. By Lemma \ref{lem:WS:histin}, there is a wevel $r$ such that $a \cfoundat r \in h$. Since $r$ is bland, $r \subseteq r \in h$, i.e.\ $r \incpot h$, so that $r \in \cpot{h} = s$. All told, $r \in k$, so $a \cfoundat r \in k$ and hence $a \in \cpot{k}$. 
	
	\emph{Right to left.} Suppose $s = \cpot{k}$. 
	It suffices to show that $k$ is a wistory. Fix $r \in k$ i.e.\ some wevel $r \in s$; it suffices to show that $r = \cpot{(r \cap k)}$. 
	To see that $\cpot{(r\cap k)} \subseteq r$, fix $a \incpot (r \cap k)$; now $a \incpot r$, and $r$ is wand-potent by Lemma \ref{lem:WS:levelpottrans}, so $a \in r$. To see that $r \subseteq \cpot{(r \cap k)}$, fix $a \in r$. By the \emph{left to right} half of this Lemma, there is a wevel $q$ with $a \cfoundat q \in r$. Now $r$ is wand-transitive by Lemma \ref{lem:WS:levelpottrans}, so $q \cfoundat r \in s$ i.e.\ $q \incpot s$; and $s = \cpot{k}$ is wand-potent by Lemma \ref{lem:WS:cpotpotent}, so that $q \in s$. All told, $q \in k$, so $a \cfoundat q \in (r \cap k)$ i.e.\ $a \in \cpot{(r \cap k)}$.
\end{proof}
\begin{proof}[\appproof{thm:WS:wo}]
	For reductio, using Lemma \ref{lem:WS:minimal} twice: let $s$ be the least wevel incomparable with some wevel; let $t$ be the least wevel incomparable with $s$. I will show, absurdly, that $s = t$. 
	
	To show that $s \subseteq t$, fix $a \in s$. By Lemma \ref{lem:WS:acc}, there is a wevel $r$ such that $a \cfoundat r \in s$. The choice of $s$ guarantees that $r$ is comparable with $t$. But if $r = t$ or $t \in r$ then $t \in s$ as $s$ is wand-transitive, contradicting our choice of $t$; so $r \in t$. Now $a \cfoundat r \in t$, so that $a \in t$ as $t$ is wand-potent. Generalizing, $s \subseteq t$.
	
	Exactly similar reasoning shows that $t \subseteq s$. So $s = t$.
\end{proof}%
\begin{proof}[\appproof{lem:WS:levof}]
	\emphref{levofexists} $\wevof{a}$ exists by \ref{WS:strat} and the well-ordering of wevels. Using \ref{WS:sep}, let $b = \Setabs{x \in \wevof{a}}{x \incpot a}$. If $a$ is not bland, then $a$ is empty by \ref{WS:notein}, so $b = \emptyset = \cpot{a}$. Suppose $a$ is bland and let $x \incpot a$; now $a \subseteq \wevof{a}$ so $x \incpot \wevof{a}$ and hence $x \in \wevof{a}$ as $\wevof{a}$ is wand-potent by Lemma \ref{lem:WS:levelpottrans}; so $b = \cpot{a}$. 
	
	\emphref{levofnotin} There is no wevel $t$ with $a \cfoundat t \in \wevof{a}$, so $a \notin \wevof{a}$ by Lemma \ref{lem:WS:acc}. 
	
	\emphref{levofquick} If $r \subseteq s$ then $s \notin r$ by the well-ordering of wevels. Conversely, if $s \notin r$, then either $r \in s$ or $r = s$ by comparability; and $r \subseteq s$ either way, as $s$ is wand-transitive.
	
	\emphref{levofidem}  By \eqref{levofnotin}, $s \notin \wevof{s}$. By  \eqref{levofquick}, $\wevof{s} \notin s$. So $s = \wevof{s}$, by comparability.
	
	\emphref{levofsubs} In this case, $\wevof{a} \notin \wevof{b}$, by definition of $\wevof{b}$, so $\wevof{b} \subseteq \wevof{a}$ by \eqref{levofquick}.
	
	\emphref{levofin} In this case, $b \notin \wevof{b}$ by \eqref{levofnotin}, so $\wevof{a} \nsubseteq \wevof{b}$, and hence $\wevof{b} \in \wevof{a}$ by \eqref{levofquick}.
	
	\emphref{levnotin} If $a \in a$ then $\wevof{a} \in \wevof{a}$ by \eqref{levofin}, contradicting well-ordering.
\end{proof}
\begin{proof}[\appmore{lem:WS:fundamental}]
	I will start by assuming that there is no last wevel. With $\emph{BigTap}$ defined as in the text, stipulate that $\emphspaced{Path}(e, a, n)$ iff: 
	\begin{listbullet}
		\item $\domain{e} = n+1$; and
		\item $\emphspaced{BigTap}(e,n) = a$; and
		\item 
		$\ordrank{(\emphspaced{BigTap}(e, i))} + 1 = 
		\ordrank{(\emphspaced{BigTap}(e, i+1))}$ if $i < n$. 
	\end{listbullet}
	Let $\emphspaced{depth}(a, n)$ be given by: $\ordrank(e(0))$ for any $e$ with $\emphspaced{Path}(e, a, n)$. By Lemma \ref{lem:WS:min}, $\emph{depth}$ is strictly descending where defined, so 	
	fix the largest $m$ where $\emph{depth}(m)$ is defined. Now $e(0)$ is bland whenever $\emph{Path}(e, a, m)$; otherwise by Lemma \ref{lem:WS:min} there would be $d$ with $\emph{Path}(d, a, m+1)$ given by $\tapop{d(1{})}(d(0)) = e(0)$ and $d(i+1) = e(i)$ for $0 < i \leq n+1$. This shows that no version of Choice is needed.
	
	The preceding argument assumes that there is no last wevel. Indeed, if the last wevel is $\wevof{a}$ and $a = \tapop{w}{b}$ with $b$ bland and $\ordrank{b} + 1 = 
	\ordrank{a}$, there will be no $e$ such that $\emph{Path}(e, a, 2)$. 
	However, as in footnote \ref{fn:WS:fundamental}, we can ask whether there are $v_0, \ldots, v_0$ and $e$ and $n$ such that $\emph{Path}(e, v_0, n)$ and $\tapop{v_4}\tapop{v_3}\tapop{v_2}\tapop{v_1}v_0 = a$, thereby ensuring that each pair $\tuple{i, e(i)}$ sits sufficiently beneath $\wevof{a}$.
\end{proof}\noindent 
I will now move on to the proofs related to \S\ref{s:WS:levelsur}.\label{proofsforlevurbase} As noted, we start by mirroring the proofs of Lemmas \ref{lem:WS:cfoundat}--\ref{lem:WS:histin}.
\begin{lem}[$(\forall x \in \urbase)x \notin x$]\label{lem:C0:potpot}
	If $\potur{\urbase}{(a)}$ exists, then $\potur{\urbase}{(a)}$ is potent$_\urbase$ and $\potur{\urbase}{(a)} \notin \urbase \subseteq \potur{\urbase}{(a)}$.
\end{lem}
\begin{proof}
	For potency$_\urbase$: fix bland $x$ and $c \notin \urbase$ such that $x \subseteq c \in \potur{\urbase}{(a)}$, i.e.\ $c \inpot{\urbase} a$. So $c$ is bland and there is $b$ with $c \subseteq b \in a$. Hence $x \subseteq b \in a$, so $x \inpot{\urbase} a$, i.e.\ $x \in \potur{\urbase}{(a)}$. 
	
	To see that $\urbase \subseteq \potur{\urbase}{(a)}$: if $x \in \urbase$ then $x \inpot{\urbase} a$, so $x \in \potur{\urbase}{(a)}$. Finally: if  $\potur{\urbase}{(a)} \in \urbase$, then $\potur{\urbase}{(a)} \in \potur{\urbase}{(a)}$, contradicting our assumption about $\urbase$.
\end{proof}
\begin{lem}[\CoreLev]
	\label{lem:C0:levelpottrans}
	If $s$ is a level$_\urbase$, then $s$ is potent$_\urbase$ and transitive$_\urbase$, and $s \notin \urbase \subseteq s$.
\end{lem}
\begin{proof}
	Let $h$ be a history$_\urbase$ such that $s = \potur{\urbase}{(h)}$. Now Lemma \ref{lem:C0:potpot} gives everything  except transitivity$_\urbase$. For that, fix $a \notin \urbase$ such that $a \in s = \potur{\urbase}{(h)}$. So $a$ is bland and $a \subseteq c \in h$ for some $c$. Now $a \subseteq c = \potur{\urbase}{(c \cap h)} \subseteq \potur{\urbase}{(h)} = s$.
\end{proof}
\begin{lem}[\CoreLev; scheme]\label{lem:C0:minimal}
	If something bland is $\phi$ and every $\phi$ is potent$_\urbase$ and no $\phi$ is in $\urbase$, then there is an $\in$-minimal $\phi$. Hence: if some level$_\urbase$ is $\phi$, some $\in$-minimal level$_\urbase$ is $\phi$.
\end{lem}
\begin{proof}
	Define $c, d$ exactly as in Lemma \ref{lem:WS:minimal}, to obtain $a$ such that $\phi(a)$ and $d \notin a$. Since $a$ is potent$_\urbase$,
	$(\forall y \notin \urbase)(d \subseteq y \lonlyif y \notin a)$. Now suppose $\phi(b)$, so that $b \notin \urbase$ by assumptions about $\phi$, and also $d \subseteq c \subseteq b$. So $b \notin a$.
\end{proof}
\begin{lem}[\CoreLev]\label{lem:C0:histin}
	Every member of a history$_\urbase$ is a level$_\urbase$.
\end{lem}
\begin{proof}
	Almost exactly as Lemma \ref{lem:WS:histin}: by a reductio assumption, there is an $\in$-minimal non-level$_\urbase$ $a \in h$. However, when considering $b \in a \cap h$ and showing that $b \subseteq a$, there are two cases to consider. Fix $x \in b = \potur{\urbase}{(b \cap h)}$. If $x \in \urbase$, then $x \in \potur{\urbase}(a \cap h) = a$; if $x \notin \urbase$, then $x$ is bland, so $x \subseteq b \in a$ since $b$ is transitive by Lemma \ref{lem:C0:levelpottrans}, so that $x \in \inpot{\urbase}(a \cap h) = a$.
	%
	%
	%
\end{proof}\noindent
I can now prove the results mentioned in the main text:
\begin{proof}[\appproof{lem:C0:acc}]
	Throughout, let $k = \Setabs{r \in s}{r\text{ is a level}_\urbase}$.
	
	\emph{Left to right.} Let $s$ be a level$_\urbase$. Suppose $a \in \potur{\urbase}{(k)}$, i.e.\ $a \inpot{\urbase} k \subseteq s$. Making free use of Lemma \ref{lem:C0:levelpottrans}: if $a \in \urbase$ then $a \in s$; alternatively $a$ is bland and $a \subseteq r \in k$, so that $a \in \potur{\urbase}(k) = s$.
	
	Conversely, suppose $a \in s$. Fix a history$_\urbase$ $h$ such that $s = \potur{\urbase}{(h)}$. If $a \in \urbase$, then $a \in \potur{\urbase}{(k)}$. If $a$ is bland and there is $r$ such that $a \subseteq r \in h$, then $r$ is a level$_\urbase$ by Lemma \ref{lem:C0:histin}; also $r \subseteq r \in h$ and $r$ is bland, so that $r \in \potur{\urbase}(h) = s$; hence $r \in k$; now $a \subseteq r \in k$, i.e.\ $a \in \potur{\urbase}{(k)}$. 
	
	\emph{Right to left.} Suppose $s = \potur{\urbase}{(k)}$. It suffices to show that $k$ is a history$_\urbase$. Fix $r \in k$ i.e.\ some level$_\urbase$ $r \in s$; it suffices to show that $r = \potur{\urbase}{(r \cap k)}$. 
	
	Fix $a \in \potur{\urbase}{(r\cap k)}$. If $a \in \urbase$ then $a \in r$ by Lemma \ref{lem:C0:levelpottrans}. Alternatively $a$ is bland and  $a \subseteq q \in r \cap k$ for some $q$, so $q \notin \urbase$ by Lemma \ref{lem:C0:levelpottrans}, and hence $a \in r$ as $r$ is potent$_\urbase$. 
	
	Conversely, fix $a \in r$. By \emph{left to right} of this Lemma, $r = \potur{\urbase}\Setabs{q \in r}{\levpred_\urbase(q)}$. So if  $a \in \urbase$, then $a \in \potur{\urbase}(r \cap k)$. Alternatively, $a$ is bland and there is a level$_\urbase$ $q$ such that $a \subseteq q \in r$. By Lemma \ref{lem:C0:levelpottrans}, $q \notin \urbase$ and $r$ is transitive$_\urbase$, so $q \subseteq r \in k$, i.e.\ $q \in \potur{\urbase}(k)= s$. All told, $q \in k$. Hence $a \subseteq q \in (r \cap k)$ i.e.\ $a \in \potur{\urbase}{(r \cap k)}$.
\end{proof}
\noindent The proof of Lemma \ref{lem:C0:wo} is then almost exactly like the proof of Theorem \ref{thm:WS:wo}. So I close with the claim about hereditarily bland sets. 
\begin{proof}[\appmore{lem:WS:hered}]
	I claimed that that $\breve{s}$ witnesses that $a$ is hereditarily bland. To see that $\breve{s}$ is transitive: if $x \in \breve{s} \subseteq s$ then $x \subseteq s$ by Lemma \ref{lem:WS:levelpottrans}, and $x$'s members are all hereditarily bland by \emph{left to right} of this Lemma, so $x \subseteq \breve{s}$.
\end{proof}

\section{Proofs related to \S\ref{s:int:toc}}
\begin{proof}[\appproof{lem:blt:helow:es}]
	For \ref{lt:ext}$^\toc$, fix hereditarily bland $a$ and $b$ and suppose that $(\forall x(x \in a \liff x \in b))^\toc$; then $\forall x(x \in a \liff x \in b)$ by Lemma \ref{lem:WS:hered}, so $a = b$ by \ref{WS:ext}. Similarly, \ref{lt:scheme}$^\toc$ follows from \ref{WS:sep} by Lemma \ref{lem:WS:hered}.
\end{proof}
\begin{lem}\label{lem:helow:bevlevtrick}
	Recalling Definition \ref{def:LTstuff}, for any wevels $r, s$: 
	\begin{listn-0}
		\item\label{bevlev:helevtrans}$\breve{s}$ is hereditarily bland, potent and transitive
		\item\label{bevlev:sin} 
		$r \in s $ iff $\breve{r} \in \breve{s}$		
		\item\label{bevlev:nochange} 
		$s = \wevof{\breve{s}}$ 
		\item\label{bevlev:bevpot} 
		$\breve{s} = \pot^{\toc}(h)$, where we let $h \coloneq \Setabs{\breve{r} \in \breve{s}}{\wevpred(r)}$
		\item\label{bevlev:levels} 
		$\breve{s}$ is a level$^\toc$
	\end{listn-0}
\end{lem}
\begin{proof}
	I use Lemma \ref{lem:WS:hered} freely and without further comment. 
	
	\emphref{bevlev:helevtrans} Clearly $\breve{s}$ is hereditarily bland; then $\breve{s}$ is potent and transitive by Lemmas \ref{lem:WS:levelpottrans}. 
	
	\emphref{bevlev:sin} \emph{Left to right}. By \eqref{bevlev:helevtrans}. \emph{Right to left}. Note that $\breve{s} \notin \breve{s}$ by Lemma \ref{lem:WS:levof}\eqref{levnotin}. Now, suppose that $\breve{r} \in \breve{s}$. Clearly $r \neq s$. Also $\breve{s} \notin \breve{r}$, since $\breve{r}$ is transitive by \eqref{bevlev:helevtrans}; so $s \notin r$ by \emph{left to right}. So $r \in s$, by well-ordering of wevels.
	
	\emphref{bevlev:nochange} By Lemma \ref{lem:WS:levof}:  $\wevof{\breve{s}} \subseteq \wevof{s} = s$. Conversely, since $\breve{s} \notin \breve{s}$, also $\breve{s} \notin s$, so that $\wevof{\breve{s}} \notin s$ and therefore $s \subseteq \wevof{\breve{s}}$. 
	
	\emphref{bevlev:bevpot} Let $a$ be hereditarily bland. If $a \in \pot^\toc{h}$, then $a \in \breve{s}$ as $\breve{s}$ is potent by \eqref{bevlev:helevtrans}.
	Conversely, if $a \in \breve{s}$, then $a \cfoundat r \in s$ for some wevel $r$ by Lemma \ref{lem:WS:acc}, and $a \subseteq \breve{r} \in \breve{s}$ by \eqref{bevlev:sin}, and so $a \in \pot^\toc(h)$. 
	
	\emphref{bevlev:levels} As $s = \pot^\toc(h)$, as in \eqref{bevlev:bevpot}, it suffices to show that $\histpred^\toc(h)$. If $\breve{r} \in h$, i.e.\ $\breve{r} \in \breve{s}$, 
	then $\breve{r} \cap^\toc h = \Setabs{\breve{q} \in \breve{r}}{\wevpred(q)}$ as $\breve{s}$ is transitive by \eqref{bevlev:helevtrans}; so $\breve{r} = \pot^\toc{(\breve{r} \cap^\toc h)}$ by \eqref{bevlev:bevpot}.
\end{proof}
\begin{proof}[\appproof{lem:blt:levhelow}]
	By Lemma \ref{lem:helow:bevlevtrick}, if $s$ is a wevel then $\breve{s}$ is a level$^\toc$ and $\wevof{\breve{s}} = s$. To complete the proof, it suffices to note that if $p$ and $q$ are distinct levels$^\toc$, then $\wevof{p} \neq \wevof{q}$. To see this: Lemma \ref{lem:blt:helow:es} yields Theorem \ref{thm:WS:wo}$^\toc$, so that $p \in q$ or $q \in p$, so that either $\wevof{p} \in \wevof{q}$ or $\wevof{q} \in \wevof{p}$ by Lemma \ref{lem:WS:levof}\eqref{levofin}.
\end{proof}
\begin{proof}[\appproof{lem:WS:helowinterpret}]
	Fix hereditarily bland $a$ and a wevel $s \supseteq a$; now $a \subseteq^\toc \breve{s}$ by Lemma \ref{lem:WS:hered}, and $\breve{s}$ is a hereditarily bland level$^\toc$ by Lemma \ref{lem:blt:levhelow}, yielding Stratification$^\toc$.  
\end{proof}

\section{Proofs related to \S\ref{s:tolt}}

\begin{proof}[\appproof{lem:LTwand:note}]
	\emph{\ref{WS:notein}$^\star$.} If $x \in_\star a$, then $\blandpred_*(a)$ and $\Delta_\star(a)$, so $\blandpred_\star(a)$. 
	
	\emph{\ref{WS:notetap}$^\star$.} Let $\taptostar{w}{a}{c}$ be witnessed by $\tuple{u,b} \in c$ and $\tuple{w,a,u,b}\in \equivstage{\conchrank{a}}$. Now $w, u$ are wands$_\star$, so $u \neq \emptyset$, i.e.\ $c$ is not bland$_\star$.  
\end{proof}
\begin{proof}[\appproof{lem:LTwand:xor}]
	By construction, either $c = \carrier{s}$ for some $s$ whose members are conches, so that $c$ is bland$_\star$, or $c$ is a $\sigma$-tap for some $\sigma$. If $c$ is a $\sigma$-tap, then there is some $\tuple{w,a} \in c \cap \domstage{\sigma}$ by \eqref{conch:c:dom}, so $w$ is a wand$_\star$ and hence $c$ is not bland$_\star$.
\end{proof}
\begin{proof}[\appproof{lem:LTwand:staristrans}]
	If $\blandpred_\star(a)$, then by definition $\Delta_\star(a)$. 
	
	If $\taptostar{w}{a}{c}$, then $\tuple{w,a} \in \domstage{\conchrank{a}}$ and $\Delta_\star(c)$, (explicitly, if $\star = \tolt$, or because $\conchrank{c} \leq \sigma$ so that $c \in \conchstage{\sigma}$). So now $\conchrank{a} \leq \conchrank{a}$, so that $\Delta_\star(a)$ (immediately if $\star=\tolt$, or because $\conchrank{a} < \sigma$). Also $\wandpred_\star(w)$, so that $\Ufrom{\emptyset}(w)$ and therefore $\Delta_\star(w)$. 
	
	If $b \in_\star a$, then $\Delta_\star(a)$; we must show that $\Delta_\star(b)$. Note that $\blandpred_*(a)$ and $b \in_* a$, so $b \in \uncarrier{a}$. By Lemma \ref{lem:LTwand:xor}, $a$ is not any $\sigma$-tap. If $\star = \tolt$, then $a$ is a conch, hence all its members are conches, so $b$ is a conch. If instead $\star = \paramint{\sigma}$, let $a \in \Ulev{\alpha}{\conchstage{\sigma}}$: if $a \in \conchstage{\sigma}$ then $b \in \uncarrier{a} \subseteq \belowstage{\sigma} \subseteq \conchstage{\sigma} = \Ulev{0}{\conchstage{\sigma}}$; if $\uncarrier{a} \subseteq \Ulev{\beta}{\conchstage{\sigma}}$ for some $\beta < \alpha$ then $b \in \Ulev{\beta}{\conchstage{\sigma}}$.
\end{proof}
\begin{proof}[\appproof{lem:LTwand:transint}]
	\emphref{n:rank} As in Lemma \ref{lem:LTwand:xor}, $\carrier{a}$ is not a $\sigma$-tap. So if $x \in a$ then $x \in \belowstage{\conchrank{\carrier{a}}}$, i.e.\ $\conchrank{x} < \conchrank{\carrier{a}}$. But also if $\beta < \conchrank{\carrier{a}}$ then $\carrier{a} \notin \conchstage{\beta}$, so there is some $x \in a$ with $x \notin \belowstage{\beta}$, i.e.\ $\conchrank{x} \geq \beta$. 
	
	\emphref{n:in} If $b \in_\star \carrier{a}$ then $b \in_* \carrier{a}$ so $b \in a$. Conversely, if $b \in a$ then $b \in_*\carrier{a}$, and we have assumed that $\Delta_\star(a)$.
	
	\emphref{n:sub} First, note that if $c \subseteq a$ then $\Delta_\star(\carrier{c})$. (This is shown as in the proof for ``$\in_\star$'' for Lemma \ref{lem:LTwand:staristrans}.) Now: suppose $c \subseteq a$, and let $b \in_\star \carrier{c}$; so $b \in c \subseteq a$ and now $b \in_\star \carrier{a}$ using \eqref{n:in} twice. Conversely: if $\carrier{c} \subseteq^\star \carrier{a}$ and $\Delta_\star(\carrier{c})$ and $b \in c$, then $b \in a$ using \eqref{n:in} twice.
\end{proof}
\begin{proof}[\appproof{lem:LTwand:extsepstar}]
	\emph{For \ref{WS:ext}$^\star$.} 
	Let $\carrier{a}$ and $\carrier{b}$ satisfy $\Delta_\star$, with $\relforall{x}{\Delta_\star}(x \in_\star \carrier{a} \liff x \in_\star \carrier{b})$. By Lemma \ref{lem:LTwand:transint}, $\forall x(x \in a \liff x \in b)$; so $a = b$ and hence $\carrier{a} = \carrier{b}$. 
	
	\emph{For \ref{WS:sep}$^\star$.}  
	Fix any formula $\phi$. Let $\carrier{a}$ satisfy $\Delta_\star$. Let $c = \Setabs{x \in a}{\phi}$ (choosing `$c$' to avoid clash with $\phi$); note that $\Delta_\star(\carrier{c})$ by Lemma \ref{lem:LTwand:transint}. Now $x \in_\star \carrier{c}$ iff $x \in c$ iff $x \in a \land \phi$ iff $x \in_\star \carrier{a} \land \phi$. 
\end{proof}
\begin{proof}[\appmore{prop:LTwand:iso}]
	\emphref{Theta:Ulev} Let $\heblandpred^\star(b)$. Suppose for induction that if $\heblandpred^\star(x)$ and $\conchrank{x} < \conchrank{b}$ then $x \in \Ulev{\conchrank{x}+1}{\emptyset}$. Since $\conchrank{b} = \lsub_{x \in \uncarrier{b}}\conchrank{x}$ by Lemma \ref{lem:LTwand:transint}, we have $\uncarrier{b} \subseteq \Ulev{\conchrank{b}}{\emptyset}$ so $b \in \Ulev{\conchrank{b}+1}{\emptyset}$. 
	
	\emphref{Theta:iso:in} $\Theta a \in^{\toc\star} \Theta b$ iff 
	$\Theta a \in_\star \Theta b \land \heblandpred^\star(\Theta b)$ iff 
	$\Theta a \in_\star \Theta b$ iff 
	$\Theta a \in_* \Theta b$ iff $a \in b$.
	
	\emphref{Theta:iso:wandpredlt} $\wandpredlt^{\toc\star}(\Theta w)$ iff 
	$\relexists{v}{\wandpredlt}\Theta w = \Theta v \land \Ufrom{\emptyset}(\Theta w)$ iff $\wandpredlt(w) \land \Ufrom{\emptyset}(\Theta w)$ iff $\wandpredlt(w)$. 
\end{proof}
\begin{proof}[\appproof{lem:LTwand:nonbland}] 
	%
	\emphref{n:tap:alpha+1} 
	Evidently $c \in \conchstage{\alpha+1}$; and $c \notin \conchstage{\beta}$ for any $\beta \leq \alpha$ by Lemma \ref{lem:LTwand:xor}. When $\star = \paramint{\sigma}$, Definition \ref{def:internalhierarchy} entails that $c \in \conchstage{\sigma}$ so that $\conchrank{c} \leq \sigma$. 
	
	\emphref{n:tap:ranklower} 
	Let $u$ and $b$ witness that $\taptostar{w}{a}{c}$; so  $\tuple{w,a,u,b}\in\equivstage{\conchrank{a}}$ and $\tuple{u,b} \in c$. Then $\conchrank{b} \leq \conchrank{a}$, and $\alpha = \conchrank{b}$ by \eqref{conch:c:rank}, and $\conchrank{c} = \conchrank{b} + 1$ by \eqref{n:tap:alpha+1}.
	
	\emphref{n:tap:memtap} Fix $\tuple{w,a} \in c$. So $\conchrank{a} = \alpha$ by \eqref{conch:c:rank} and $\tuple{w,a} \in \domstage{\alpha}$ by \eqref{conch:c:dom}. So $\wandpred_\star(w)$ so that $\tuple{w,a,w,a} \in \equivstage{\alpha}$. Now $w$ and $a$ themselves witness that $\taptolt{w}{a}{c}$. When $\star = \paramint{\sigma}$, note that $\conchrank{c} = \alpha+1$ by \eqref{n:tap:alpha+1}, and that $c \in \conchstage{\sigma}$ since $\Delta_{\paramint{\sigma}}(c)$, so that $\conchrank{a} < \conchrank{c} \leq \sigma$. 
\end{proof}
\begin{proof}[\appmore{lem:LTwand:helpercum}]
	\emph{When $c$ is bland$_\star$.} 
	Now $c \cfoundat^\star \carrier{\belowstage{\alpha}}$ iff $c \subseteq^\star \carrier{\belowstage{\alpha}}$, 
	iff $\uncarrier{c} \subseteq \belowstage{\alpha}$
	by Lemma \ref{lem:LTwand:transint}, iff  $\conchrank{c} \leq \alpha$. 
	
	\emph{If $c$ is some $\sigma$-tap.} 
	Now $\tuple{w,a} \in c$ for some $a$, and $\taptostar{w}{a}{c}$ by Lemma \ref{lem:LTwand:nonbland}\eqref{n:tap:memtap}. Since $\conchrank{a} < \conchrank{c}$, we have $a \in \belowstage{\conchrank{c}}$, so that $a \in_\star \carrier{\belowstage{\conchrank{c}}}$ by Lemma \ref{lem:LTwand:transint}. Also, $w$ is a wand$_\star$ by \eqref{conch:c:dom} and so $\Delta_\star(w)$. Assembling all this, $c \cfoundat^\star \carrier{\belowstage{\conchrank{c}}}$. To see this is minimal, suppose $\taptostar{u}{b}{c}$. Then $\conchrank{c} \leq \conchrank{b} + 1$ by Lemma \ref{lem:LTwand:nonbland}\eqref{n:tap:ranklower}. So if $\beta < \conchrank{c}$, then $b \notin \belowstage{\beta}$ i.e.\ $b \notin^\star \carrier{\belowstage{\beta}}$; hence $c \ncfoundat^\star \carrier{\belowstage{\beta}}$. 
\end{proof}
\begin{proof}[\appmore{lem:LTwand:levelnu}]
	Let $c$ satisfy $\Delta_\star$. 
	Suppose for induction that $\carrier{\Ulev{\beta}{\urbase}}$ is the $\beta^\text{th}$ level$^\star_{\carrier{\urbase}}$ for all $\beta < \alpha$. By Lemmas \ref{lem:LTwand:transint} and \ref{lem:LTwand:helpercum}, and the definition of $\Ulev{\alpha}{\urbase}$, we have: 
\begin{align*}
	c \in_\star \carrier{\Ulev{\alpha}{\urbase}} \text{ iff }& 
	\big(\blandpred_\star(c) \land \relexists{r}{\levpred^\star_{\carrier{\urbase}}}(c \subseteq^\star r \in_\star \carrier{\Ulev{\alpha}{\urbase}} )\big)	\text{ or }c \in_\star \carrier{\urbase}
\end{align*}
So $\Ulev{\alpha}{\urbase}$ is the $\alpha^\text{th}$ level$^\star_{\carrier{\urbase}}$ by Lemma \ref{lem:C0:acc}$^\star$. Now use induction on level$^\star_{\carrier{\urbase}}$, i.e.\ Lemma \ref{lem:C0:wo}$^\star$.
\end{proof}
\begin{proof}[\appmore{lem:LTwand:lastbits}]	
I start by spelling out the very last part of the proof in the main text, which establishes that $c = d$. Let $\conchrank{a_1} = \conchrank{a_2} = \alpha$. Fix $\tuple{u,b} \in c$; now $\tuple{w_1, a_1, u, b} \in \equivstage{\alpha}$ by \eqref{conch:c:equiv}, and $\conchrank{b} = \alpha$ by \eqref{conch:c:rank}, so $\tuple{w_2, a_2, u, b} \in \equivstage{\alpha}$ by \eqref{conch:c:equiv}, and now $\tuple{u,b} \in d$ by \eqref{conch:c:equiv}. So $c \subseteq d$. Similarly $d \subseteq c$.

I now consider the other cases, omitted from the main text.

\emph{\ref{WS:make}$^\tolt$.} \emph{Right to left.} Suppose $\taptolt{w}{a}{c}$. So there is $\tuple{u,b} \in c$ with  $\churchequivtolt{w}{a}{u}{b}$; now $\wanddomtolt{u}{b}$ by  \eqref{conch:c:dom}, so $\wanddomtolt{w}{a}$ by \domref{dom:default}$^\tolt$. \emph{Left to right.} Suppose $\wanddomtolt{w}{a}$. Let $\sigma$ be minimal such that $\exists u \exists b (\conchrank{b} = \sigma \land \churchequivtolt{w}{a}{u}{b})$; note that $\sigma$ exists by \equivref{equiv:default}$^\tolt$. Define $c = \Setabs{\tuple{u,b} \in \domstage{\sigma}}{\churchequivtolt{w}{a}{u}{b}}$. 
It is easy to verify using \domref{dom:default}$^\tolt$ and \domref{dom:wand}$^\tolt$ that $c$ is a $\sigma$-tap, so that $\taptolt{w}{a}{c}$. 

\emph{\ref{WS:eq}$^\tolt$.} Suppose $\taptolt{w}{a}{c}$. Fix $\tuple{w',a'} \in c$ such that $\churchequivtolt{w}{a}{w'}{a'}$ and $c$ is an $\conchrank{a'}$-tap. 
\emph{Left to right.} Suppose $\taptolt{u}{b}{c}$. Fix $\tuple{u',b'} \in c$ with $\churchequivtolt{u}{b}{u}{b'}$. Now $\churchequivtolt{w'}{a'}{u'}{b'}$ by \eqref{conch:c:equiv}, so $\churchequivtolt{w}{a}{u}{b}$. \emph{Right to left.} Suppose $\churchequivtolt{w}{a}{u}{b}$. By \ref{WS:make}$^\tolt$, $\wanddomtolt{w}{a}$, so that $\wanddomtolt{u}{b}$ by \domref{dom:default}$^\tolt$. Also $\churchequivtolt{w'}{a'}{u}{b}$ so that $\conchrank{a'} \leq \conchrank{b}$ by \eqref{conch:c:rank}. So $w',a'$ witness that $\taptolt{u}{b}{c}$.
\end{proof}

%

%

\section{Proofs related to \S\ref{s:app:CUS}}\label{s:app:CUSstuff}
\begin{proof}[\appproof{lem:CUS:00elim}]
Since $\tapop{0}{\tapop{0}{a}}$ exists, $\relforall{b}{\blandpred}\tapop{0}{a} \neq \tapop{0}{b}$ by \ref{CUS:make}, so $a$ is not bland. Also $\tapop{0}{a}$ exists, so again $\relforall{b}{\blandpred}a \neq \tapop{0}{b}$. Now use \ref{CUS:0}.
\end{proof}
\begin{proof}[\appproof{lem:CUS:ma-0ma}]
Here, $\tapop{0}{\tapop{m}{a}}$ exists by \ref{CUS:n:0b} and \ref{CUS:make}. So $\tapop{0}{\tapop{0}{\tapop{m}{a}}}$ exists by \ref{CUS:n:0ub} and \ref{CUS:make}. 
\end{proof}
\begin{proof}[\appmore{lem:CUS:compin}]
\emph{Case (\ref{CCC:bland}), when $a = c$ for bland $c$.} Then $\tapop{0}{a}$ exists by \ref{CUS:make} and \ref{WS:notetap}, so $\widetapop{0}{a} = \tapop{0}{a}$; now inspect the definition

\emph{Case (\ref{CCC:one}.0), when $a = \tapop{0}{c}$ for bland $c$.} Now $\widetapop{0}{a} = c$ by \ref{CUS:make}, so we are asking whether $x \varin \tapop{0}{c}$ iff $x \varnotin c$; this reduces to case (\ref{CCC:bland}).

\emph{Case (\ref{CCC:one}.$n$), when $a = \tapop{n}{c}$ for bland $c$ and $n > 0$.} Now $\widetapop{0}{a} = \tapop{0}{\tapop{n}{c}}$ by Lemma \ref{lem:CUS:ma-0ma}; now inspect the definition. 

\emph{Case (\ref{CCC:compn}), when $a = \tapop{0}{\tapop{n}}{c}$ for bland $c$ and $n > 0$.} 
Now $\tapop{0}{\tapop{0}{\tapop{n}{c}}} = \tapop{n}{c}$ by Lemmas \ref{lem:CUS:00elim}--\ref{lem:CUS:ma-0ma}. 
So $\widetapop{0}{a} = \tapop{n}{c}$, which reduces to case (\ref{CCC:one}.$n$).
\end{proof}\noindent
I will now prove that \CUS and $\CUS^\bullet$ are synonymous. I start with a simple result:\label{app:CUSCUSbullet}
\begin{prop}[\CUS]\label{prop:CUS:synonymy}
$\bullet\circ$ is equivalent to identity, in that: $\blandpred(a) \liff  \blandpred^{\bullet\circ}(a)$; and $x \in a \liff x \in^{\bullet\circ}a$; and $\wandpred(n)\liff  \wandpred^{\bullet\circ}(n)$; and $\tapto{n}{a}{c}\liff \taptoint{n}{a}{c}{\bullet\circ}$.
\end{prop}
\begin{proof}
By Lemma \ref{lem:CUS:generalext}, $\emptyset$ is the object such that $\forall z\ z \varnotin \emptyset$. 

\emph{Concerning $\blandpred$.} Unpacking the definitions, we have:
\begin{align*}
	\blandpred^{\bullet\circ}(a)&\text{ iff }
	a \varnotin a \land 
	\exists b\forall x(x \varin b \liff (x \varin a \lor z = \emptyset))
\end{align*}
If $a$ is bland, this clearly holds, with bland $b = a \cup \{\emptyset\}$ as witness. When $a$ is not bland, by Lemma \ref{lem:CUS:kinds} there is some bland $c$ such that one of these three situations holds, and in each case, $\lnot\blandpred^{\bullet\circ}(a)$.
\begin{listbullet}
	\item \emph{$a = \tapop{0}{c}$}: now $\ordrank{c} < \ordrank{a}$ by Lemma \ref{lem:cus:nicelyrank}, so $a \notin c$ and hence $a \varin a$.
	\item \emph{$a = \tapop{n}{c}$ with $n > 0$}: by Lemma \ref{lem:CUS:kinds}, there is no $b$ whose members$_{\varin}$ are exactly $\emptyset$ and everything $n$-equinumerous with $c$. 
	\item \emph{$a = \tapop{0}{\tapop{n}{c}}$ with $n > 0$}: now $a \napprox_n c$, so $a \varin a$. 
\end{listbullet}

\emph{Concerning $\in$.} Unpacking the definitions, and using the above, we get:
\begin{align*}
	x \in^{\bullet\circ} a 
	&\text{ iff }x \varin a \land \blandpred(a)
\end{align*}
Now if $x \varin a$ and $\blandpred(a)$, then $x \in a$ by definition. And if $x \in a$ then $\blandpred(a)$ by \ref{WS:notein}, so $x \varin a$ by definition.

\emph{Concerning $\wandpred$.} By \CUS's axiom \ref{CUS:omega}, we have: $\wandpred(x)$ iff $x \in \omega$. Now $x \in \omega$ iff $(x \in \omega)^{\bullet\circ}$, since these notions are defined only using $\blandpred$ and $\in$. 

\emph{Concerning $\tappred$.} Note that $x \approx_n y$ is defined using only $\blandpred$ and $\in$; so $x \approx_n y$ iff $x \approx^{\bullet\circ}_n y$. Similarly, notions like $(0 < n < \omega)^{\bullet\circ}$ translate verbatim. Using all this, unpacking the definitions gives us:
\begin{align*}
	\taptoint{n}{a}{c}{\bullet\circ} 
	&\text{ iff }
	\big(n = 0 \land 
	\relforall{b}{\blandpred}\exists x(x \varin b \liff x \varin a) \land 		
	\forall x(x \varin c \liff x \varnotin a)\big) \lor {}\\
	&\phantom{\text{ iff }}\big(0 < n < \omega \land a \approx_n a \land \forall x(x \varin c \liff x \approx_n a)\big)
	\intertext{Then applying \ref{WS:ext} and the definition of $\varin$:}
	\taptoint{n}{a}{c}{\bullet\circ}&\text{ iff }\big(n = 0 \land \relforall{b}{\blandpred}a \neq \tapop{0}{b} \land \widetapop{0}{a} = c\big) \lor (0 < n < \omega \land \tapop{n}{a} = c)
\end{align*}
which is equivalent to $\tapto{n}{a}{c}$ by \ref{CUS:make}.
\end{proof}\noindent
So $\circ : \CUS^\bullet \functionto \CUS$ is an interpretation, and indeed half of a synonymy. For the other half, we must prove some basic results about our ``wands'':
\begin{lem}[$\CUS^\bullet$]\label{lem:CUSbullet:help}\phantom{.}
\begin{listn-0}
	\item\label{bullet:blandin} If $\blandpred_\bullet(a)$, then $x \in a$ iff $x \in^{\circ\bullet} a$
	\item\label{bullet:empty} $\emptyset$ is (the empty set)$^{\circ\bullet}$
	\item\label{bullet:hebland} If $\heblandpred^\bullet(a)$ and $\forall x(x \in a \liff x \in b)$, then $a = b$
	\item\label{bullet:ord}
	$\omega$ is (the set of finite ordinals)$^{\circ\bullet}$
\end{listn-0}
\end{lem}
\begin{proof}
\emphref{bullet:blandin}
Note that $x \in^{\circ\bullet} a$ is a disjunction of the form: $$((\blandpred_\bullet(a) \land x \in_\bullet a) \lor (\lnot \blandpred_\bullet(a) \land \ldots))$$
Given $\blandpred_\bullet(a)$, this immediately reduces to $x \in_\bullet a$, i.e.\ $x \in a$. 

\emphref{bullet:empty} Using Lemma \ref{lem:CUS:generalext}$^\bullet$, there is a unique $z$ such that $\forall x\ x \notin^{\circ\bullet} z$, and this $z$ is $\blandpred_\bullet$. Now use \eqref{bullet:blandin}.

\emphref{bullet:hebland} Fix such $a$ and $b$. I claim that $b$ is also bland$_\bullet$. Since $a$ is bland$_\bullet$ by assumption, there is $c$ whose only members are $a$'s members (i.e.\ $b$'s members) and anything empty. For reductio, suppose $b \in b$; then $b \in a$, so $b \in_\bullet a$, so that $\heblandpred_\bullet(b)$ by Lemma \ref{lem:WS:hered}$^\bullet$, so $b$ is bland$_\bullet$ and hence $b \notin b$ a contradiction. So $b\notin b$ and hence $b$ is bland$_\bullet$. Now by \eqref{bullet:blandin} and our assumption 	$\forall x(x \in^{\circ\bullet} a \liff x \in^{\circ\bullet} b)$, so that $a = b$ by Lemma \ref{lem:CUS:generalext}$^\bullet$.	

\emphref{bullet:ord} Let $o$ be (the set of finite ordinals)$^{\circ\bullet}$. Using \eqref{bullet:blandin}, $o$ witnesses its own hereditarily blandness$^\bullet$. 
%
So, using \eqref{bullet:empty}--\eqref{bullet:hebland}, $o$'s definition reduces to the definition of the set of finite ordinals, i.e.\ $\omega$.
%
%
%
%
%
%
\end{proof}
\begin{prop}[$\CUS^\bullet$]\label{prop:CUSbullet:synonymy}
$\circ\bullet$ is equivalent to identity, in that: $x \in a$ iff $x \in^{\circ\bullet} a$
\end{prop}
\begin{proof}
Using Lemma \ref{lem:CUSbullet:help} and Lemma \ref{lem:CUS:kinds}$^\bullet$, to say that $x \in^{\circ\bullet} a$ is to say that one of these four disjuncts obtains:
\begin{listn-0}
	\item[(\ref{CCC:bland})] $\blandpred_\bullet(a)$ and $x \in a$; or
	
	
	\item[(\ref{CCC:one}i)] $\forall z(z \in a \liff z \notin c)$ and $x \notin c$, for some $\blandpred_\bullet(c)$, so that $x \in a$; or
	
	
	\item[(\ref{CCC:one}ii)] $\forall z(z \in a \liff z \approx_n^\bullet c)$ and $x \approx_n^\bullet c \approx_n^\bullet c$, for some $\blandpred_\bullet(c)$ and some $ 0 \in n \in \omega$, so that $x \in a$; or	
	
	
	\item[(\ref{CCC:compn})] $\forall z(z \in a \liff z \napprox_n^\bullet c)$ and $x \napprox_n^\bullet c \approx_n^\bullet c$, for some $\blandpred_\bullet(c)$ and some $ 0 \in n \in \omega$, so that $x \in a$. \qedhere
\end{listn-0}
\end{proof}\noindent 
Together, Propositions \ref{prop:CUS:synonymy} and \ref{prop:CUSbullet:synonymy} show that \CUS and $\CUS^\bullet$ are synonymous.
\begin{proof}[\appmore{lem:WSCCUS:DE}]
\emph{Case via disjunct \eqref{eq:b=0ma} both times.} 
In this case: $0 = n = i < m$, there is $d \approx_m a$ such that $b = \tapop{0}{\tapop{m}{d}}$; and there is $e \approx_m a$ such that $c = \tapop{0}{\tapop{m}{e}}$. Now $d \approx_m e$, so that $\tapop{m}{d} = \tapop{m}{e}$. (In \CUS, this holds by \ref{CUS:n}; in \WSC, note that $\dummyequiv{m}{d}{m}{e}$ by \eqref{eq:n}, and hence $\tapop{m}{d} = \tapop{m}{e}$ by \ref{WS:eq}.) So $b = c$, so that $\dummyequiv{n}{b}{i}{c}$ via \eqref{eq:a=b}.	

\emphref{n:DE:D} Fix wands $m,n$ and $a, b \cfoundat s$, and suppose that $\dummydom{m}{a}$ and $\dummyequiv{m}{a}{n}{b}$. 
If \eqref{eq:a=b}, 	\eqref{eq:n} or \eqref{eq:a=0nb} then trivially $\dummydom{n}{b}$. 
If \eqref{eq:b=0ma}, then $0 = n < m$ and there is $d \approx_m a$ with  
$b = \tapop{0}{\tapop{m}{d}}$; now if $b = \tapop{0}{\tapop{m}{d}} = \tapop{0}{c}$ for some $c$, then $c = \tapop{m}{d}$ by \ref{CUS:0:inj}, so that $c$ is not bland; hence $\relforall{x}{\blandpred}(b \neq \tapop{0}{x})$, and therefore $\dummydom{n}{b}$. 
\end{proof}
\begin{proof}[\appmore{lem:cus:CUS:n}]
\emph{For \ref{CUS:n}. Left to right.} Suppose $\tapop{m}{a} = \tapop{n}{b}$. Now $\wanddom{m}{a}$ and $\wanddom{n}{b}$ by \ref{WS:make}, so that $a \approx_m a$ and $b \approx_n b$. Also $\churchequiv{m}{a}{n}{b}$ by \ref{WS:eq}, so either \eqref{eq:a=b} or \eqref{eq:n} holds, and hence $m = n$ and $a\approx_m b$. \emph{Right to left.} Suppose $m=n$ and $a \approx_m b$; then $\churchequiv{m}{a}{n}{b}$ by \eqref{eq:n}, and $\tapop{m}{a} = \tapop{n}{b}$ by \ref{WS:make} and \ref{WS:eq}.

\emph{For \ref{CUS:n:0b}.} Suppose $\tapop{m}{a} = \tapop{0}{b}$. Then $\churchequiv{m}{a}{0}{b}$ by \ref{WS:eq}, so $b=\tapop{0}{\tapop{m}{d}}$ for some $d$ by \eqref{eq:b=0ma}, and hence $b$ is not bland. 
\end{proof}

\stopappendix

\printbibliography

 \end{document}